\documentclass{amsart}
\usepackage{graphicx,verbatim}
\usepackage{hyperref}  
\usepackage{tikz}
  \usepackage{tkz-graph}
\usetikzlibrary{matrix,arrows,decorations.markings}
\tikzset{->-/.style={decoration={markings,mark=at position #1 with {\color{black}\arrow{>}}},postaction={decorate,very thick}}}
\tikzstyle{vertex}=[circle, draw, inner sep=0pt, minimum size=6pt]
\newcommand{\vertex}{\node[vertex]}
\usepackage{rotating}
\usepackage{changepage}

\pgfrealjobname{ellipticgraphs}
\usepackage{tikz}\usetikzlibrary{shapes,arrows}
\usetikzlibrary{matrix}
\usetikzlibrary{decorations.pathreplacing}

\tikzstyle{overbrace text style}=[font=\tiny, above, pos=.5, yshift=3mm]
\tikzstyle{overbrace style}=[decorate,decoration={brace,raise=2mm,amplitude=3pt}]
\tikzstyle{underbrace style}=[decorate,decoration={brace,raise=2mm,amplitude=3pt,mirror},color=black]
\tikzstyle{underbrace text style}=[font=\normalsize, below, pos=.5, yshift=-3mm]
\newtheorem{theorem}{Theorem}[section]
\newtheorem{proposition}[theorem]{Proposition}
\newtheorem{lemma}[theorem]{Lemma}
\newtheorem{corollary}[theorem]{Corollary}
\theoremstyle{definition}
\newtheorem{definition}[theorem]{Definition}
\newtheorem{example}[theorem]{Example}

\newtheorem{remark}[theorem]{Remark}
\newcommand{\Q}{\mathbb{Q}}
\newcommand{\Z}{\mathbb{Z}}

\newcommand{\rk}{\mathrm{rk}}

\newcommand{\C}{\mathbb{C}}

\newcommand{\N}{\mathbb{N}}

\newcommand{\E}{\mathcal{E}}
\newcommand{\F}{\mathcal{F}}
\newcommand{\Line}{\mathcal{L}}
\newcommand{\G}{\mathcal{G}}
\newcommand{\Fq}{\mathbb{F}_q}
\numberwithin{equation}{section}
\usepackage[matrix,arrow]{xy}

\newcommand{\ndiv}{\hspace{-4pt}\not|\hspace{2pt}}

\DeclareMathOperator{\Spec}{\mathrm{Spec}}

\allowdisplaybreaks

\begin{document}

\title{Hall algebras and graphs of Hecke operators for elliptic curves}


\author{Roberto  Alvarenga}
\address{Universidade de S\~{a}o Paulo - ICMC, S\~{a}o Carlos, Brazil}
\curraddr{}
\email{alvarenga@icmc.usp.br}





\maketitle

\noindent
\textbf{Abstract.} The graph of a Hecke operator encodes all information about the action of this operator on automorphic forms over a global function field. These graphs were introduced by Lorscheid in \cite{oliver-graphs} for $\text{PGL}_{2}$ and generalized to $\text{GL}_{n}$ in \cite{roberto-graphs}.
After reviewing some general properties, we explain the connection to the Hall algebra of the function field. In the case of an elliptic function field, we can use structure results of Burban-Schiffmann \cite{olivier-elliptic1} and Fratila \cite{dragos} to develop an algorithm which explicitly calculates these graphs. We apply this algorithm to determine some structure constants and provide explicitly the rank two case in the last section.


\tableofcontents

\section*{Introduction}

This work is concerned with the graphs of Hecke operators. These graphs are defined from the action of Hecke operators on automorphic forms over a global function field. Motivated by questions of Zagier (\cite{zagier}) about unramified toroidal automorphic forms (we refer to \cite{oliver-thesis} for full details), Lorscheid develops in \cite{oliver-graphs} a theory of graphs of Hecke operators for $\mathrm{PGL}_2$ over a global function field. This theory plays an important role in the proofs of his main theorems. In \cite{oliver-elliptic}, Lorscheid analyses this theory for elliptic function fields and answers some of Zagier's questions. In \cite{roberto-graphs}, we extend the definition of these graphs from $\mathrm{PGL}_2$ to $\mathrm{GL}_n$, generalize some of Lorscheid's results and describe how to obtain these graphs for a rational function field. In this paper, we aim to describe these graphs when the global function field is elliptic. Using the theory of Hall algebras, we exhibit an algorithm to calculate these graphs. 

Before describing the reformulation in terms of coherent sheaves, let us review the original definition of the graph of a Hecke operator.

\subsection*{Graph of a Hecke operator}
Let $F$ be the function field of a smooth projective and geometrically  irreducible curve $X$ over $\mathbb{F}_q$, $\mathbb{A}$ its ring of the adelic integers and $K = \mathrm{GL}_n(\mathcal{O}_{\mathbb{A}}),$ where $\mathcal{O}_{\mathbb{A}}$ is the set of the adelic integers.  For any right $K$-invariant Hecke operator $\Phi$, there are unique $m_1, \ldots,m_r \in \C^{*}$ and pairwise distinct $[g_1], \ldots, [g_r] \in \mathrm{GL}_n(F) \setminus \mathrm{GL}_n(\mathbb{A}) / K $ such that for all automorphic forms $f$
$$\Phi(f)(g) = \sum_{i=1}^{r} m_i f(g_i)$$
(see  \cite{roberto-graphs} Proposition 1.6). For $[g],[g_1], \ldots,[g_r] \in \mathrm{GL}_n(F) \setminus \mathrm{GL}_n(\mathbb{A}) / K$ as before, we write  $\mathcal{V}_{\Phi,K}([g]) := \{([g],[g_i],m_i)\}_{i=1, \ldots, r}.$
From that, we define the graph $\mathcal{G}_{\Phi,K}$ of $\Phi$ relative to $K$, whose vertices are
$$\mathrm{Vert} \;\mathcal{G}_{\Phi,K} = \mathrm{GL}_n(F) \setminus \mathrm{GL}_n(\mathbb{A}) / K$$
and whose oriented weighted edges 
$$\mathrm{Edge}\; \mathcal{G}_{\Phi,K} = \bigcup_{[g] \in \mathrm{Vert} \mathcal{G}_{\Phi,K}} \mathcal{V}_{\Phi,K}([g]).$$
The classes $[g_i]$ are called the $\Phi-$neighbors of $[g]$ (relative to $K$). 

Fix an integer $n \geq 1$. The set of bi-$K$-invariant Hecke operators has the structure of a $\C$-algebra isomorphic to
\[\C[\Phi_{x,1}, \ldots, \Phi_{x,n}, \Phi_{x,n}^{-1}]_{x \in |X|},\]
where $x$ is a place of $F$ and $\Phi_{x,r}$ is the characteristic function of
$$K\left( \begin{array}{cc}
\pi_x I_r &  \\ 
 & I_{n-r}
\end{array}    \right)K,$$  
where $I_k$ is the $k \times k$ identity matrix. We denote $\mathcal{V}_{\Phi_{x,r},K}([g])$ by $\mathcal{V}_{x,r}([g])$ and  $\mathcal{G}_{\Phi_{x,r},K}$ by $\mathcal{G}_{x,r}.$

\subsection*{Graph of a Hecke operator: algebraic geometric perspective} 
Our goal is to describe the graphs of Hecke operators for an elliptic function field. Throughout this work, we consider the geometric point of view of graphs of Hecke operators. By a theorem due to Weil, there is a bijection between $ \mathrm{GL}_n(F) \setminus \mathrm{GL}_n(\mathbb{A}) / K$ with the set $\mathrm{Bun}_n X$ of isomorphism classes of rank $n$ vector bundles on  $X$. This theorem allows us to determine the action of an unramified Hecke operator $\Phi_{x,r}$ in terms of the equivalence classes of short exact sequences of coherent sheaves on $X.$ Namely, we consider exact sequences of the form
$$0 \longrightarrow \E' \longrightarrow \E \longrightarrow \mathcal{K}_{x}^{\oplus r} \longrightarrow 0$$  
where $\E',\E$ are rank $n$-vector bundles, $x$ is a closed point of $X,$ and $\mathcal{K}_{x}^{\oplus r}$ is the skyscraper sheaf on $x$ whose stalk is $\kappa(x)^{\oplus r}$. Let $m_{x,r}(\E,\E')$ be the number of isomorphism classes of exact sequences 
$$0 \longrightarrow \E'' \longrightarrow \E \longrightarrow \mathcal{K}_{x}^{\oplus r} \longrightarrow 0$$
with fixed $\E$ such that $\E'' \cong \E'.$ Equivalently, $m_{x,r}(\E,\E')$ is the number of subbundles $\E''$ of $\E$ that are isomorphic to $\E'$ and  which the quotient $\E/\E''$ is isomorphic to $\mathcal{K}_{x}^{\oplus r}$. We denote by  $\mathcal{V}_{x,r}(\E)$ the set of triples $\big(\E,\E',m_{x,r}(\E,\E')\big)$ such that there exists an exact sequence of the type as  above, i.e.\  $m_{x,r}(\E,\E') \neq 0$.  Then,
$$\mathrm{Vert}\;\mathcal{G}_{x,r} = \mathrm{Bun}_n X \hspace{0.3cm}\text{ and } \hspace{0.3cm} \mathrm{Edge}\; \mathcal{G}_{x,r} = \coprod_{\E \in \mathrm{Bun}_n X} \mathcal{V}_{x,r}(\E),$$
see Theorem 3.4 in \cite{roberto-graphs}.

\subsection*{Hall algebras}
The Hall algebra $\mathbf{H}_X$ of coherent sheaves on a smooth projective curve $X,$ as introduced by Kapranov in \cite{kapranov}, encodes the extensions of coherent sheaves. Let $v$ be a square root of $q^{-1}$. The Hall algebra of $X$ is the  vector space $\mathbf{H}_X := \bigoplus_{\F \in \mathrm{Coh}(X)} \C \F$ with the product
$$\F \;\G = v^{-\langle F,G \rangle} \sum_{\mathcal{H}} h_{\F,\G}^{\mathcal{H}} \mathcal{H}$$
where $\langle F,G \rangle := \dim_{\Fq} \mathrm{Ext}^{0}(\F,\G) -   \dim_{\Fq} \mathrm{Ext}^{1}(\F,\G)$ and 
\[ h_{\F,\G}^{\mathcal{H}} := \frac{\#\big\{ 0 \longrightarrow \G \longrightarrow \mathcal{H} \longrightarrow \F \longrightarrow 0\big\}}{\#\mathrm{Aut}(\F) \; \#\mathrm{Aut}(\G)} \]

The main observation, which links the theory of Hall algebras with graphs of Hecke operators, is that we can recover the multiplicities $m_{x,r}(\E,\E')$ from the product $\mathcal{K}_{x}^{\oplus r}\;\E'$ in the Hall algebra of $X$. As we show in the Lemma \ref{lemmaconnection}, the quantities $m_{x,r}(\E,\E')$ and  $h_{\mathcal{K}_{x}^{ r},\E'}^{\E}$ are equals. Thus, for a fixed $n$, the graphs of Hecke operators can be described by calculating explicitly the products $\mathcal{K}_{x}^{\oplus r}\E'$ where $\E'$ runs through the set of rank $n$ vector bundles on $X.$

By what we have explained in the previous paragraph, our problem is reduced to calculating the products $\mathcal{K}_{x}^{\oplus r}\E$ for all $\E \in \mathrm{Bun}_n X.$ To do so, we use some structure results from Burban and Schiffmann \cite{olivier-elliptic1} and Fratila \cite{dragos}  about the elliptic Hall algebra. Our strategy is to make a ''base change'' and write the product $\mathcal{K}_{x}^{\oplus r}\E$ in terms of elements in some subalgebras of the whole Hall algebra, called \textit{twisted spherical Hall algebras} (see Definition \ref{deftwistedsphericalalgebra}). Since the twisted spherical  Hall algebras are well understood and have a characterization in terms of  \textit{path algebras} (see Theorem \ref{dragosmaintheorem}), the base change allows us to explicitly calculate these products. However, these calculations depend closely on our initial data, i.e.\ the degree of $x$, the choice of $r$ and the vector bundle $\E.$ Hence, what we do in the section \ref{sectionthealgorithm}, is to develop an \textit{algorithm}, which  calculates these products step-by-step.  We outline the algorithm in the following.

\subsection*{The algorithm for elliptic function fields}
 
From this point on $X$ denotes an elliptic curve. In order to calculate the products $\mathcal{K}_{x}^{\oplus r}\E,$ we will use the twisted spherical Hall algebras, which are subalgebras of the whole Hall algebra. Let $\mathbf{Z} := \{ (r,d) \in \Z^2 \;|\; r >0 \text{ or } r=0 \text{ and } d> 0 \}$ and define for $(r,d) \in \mathbf{Z}$  the slope $\mu(r,d)$ of $(r,d)$ as $d/r$ if $r \neq 0$ and as $\infty$ if $r=0$.
Let $\sigma, \overline{\sigma}$ be the eigenvalues of the Frobenius automorphism acting on $H^1(\overline{X},\overline{\Q_\ell})$ and $d$ be a positive integer. The twisted spherical Hall algebra $\mathbf{E}_{\sigma, \overline{\sigma}}^{d}$ (cf. Definition \ref{defpathalgebra}) has an explicit description in terms of generators $T_{\mathbf{v}}^{\widetilde{\rho}}$ (cf. Definition \ref{deftwistedsphericalalgebra}) and relations (see Theorem \ref{dragosmaintheorem}), where $\mathbf{v} \in \mathbf{Z}$, $\widetilde{\rho}$ is a character on $\mathrm{Pic}^{0}X_n$ modulo the Frobenius action, and $X_n = X \times_{\mathrm{Spec}\; \Fq} \mathrm{Spec}\;\mathbb{F}_{q^n}$. These algebras were introduced by Burban and Schiffmann in \cite{olivier-elliptic1} and generalized by Fratila in \cite{dragos}. For further and precise definitions see section \ref{section3}.\\

\noindent 
\textbf{Input.}
Let $x \in X$ be a closed point, $\E \in \mathrm{Bun}_n X$ and $r$ an integer such that $1 \leq r \leq n.$ \\

\noindent 
\textbf{Base change.}
\begin{enumerate}
\item[(i)] Let $\mathsf{C}_{\mu}$ be the category of semistable coherent sheaves of slope $\mu.$  Write
$$\E = \E_1 \oplus \cdots \oplus \E_s$$  
where $\E_i \in \mathsf{C}_{\mu_i}$ and
$\mu_1 < \cdots < \mu_s$ (Harder-Narasimhan decomposition).  
\item[(ii)] Let $\mathrm{Tor}(X)$ be the category of torsion sheaves. For each $\E_i \in \mathsf{C}_{\mu_i}$, use the equivalence $\mathsf{C}_{\mu_i} \equiv \mathrm{Tor}(X)$ (cf. Theorem \ref{atiyahtheorem}) to associate $\mathcal{T}_i \in \mathrm{Tor}X$ with $\E_i$ and write $\mathcal{T}_i = \mathcal{T}_{x_{i1}} \oplus \cdots \oplus \mathcal{T}_{x_{im}}$ with $\mathcal{T}_{x_{ij}} \in \mathrm{Tor}_{x_{ij}}$ and $x_{ij} \neq x_{i{j'}} $ for $j \neq j'$  
\item[(iii)]Each $\mathcal{T}_{x_{ij}}$ corresponds to a Hall-Littlehood symmetric function $P_{\lambda_{ij}}$ in the Macdonald ring of symmetric functions. Write $P_{\lambda_{ij}}$ as sum of products of power-sums functions.
\item[(iv)] Taking the inverse image of the power-sums yields an expression of $\mathcal{T}_{x_{ij}}$ in terms of sums of products of 
$$T_{(0,m),x_i} := \frac{[m]|x_i|}{m} \sum_{|\lambda|=m/|x_i|} n_{u_{x_i}}(l(\lambda)-1) \mathcal{K}_{x_i}^{(\lambda)}$$
where  $\mathcal{K}_{x_i}^{(\lambda)}$ is the unique torsion sheaf with support at $x_i$ associated to the partition $\lambda.$ See Section \ref{section3} for the definitions of the constants in the above definition. 
By Proposition \ref{inversionlemma}, we can write $\mathcal{T}_{x_{ij}}$ as a sum of products of elements $$T_{\mathbf{v}}^{\widetilde{\rho}} := \sum_{x \in X} \widetilde{\rho}(x) T_{\mathbf{v},x}$$  
where $\widetilde{\rho}$ is a character in $\mathrm{Pic}^{0}X_d$ modulo the Frobenius action, for some extension of base field $X_d$ of $X$ and $\mu(\mathbf{v}) = \infty.$ 
\item[(v)] Since $\E = a \; \E_1 \cdots \E_s$ for some $a \in \C$ and $\mathcal{T}_i = \mathcal{T}_{x_{i1}}  \cdots  \mathcal{T}_{x_{im}}$, we may write  $$\E =  \sum_{\widetilde{\rho}_{ij_k}}^{} a_{ij_k} \;T_{\mathbf{v}_{1j_1}}^{\widetilde{\rho}_{1j_1}} \cdots T_{\mathbf{v}_{1j_{k_j}}}^{\widetilde{\rho}_{1 j_{k_j}}} \cdots T_{\mathbf{v}_{sj_1}}^{\widetilde{\rho}_{s j_1}} \cdots T_{\mathbf{v}_{s j_{k_j}}}^{\widetilde{\rho}_{s j_{k_j}}},$$
for some $a_{ij_k} \in \C$, $\widetilde{\rho}_{ij_k}$ orbits of characters in some Picard groups.
\item[(vi)] Therefore,   
\[ \mathcal{K}_{x}^{\oplus r} \; \E =  \sum_{i=1}^{m} a_{i}\; T_{\mathbf{v}_{i_0}}^{\widetilde{\rho}_{i_0}} T_{\mathbf{v}_{i_1}}^{\widetilde{\rho}_{i_1}} \cdots T_{\mathbf{v}_{i_{\ell}}}^{\widetilde{\rho}_{i_{\ell}}}\]
for some $a_i \in \C$ where $\widetilde{\rho}_{i_j}$ are orbits of characters in some Picard groups  and where  $\mathbf{v}_{i_j} \in \mathbf{Z}$ are such that $\mu(\mathbf{v}_{i_0}) = \infty$ and $\mu(\mathbf{v}_{i_1}) \leq \cdots \leq \mu(\mathbf{v}_{i_\ell})$. 
\end{enumerate}

\noindent 
\textbf{Order by slopes}
\begin{enumerate}
\item[(i)] Aiming to use the structure of $\mathbf{E}_{\sigma, \overline{\sigma}}^{d}$ (cf. Definition \ref{defpathalgebra}) to write the above product in increasing order of slopes, we consider the commutator  $\big[\mathcal{K}_{x}^{\oplus r } , \E\big]$ since $\pi^{\mathrm{vec}}(\mathcal{K}_{x}^{\oplus r} \; \E) = \pi^{\mathrm{vec}}\big[\mathcal{K}_{x}^{\oplus r } , \E\big].$ Where $\pi^{\mathrm{vec}}$ means that we are considering the vector bundles which appear in the product.   
\item[(ii)] Observe that $\big[ T_{\mathbf{v}_{i}}^{\widetilde{\rho}_{i}}, T_{\mathbf{v}_{j}}^{\widetilde{\rho}_{j}} \big] =0$ unless  $T_{\mathbf{v}_{i}}^{\widetilde{\rho}_{i}}, T_{\mathbf{v}_{j}}^{\widetilde{\rho}_{j}}$ belong to same twisted spherical Hall algebra (cf. Definition \ref{deftwistedsphericalalgebra}). By Theorem \ref{dragosmaintheorem}, the problem of ordering  the slopes is reduced to $\mathbf{E}_{\sigma, \overline{\sigma}}^{d}.$
\item[(iii)] Using the structure of $\mathbf{E}_{\sigma, \overline{\sigma}}^{d}$, we can calculate  $\big[ T_{\mathbf{v}}^{\widetilde{\rho}}, T_{\mathbf{w}}^{\widetilde{\rho}} \big] $ explicitly. If $\mu(\mathbf{v}) > \mu(\mathbf{w})$ we replace $ T_{\mathbf{v}}^{\widetilde{\rho}} T_{\mathbf{w}}^{\widetilde{\rho}} $ by $ \big[T_{\mathbf{w}}^{\widetilde{\rho}}, T_{\mathbf{v}}^{\widetilde{\rho}}  \big] + T_{\mathbf{w}}^{\widetilde{\rho}} T_{\mathbf{v}}^{\widetilde{\rho}}.$
\item[(iv)] Combining the previous steps we may write 
\[ \mathcal{K}_{x}^{\oplus r} \; \E =  \sum_{i=1}^{m} a_{i}\;  T_{\mathbf{v}_{i_1}}^{\widetilde{\rho}_{i_1}} \cdots T_{\mathbf{v}_{i_{\ell}}}^{\widetilde{\rho}_{i_{\ell}}}\]
where $\mu(\mathbf{v}_{i_1}) \leq \cdots \leq \mu(\mathbf{v}_{i_\ell}).$
\end{enumerate}

\noindent 
\textbf{Base change back and calculations}
\begin{enumerate}
\item[(i)]  We replace $ T_{\mathbf{v}_{i_j}}^{\widetilde{\rho}_{i_j}}$ by its definition $\sum_{x'} \widetilde{\rho}_{i_j}(x') T_{\mathbf{v}_{i_j}, x'}.$ The elements  $T_{\mathbf{v}_{i_j}, x'}$ are given explicitly by sums of semistable sheaves of slope $\mu(\mathbf{v}_{i_j})$ via Atiyah's classification (cf. Theorem \ref{atiyahtheorem}). Therefore the problem reduces to calculating the product $T_{\mathbf{v}_1,x_1} \; T_{\mathbf{v}_2,x_2}$ with $\mu(\mathbf{v}_1) \leq \mu(\mathbf{v}_2).$  
\item[(ii)] If $\mu(\mathbf{v}_1) < \mu(\mathbf{v}_2)$ or $x_1 \neq x_2$, then 
\[T_{\mathbf{v}_1,x_1} \; T_{\mathbf{v}_2,x_2} = T_{\mathbf{v}_1,x_1} \oplus T_{\mathbf{v}_2,x_2}\]
and we are able to calculate explicitly the product $\pi^{\mathrm{vec}}(\mathcal{K}_{x}^{\oplus r} \; \E).$
\item[(iii)] If $\mu(\mathbf{v}_1) = \mu(\mathbf{v}_2)$ and $x_1 = x_2$, the product $T_{\mathbf{v}_1,x_1} \; T_{\mathbf{v}_2,x_2}$ can be calculated as the product of power-sums in the Macdonald ring of symmetric functions (cf. Proposition \ref{prop-macdonald}). In this case, we write the product $T_{\mathbf{v}_1,x_1} \; T_{\mathbf{v}_2,x_2}$ in the Macdonald ring of symmetric functions as linear combination of Hall-Littlehood functions $\sum_i P_{\lambda_i}$ and come back to the Hall algebra via Proposition \ref{prop-macdonald}, which associates $P_{\lambda_i}$ with the unique semistable sheaf in $\mathsf{C}_{\mu(\mathbf{v}_1)}$ that corresponds to $\mathcal{K}_{x_1}^{(\lambda_i)}$ via $\mathsf{C}_{\mu(\mathbf{v}_1)} \equiv \mathrm{Tor} X.$ 
\end{enumerate}

Let us briefly describe the content of this article. In the first section, we set up the notation, recall Atiyah's classification of coherent sheaves on an elliptic curve and the geometric definition of graphs of Hecke operators. We establish in section \ref{section2} the relation between graphs of Hecke operators and the Hall algebra of an arbitrary smooth projective curve over a finite field.  
In the third section, we specialise the theory of the previous section to elliptic curves, derive some properties for the graphs of Hecke operators  and review useful results by Fratila, Burban and Schiffmann. Section \ref{sectionthealgorithm} is the main part of this paper. In this section, we prove some structure theorems and describe the algorithm to calculate the graphs of Hecke operators for an elliptic curve. In section \ref{section-structureconstants}, we use our algorithm to deduce some general structure results. In section \ref{rank2}, we determine the complete graph $\mathcal{G}_{x,1}$ for $n=2$ and $|x|=1.$ Some applications are indicated. 
\\

\section{Background}\label{section1}

\paragraph*{\textbf{Notation.}}\label{paragraphnotation}

Let $\mathbb{F}_q$ be a finite field with $q$ a prime power. Throughout the article, $X$ denotes an elliptic curve defined over $\mathbb{F}_q;$ that is, $X$ is a smooth and geometrically irreducible projective curve of genus one having a rational point. Note that, by Hasse-Weil's inequality (we refer to \cite{silverman} for a proof), we have $|\#X(\mathbb{F}_q) - (q+1)| \leq 2 \sqrt{q},$
hence any smooth projective curve of genus one has such a point, say $x_0 \in X(\mathbb{F}_q)$. We denote $\#X(\mathbb{F}_{q^d})$ by $N_d.$

Let $F$ be the function field of $X$. We denote by $|X|$ the set of closed points of $X$ or, equivalently, the set of places of $F$. For $x \in |X|$, we denote by $F_x$ the completion of $F$ at $x$, by $\mathcal{O}_x$ its integers, by $\pi_x \in \mathcal{O}_x$ a uniformizer (we can suppose $\pi_x \in F$)  and by $q_x$ the cardinality of the residue field $\kappa(x):=\mathcal{O}_x/(\pi_x) \cong \mathbb{F}_{q_x}.$ Let $|x|$ be the degree of $x$ which is defined by the extension field degree $[\kappa(x) : \Fq]$; in other words, $q_x = q^{|x|}$.

For an extension of finite fields $\mathbb{F}_{q^n}$ of $\mathbb{F}_q$ we will denote by $X_n$ the fiber product $X \times_{\mathrm{Spec}(\mathbb{F}_q)} \mathrm{Spec}(\mathbb{F}_{q^n}).$ Using our choice of a rational point $x_0 \in X(\Fq)$, the Riemann-Roch theorem yields a bijective map 

\begin{eqnarray*}
X(\mathbb{F}_{q^n}) = X_n(\mathbb{F}_{q^n}) & \stackrel{1:1}{\longrightarrow}  & \mathrm{Pic}^{0}(X_n) \\
x  & \longmapsto & \mathcal{O}_{X_n}(x-x_0) 
\end{eqnarray*}
Therefore we can transport the group  structure  from $\mathrm{Pic}^{0}(X_n) $ to $ X_n(\mathbb{F}_{q^n})$  such that $x_0$ is the neutral element. Moreover, the inclusions $ X(\mathbb{F}_{q^m}) \subseteq  X(\mathbb{F}_{q^n})$ for $m|n$ are compatible with the respective group structure. To avoid confusion with the addition of divisors,  we will denote the sum of $x,y \in X(\mathbb{F}_{q^n})$ with respect to the group law by $x \oplus y \in  X(\mathbb{F}_{q^n}).$

Let $Y$ be a scheme over $\Fq.$ We call the morphism $\mathrm{Fr}_Y : Y \rightarrow Y$ induced by the sheaf homomorphism $\mathcal{O}_Y \rightarrow \mathcal{O}_Y : a \mapsto a^q$ the Frobenius endomorphism of $Y.$ This morphism is the identity on the underlying topological spaces, and it raises functions to the $q$-th power (\cite{liu} 3.2.4). 

The Frobenius endomorphism acts on the set $X(\mathbb{F}_{q^n})$ of $\mathbb{F}_{q^n}$-rational points of $X$ by composition $\mathrm{Fr}_X \circ x : \mathrm{Spec}(\mathbb{F}_{q^n}) \rightarrow X,$ where $x \in X(\mathbb{F}_{q^n}).$ Moreover, this action is compatible with the group structure of $X(\mathbb{F}_{q^n}).$ If $m|n$ are positive integers, then we have an obvious identification $X(\mathbb{F}_{q^m}) = X(\mathbb{F}_{q^n})^{\mathrm{Fr}_{X}^m}.$ In particular $X(\mathbb{F}_{q^n}) = X(\overline{\Fq})^{\mathrm{Fr}_{X}^{n}}.$ We have an identification of sets $|X| \cong X(\overline{\Fq})\big/ \mathrm{Fr}_X$ where the quotient means that two $\overline{\Fq}$-rational points of $X$ are identified if they have the same orbit under the Frobenius action. Similarly, we have $|X_n| \cong X(\overline{\Fq})\big/ \mathrm{Fr}_{X}^n.$ These two equalities allow us to define an action of $\mathrm{Fr}_X$ on the closed points of $X_n.$ We denote this action by $\mathsf{frob}_n: |X_n| \rightarrow |X_n|.$ In particular $|X_n|\big/\mathsf{frob}_n = |X|.$

For a finite abelian group $G$ we denote by $\widehat{G}$ its group of characters. 
From the isomorphism $X(\mathbb{F}_{q^n})  \cong \mathrm{Pic}^{0}(X_n)$ that we fixed before, we obtain the following commutative diagram
\[ \xymatrix@R5pt@C7pt{ \mathrm{Pic}^{0}(X_n) \ar[dd]^{\wr} \ar[rr]^{\mathrm{Fr}_{X,n}^{*}} & & \mathrm{Pic}^{0}(X_n) \ar[dd]^{\wr}\\ 
& & \\
X(\mathbb{F}_{q^n}) \ar[rr]^{\mathrm{Fr}_{X}} & & X(\mathbb{F}_{q^n}).} \]
The Frobenius $\mathrm{Fr}_{X,n}^{*}$ acts by duality on each group $\widehat{\mathrm{Pic}^{0}(X_n)}$ and we will denote this action simply by $\mathrm{Fr}_{X,n}.$ We will  denote by $\widetilde{\mathrm{Pic}^{0}(X_n)}$ (or in some cases only $\mathrm{P}_n$) the quotient $\widehat{\mathrm{Pic}^{0}(X_n)} \big/ \mathrm{Fr}_{X,n}.$

Let $m,n$ be positive integers such that $m|n.$ We define the relative norm map
\[ \mathrm{Norm}_{m}^{n}:  \mathrm{Pic}^{0}(X_n) \longrightarrow \mathrm{Pic}^{0}(X_m) \] 
by
\[ \mathrm{Norm}_{m}^{n}(\Line) := \bigotimes_{i=0}^{n/m -1} (\mathrm{Fr}_{X,n})^{mi}(\Line) \]  
The fact that this map is well defined follows from Galois descent i.e.\ a line bundle on $\overline{X}$ that it is isomorphic to its Frobenius conjugates, descends to a line bundle on $X.$ For more details see Proposition 3.4 in \cite{arason}. 
By dualizing we obtain relative norm maps between the character groups for which we will use the same notation
\[ \mathrm{Norm}_{m}^{n}:  \widehat{\mathrm{Pic}^{0}(X_n)} \longrightarrow \widehat{\mathrm{Pic}^{0}(X_m)}. \] 

Recall that we fixed a rational point $x_0$ on $X$ and therefore for any integers $n \geq 1$ and $d \in \Z$ we can identify $\mathrm{Pic}^{d}(X_n)$ with $\mathrm{Pic}^{0}(X_n)$ by subtraction $d x_0.$ This allows us to extend any character $\rho$ of $\mathrm{Pic}^{0}(X_n)$ to a character of $\mathrm{Pic}(X_n)$ by putting $\rho (\mathcal{O}_{X_n}(x_0))=1.$ Unless otherwise specified we will view the characters of $\mathrm{Pic}^{0}(X_n)$ as characters of $\mathrm{Pic}(X_n).$ 

Denote $\widetilde{X} := \coprod_n \widetilde{\mathrm{Pic}^{0}(X_n)}.$ By abuse of language, we call the elements of $\widetilde{X}$ characters even if they are actually orbits of characters. For a character $\rho \in \widetilde{X}$, we say that it is a character of degree $n$ if $\rho \in \widetilde{\mathrm{Pic}^{0}(X_n)}.$ \\ 

\paragraph*{\textbf{Coherent sheaves on elliptic curve}}\label{atiyahparagraph}

We denote by $\mathrm{Coh}(X)$  the category of coherent sheaves on $X.$ Let us first outline, following Atiyah, the classification of coherent sheaves on elliptic curves (in \cite{atiyah-elliptic} it is assumed that the ground field is algebraically closed, but the proof applies to an arbitrary field). We denote by $\rk (\F)$ (resp. $\deg(\F)$) the rank (resp. degree) of a coherent sheaf $\F$. Recall that the \textit{slope} of a nonzero sheaf $\F \in \mathrm{Coh}(X)$ is $\mu(\F) = \deg(\F)\big/  \rk(\F)$ and that a sheaf $\F$ is \textit{semistable} (resp. \textit{stable}) if for any proper nontrivial subsheaf $\G \subseteq \F$ we have $\mu(\G) \leq \mu(\F)$ (resp.$ \mu(\G) < \mu(\F)$). The full subcategory $\mathsf{C}_{\mu}$ of $\mathrm{Coh}(X)$ consisting of all semistable sheaves of a fixed slope $\mu \in \Q \cup \{\infty\}$ is abelian, artinian and closed under extension. Moreover, if $\F,\G$ are semistable with $\mu(F) < \mu(G)$, then $\mathrm{Hom}(\G,\F)= \mathrm{Ext}(\F,\G)=0.$ Any sheaf $\F$ possesses a unique filtration (the Harder-Narasimhan filtration, or HN-filtration see \cite{hardernarasimhan})
$$0 = \F^{r+1} \subset \F^{r} \subset \cdots \subset \F^{1}=\F$$ 
for which $\F^{i}\big/\F^{i+1}$ is semistable of slope $\mu_i$ and $\mu_1 < \cdots <\mu_r.$ Observe that $\mathsf{C}_{\infty}$ is the category of torsion sheaves and hence equivalent to the product category $\prod_{x} \mathrm{Tor}_x$ where $x$ runs through the set of closed points of $X$ and $\mathrm{Tor}_x$ denotes the category of torsion sheaves supported at $x.$ Since $\mathrm{Tor}_x$ is equivalent to the category of finite length modules over the local ring $\mathcal{O}_{X,x}$ of the point $x,$ there is a unique simple sheaf $\mathcal{K}_x$ in $\mathrm{Tor}_x,$ namely the skyscraper in $x$ with stalk $\kappa(x).$

\begin{theorem}[\cite{atiyah-elliptic} Theorem 7] \label{atiyahtheorem}The following hold.

$\mathrm{(i)}$ The HN-filtration of any coherent sheaf splits (noncanonically). In particular, any indecomposable coherent sheaf is semistable.

$\mathrm{(ii)}$ The set of stable sheaves of slope $\mu$ is the class of simple objects of $\mathsf{C}_{\mu}.$

$\mathrm{(iii)}$ The choice of any rational point $x_0 \in X(\mathbb{F}_q),$ induces an exact equivalence of abelian categories $\epsilon_{\nu,\mu}: \mathsf{C}_{\mu} \rightarrow \mathsf{C}_{\nu},$ for any $\mu,\nu \in \Q \cup \{\infty\}.$
\end{theorem}

The proof of Atiyah also provides an algorithm to compute the equivalences in Theorem \ref{atiyahtheorem} (iii). His proof is for an algebraically closed field (of any characteristic). In \cite{olivier-elliptic1}, Burban and Schiffmann prove these equivalences for the finite field case. We briefly recall their treatment; for complete details see \cite{olivier-elliptic1}, section 1 (paragraph 1.2) and Appendix A. We also refer \cite{burban} for a introductory discussion of the following statements.

The Grothendieck group $K_0(X)$ of $\mathrm{Coh}(X)$ is isomorphic to $\Z \oplus \mathrm{Pic}(X)$ and the isomorphism is given by $\F \mapsto \big(\mathrm{rk}(\F), \mathrm{det}(\F)\big)$. Moreover, if we compose this morphism with the one sending a line bundle to its degree, then we get a group homomorphism 
\begin{eqnarray*}
\varphi : K_0(X) &\longrightarrow& K_{0}'(X):=\Z^2 \\
\F &\longmapsto& \big(\mathrm{rk}(\F), \deg(\F)\big)
\end{eqnarray*}
The group $K_{0}'(X)$ is called the numerical Grothendieck group and for a sheaf $\F,$ we  denote by $\overline{\F}$ its image in this group. 
The Euler bilinear form $\langle ,\rangle: K_0(X) \otimes K_0(X) \rightarrow \Z$ is defined by the formula 
$$\overline{\F} \otimes \overline{\G} \mapsto \dim \mathrm{Hom}(\F,\G) - \dim \mathrm{Ext}(\F,\G).$$
By the Riemann-Roch theorem, one has 
$$\langle \overline{\F}, \overline{\G} \rangle = \rk(\F)\deg(\G) - \rk(\G) \deg(\F),$$
in particular the Euler form is skew-symmetric.

The kernel of $\varphi$ is the radical of the Euler form, which is given by 
$$\mathrm{rad}\langle -,-\rangle := \{\F \in K_{0}(X) | \rk(\F)=\deg(\F)=0\}.$$ 
Hence we obtain an isomorphism
\[\overline{\varphi} =(\rk, \deg) : K_{0}(X)\big/ \mathrm{rad}\langle -,-\rangle \stackrel{\sim}{\longrightarrow} \Z^2, \]
and the Euler form is well defined on classes in $K_{0}'(X).$

The group of exact auto-equivalences $\mathrm{Aut}(\mathcal{D}^{b}(\mathrm{Coh}(X)))$ of the bounded derived category of coherent sheaves on $X$ acts on $K_{0}(X)$ by automorphism, since $K_0(X)$ is  an invariant of the derived category. This action preserves the Euler form and its radical. Hence, we get a group homomorphism
$$\pi: \mathrm{Aut}(\mathcal{D}^{b}(\mathrm{Coh}(X))) \longrightarrow \mathrm{SL}_2(\Z)$$
which sends an auto-equivalence $\mathbb{F}$ to the upper map in the following commutative diagram
\[\xymatrix@R5pt@C7pt{ \Z^2 \ar[rr]^{\pi(\mathbb{F})} \ar[dd]^{\varphi^{-1}}_{\cong} & & \Z^2  \\
& & \\
K_{0}(X)\big/ \mathrm{rad}\langle -,-\rangle \ar[rr]_{\mathbb{F}_{*}}^{\cong} & & \ar[uu]^{\cong}_{\varphi} K_{0}(X)\big/ \mathrm{rad}\langle -,-\rangle
}.\]
which actually lies in $\mathrm{SL}_2(\Z)$.

Let $\F$ be indecomposable in $\mathrm{Coh}(X)$ with $\overline{\F} = \varphi(\F)=(n,d) \in \Z^2.$  Let $g := \gcd(n,d)$ be the greatest common divisor. If $d \neq 0$, then there is a matrix $f \in \mathrm{SL}_2(\Z)$ with $f(n,d) = (0,g).$ If $d=0$, then let $f$ be the matrix ${\tiny \begin{pmatrix}
0 & -1 \\ 
1 & 0
\end{pmatrix}} $   in $\mathrm{SL}_2(\Z)$ which flips the coordinate axes.  Since $\pi: \mathrm{Aut}(\mathcal{D}^{b}(\mathrm{Coh}(X))) \longrightarrow \mathrm{SL}_2(\Z)$ is surjective, we can lift $f$ to an auto-equivalence $\mathbb{F} \in \mathrm{Aut}(\mathcal{D}^{b}(\mathrm{Coh}(X))).$ Since $\mathbb{F}$ is an equivalence of categories the image $\mathbb{F}(\F)$ is indecomposable again. The object $\mathbb{F}(\F)$ is then a coherent sheaf of rank $0.$ Therefore, $\mathcal{T}=\mathbb{F}(\F)$ is a torsion sheaf and $\F = \mathbb{F}^{-1}(\mathcal{T}).$ 

The above discussion shows us how Atiyah's theorem works for the classification of indecomposable vector bundles on $X,$ and more generally the equivalences $\mathsf{C}_{\mu} \cong \mathsf{C}_{\nu}$, for any $\mu,\nu \in \Q \cup \{\infty \}.$ \\

\paragraph*{\textbf{Reminder on graphs of Hecke operators }}\label{remindergraphs}

We recall the definition of the graphs $\mathcal{G}_{x,r}$ of an unramified Hecke operator as introduced in \cite{oliver-graphs} for $\mathrm{PGL}_2$ and generalized to $\mathrm{GL}_n$ in \cite{roberto-graphs}. These graphs encode the action of the unramified Hecke operators $\Phi_{x,r}$  on the space of automorphic forms for $\mathrm{GL}_n$ over a global function field. We concentrate in this resume on the geometric point of view. For the translation into adelic language, see section 1 in \cite{roberto-graphs}.

We say that two exact sequences of sheaves
$$0 \longrightarrow \F_1 \longrightarrow \F \longrightarrow \F_{2} \longrightarrow 0 \hspace{0.2cm}\text{ and } \hspace{0.2cm} 0 \longrightarrow \F_{1}' \longrightarrow \F \longrightarrow \F_2' \longrightarrow 0$$  
are \textit{isomorphic with fixed $\F$} if there are isomorphism $\F_1 \rightarrow \F_1'$ and $\F_2 \rightarrow \F_2'$ such that 
$$\xymatrix@R5pt@C7pt{ 0 \ar[rr] && \F_1 \ar[rr] \ar[dd]^{\cong} && \F \ar[rr] \ar@{=}[dd] && \F_2 \ar[rr] \ar[dd]^{\cong} && 0 \\
&& && && && \\
0 \ar[rr] && \F_1' \ar[rr] && \F \ar[rr] && \F_2' \ar[rr] && 0 } $$
commutes. 
Let $\mathcal{K}_{x}$ be the torsion sheaf that is supported at $x$ and has stalk $\kappa(x)$ at $x,$ i.e. the skyscraper torsion sheaf at $x$. Fix $\E \in \mathrm{Bun}_n X$. For $r \in \{1, \ldots,n\},$ and $\E' \in \mathrm{Bun}_n X$ we define $m_{x,r}(\E,\E')$ as the number of isomorphism classes of exact sequences
$$0 \longrightarrow \E'' \longrightarrow \E \longrightarrow \mathcal{K}_{x}^{\oplus r} \longrightarrow 0$$
with fixed $\E$ and with $\E'' \cong \E'.$ We will denote $\mathcal{K}_{x}^{\oplus r}$ simply by $\mathcal{K}_{x}^{ r}$. 

\begin{definition} \label{defconection} Let $x \in |X|$. For a vector bundle $\E \in \mathrm{Bun}_n X$ we define
$$\mathcal{V}_{x,r}(\E) := \{(\E, \E', m) | m = m_{x,r}(\E,\E') \neq 0\},$$
and we call $\E'$ a $\Phi_{x,r}$-neighbour of $\E$ if $m_{x,r}(\E,\E') \neq 0$, and $m_{x,r}(\E,\E')$ its multiplicity.
\end{definition}

We define the graph $\mathcal{G}_{x,r}$ by 
\[ \mathrm{Vert}\mathcal{G}_{x,r} = \mathrm{Bun}_n X \quad \text{ and } \quad
  \mathrm{Edge} \mathcal{G}_{x,r} = \coprod_{\E \in \mathrm{Bun}_n X} \mathcal{V}_{x,r}(\E). \] 

If $\mathcal{V}_{x,r}(\E) = \{(\E, \E_1, m_1), \ldots, (\E, \E_r, m_r)\}$, we make the following drawing conventions to  illustrate the graph $\mathcal{G}_{x,r}$: vertices are represented by 
labelled dots, and an edge $(\E,\E',m)$ together with its origin $\E$ and its terminus $\E'$ is drawn as
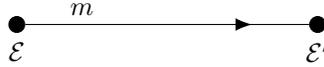
\[
  \beginpgfgraphicnamed{tikz/fig0}
  \begin{tikzpicture}[>=latex, scale=2]
        \vertex[circle,fill,label={below:$\E$}](00) at (0,0) {};
        \vertex[circle,fill,label={below:$\E'$}](10) at (2,0) {};
    \path[-,font=\scriptsize]
    (00) edge[->-=0.8] node[pos=0.2,auto,black] {\normalsize $m$} (10)
    ;
   \end{tikzpicture}
 \endpgfgraphicnamed
\]
The $\Phi_{x,r}$-neighbourhood of $\E$ is thus illustrated as
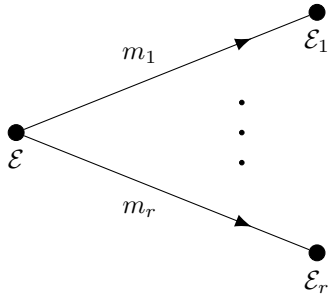
\[
  \beginpgfgraphicnamed{tikz/fig}
  \begin{tikzpicture}[>=latex, scale=2]
        \vertex[circle,fill,label={below:$\E$}](00) at (0,0) {};
        \vertex[circle,fill,label={below:$\E_1$}](11) at (2,0.8) {};
        \vertex[circle,fill,label={below:$\E_r$}](10) at (2,-0.8) {};
\draw (1.5,0.2) circle (0.015cm);
\fill  (1.5,0.2) circle (0.015cm);      
\draw (1.5,0) circle (0.015cm);
\fill  (1.5,0) circle (0.015cm);  
\draw (1.5,-0.2) circle (0.015cm);
\fill  (1.5,-0.2) circle (0.015cm);   
            \path[-,font=\scriptsize]
    (00) edge[->-=0.8] node[pos=0.5,auto,black,swap] {\normalsize $m_r$} (10)
    (00) edge[->-=0.8] node[pos=0.5,auto,black] {\normalsize  $m_1$} (11)
    ;
   \end{tikzpicture}
 \endpgfgraphicnamed
\]



\section{Hall algebras and graphs of Hecke operators}\label{section2}

\noindent
Let $Y$ be a smooth projective and geometric irreducible curve over a finite field $\mathbb{F}_q.$ We start with the definition of the Hall algebra of $Y$ (following \cite{kapranov} and \cite{ringel1}). Fix a square root $v$ of $q^{-1}$. Let $\mathbf{H}_Y$ be the $\C-$vector space 
$$\mathbf{H}_Y : = \bigoplus_{\F} \C\; \F$$
where $\F$ runs through the isomorphism classes of objects in $\mathrm{Coh}(Y).$ 
Let  $(\F,\G,\mathcal{H})$ be a triple of coherent sheaves on $Y.$ We denote by $P_{\F,\G}^{\mathcal{H}}$ the cardinality of the set of short exact sequences
$$0 \longrightarrow \G \longrightarrow \mathcal{H} \longrightarrow \F \longrightarrow 0.$$ 
Note that $P_{\F,\G}^{\mathcal{H}}$ is finite since $\mathrm{Coh}(Y)$ is a finitary category. Define
$$h_{\F,\G}^{\mathcal{H}} := \frac{P_{\F,\G}^{\mathcal{H}}}{\#\mathrm{Aut}(\F) \#\mathrm{Aut}(\G)}.$$ 
The following product defines on $\mathbf{H}_Y $ the structure of an associative algebra:
$$\F\; \G := v^{-\langle \F,\G \rangle}\sum_{\mathcal{H}} h_{\F,\G}^{\mathcal{H}} \;\mathcal{H}.$$


Another way of defining the Hall algebra is as follows. Consider the space $\mathrm{Fun}_0(\mathcal{M}_Y,\C)$ of $\C$-valued functions with finite support on the set $\mathcal{M}_Y$ of isomorphism classes of objects of $\mathrm{Coh(Y)}.$ This space of functions identifies naturally, as a vector space, with $\mathbf{H}_Y.$ We can endow $\mathrm{Fun}_0(\mathcal{M}_Y,\C)$ with a convolution product. For $f,g \in \mathrm{Fun}_0(\mathcal{M}_Y,\C)$ we define 
\[(f * g)(\F) = \sum_{\G \subseteq \F} v^{-\langle \F/\G, \G \rangle } f(\F/\G) g(\G). \]
Using the identification of $\mathbf{H}_Y$ with $\mathrm{Fun}_0(\mathcal{M}_Y,\C)$, we obtain the same structure of an algebra on $\mathbf{H}_Y$ as defined in the previous paragraph. See Lecture 1, section 1.3 in \cite{olivier-hall} for a proof of this fact and more details. 

The following lemma is the main connection between the theory of graphs of Hecke operators and the theory of Hall algebras of a smooth and projective curve.

\begin{lemma}\label{lemmaconnection} Let $Y$ be a smooth projective curve over $\mathbb{F}_q.$ If $\E,\E' \in \mathrm{Bun}_n Y,$ then
$$h_{\mathcal{K}_{x}^{ r},\E'}^{\E} = m_{x,r}(\E, \E')$$ 

\begin{proof} Let $A := \{(0 \longrightarrow \E'' \longrightarrow \E \longrightarrow \mathcal{K}_{x}^{\oplus r} \longrightarrow 0) \big| \E'' \cong \E'\}$ be the set of representatives classes of short exact sequences as in
\ref{remindergraphs}.  Note that $\# A = m_{x,r}(\E, \E').$ Define 
$$B := \big\{\E'' \subseteq \E \big| \E'' \cong \E' \text{ and } \E / \E'' \cong \mathcal{K}_{x}^{ r}\big\}.$$
Since $\mathrm{Aut}(\E') \times \mathrm{Aut}(\mathcal{K}_{x}^{ r})$ acts freely on set of exact sequences 
$$0 \longrightarrow \E' \longrightarrow \E \longrightarrow \mathcal{K}_{x}^{ r} \longrightarrow 0$$ 
we have $h_{\mathcal{K}_{x}^{ r},\E'}^{\E} = \# B.$ We claim that the map $A \rightarrow B$ given by
\[\big(0 \longrightarrow \E'' \longrightarrow \E \longrightarrow \mathcal{K}_{x}^{ r} \longrightarrow 0\big) \mapsto \E''\]
is a bijection. Indeed, if $\big(0 \longrightarrow \E'' \longrightarrow \E \longrightarrow \mathcal{K}_{x}^{ r} \longrightarrow 0\big)$ and $\big(0 \longrightarrow \E''' \longrightarrow \E \longrightarrow \mathcal{K}_{x}^{ r} \longrightarrow 0\big)$ have the same image via the above map, then $\E'' \cong \E' \cong \E'''$ and 
$$\xymatrix@R5pt@C7pt{ 0 \ar[rr] && \E'' \ar[rr] \ar[dd]^{\cong} && \E \\
&& && \\
0 \ar[rr] && \E''' \ar[rruu]&& }$$
commutes. This implies that these two short exact sequences are the same element in $A$,  thus we have a injection. The surjectivity is trivial. Therefore $\# A = \# B.$
\end{proof}
\end{lemma}

The Hall algebra $\mathbf{H}_Y$ has a natural grading in the subset 
$$\mathbf{Z} := \big\{(n,d) \in K_{0}'(Y) \big| n> 0 \text{ or } n=0 \text{ and } d \geq 0\big\}$$
of the numerical Grothendieck group $K_{0}'(Y)$,  given by $\mathbf{H}_Y[\mathbf{v}] := \bigoplus_{\overline{\F}=\mathbf{v}} \C [\F]$, for $\mathbf{v} \in \mathbf{Z}.$ For $\mathbf{v}=(n,d) \in \mathbf{Z}$, $\gamma(\mathbf{v})$ denotes the $\gcd(n,d).$

We will denote by $\mathbf{H}_{Y}^{\mathrm{tor}} := \bigoplus_{d \geq 0} \mathbf{H}_Y[0,d]$ the subalgebra of torsion sheaves and by $\mathbf{H}_{Y}^{\mathrm{vec}} := \bigoplus_{(n,d) \in \mathbf{Z}, n>0} \mathbf{H}_Y[n,d]$ the subalgebra of vector bundles. We denote by $\pi^{\mathrm{vec}}$ the projection map $\mathbf{H}_Y \rightarrow \mathbf{H}_{Y}^{\mathrm{vec}}.$ 

\begin{remark} \label{propconnection} It follows from the definition of the Hall algebra, the above lemma and the previous paragraph that the graph of the Hecke operator $\mathcal{G}_{x,r}$ is completely determined by $v^{-nr|x|}\pi^{\mathrm{vec}}(\mathcal{K}_{x}^{ r}\text{ } \E)$   where $\E$ runs through $\mathrm{Bun}_n Y.$
The product $v^{-nr|x|}\pi^{\mathrm{vec}}(\mathcal{K}_{x}^{ r}\text{ } \E)$ in the Hall algebra give us all the edges in $\mathcal{G}_{x,r}$ that arrive in the vertex given by $\E.$ 
\end{remark}

\begin{example} The following is an easy example for $Y= \mathbb{P}_{\mathbb{F}_q}^1$ that exemplifies the connection between the Hall algebra for $Y$ and the graph of a certain Hecke operator. In Theorem 13 \cite{kassel}  Baumann and Kassel prove that 
$$\mathcal{K}_{x}^{\oplus r} \text{ }\mathcal{O}_{\mathbb{P}^1}(d) = \mathcal{O}_{\mathbb{P}^1}(d+ |x|) \oplus \mathcal{K}_{x}^{\oplus r-1} + q_{x}^{r}\; \big(\mathcal{O}_{\mathbb{P}^1}(d) \oplus \mathcal{K}_{x}^{\oplus r}\big),$$ 
which  gives us 
\[
  \beginpgfgraphicnamed{tikz/fig1}
  \begin{tikzpicture}[>=latex, scale=2]
        \vertex[circle,fill,label={below:$\mathcal{O}_{\mathbb{P}^1}(d+ |x|)$}](00) at (0,0) {};
        \vertex[circle,fill,label={below:$\mathcal{O}_{\mathbb{P}^1}(d)$}](10) at (2,0) {};
    \path[-,font=\scriptsize]
    (00) edge[->-=0.8] node[pos=0.2,auto,black] {\normalsize $1$} (10)
    ;
   \end{tikzpicture}
 \endpgfgraphicnamed
\]
for $r=1$. Compare this with section 5 (Figure 5) in \cite{roberto-graphs}.
\end{example}

\begin{lemma}\label{lemmacommutator} Let $x \in |Y|$. If $\E \in \mathrm{Bun}_n Y$, then 
$$\pi^{\mathrm{vec}}(\mathcal{K}_{x}^{ r}\text{ } \E) = \pi^{\mathrm{vec}}[\mathcal{K}_{x}^{r}, \E],$$
where the commutator is taken in the Hall algebra $\mathbf{H}_Y.$ 

\begin{proof} This follows from the fact that none of the elements in $\mathrm{Ext}(\E, \mathcal{K}_{x}^{r})$ belongs to $\mathbf{H}_{Y}^{\mathrm{vec}}.$ 
\end{proof}
\end{lemma}

\begin{lemma}[\cite{dragos}, Lemma 2.10] For $r=1$ in the previous lemma,
\[\pi^{\mathrm{vec}}(\mathcal{K}_{x}^{}\text{ } \E) = [\mathcal{K}_{x}^{}, \E],\]
where the commutator is taken in the Hall algebra $\mathbf{H}_Y.$ 
\end{lemma}

We conclude this section with two applications of our results in \cite{roberto-graphs} to $\mathbf{H}_Y.$ The first one is the following corollary and the second one is stated in the next remark.

\begin{corollary} $\displaystyle\sum_{\E ' \in \mathrm{Bun}_n X} h_{\mathcal{K}_{x}^{r},\E '}^{\E} = \dfrac{(q_{x}^{n} -1)(q_{x}^{n-1} - 1)(q_{x}^{n-2} - 1) \cdots (q_{x}^{r+1} - 1)}{(q_{x}^{n-r} -1)(q_{x}^{n-r-1} - 1)(q_{x}^{n-r-2}-1) \cdots (q_{x} -1)}.$ 

\begin{proof} This follows from  \cite[Thm. 2.7]{roberto-graphs} where we prove that for a fix $x \in |X|$  
$$\sum_{\E' \in \mathrm{Bun}_n Y} m_{x,r}(\E,\E') = \# \mathrm{Gr}(n-r,n)(\kappa(x)),$$
and the Lemma \ref{lemmaconnection}. 
\end{proof}
\end{corollary}

\begin{remark} The products $\pi^{\mathrm{vec}}(\mathcal{K}_{x}^{\oplus r} \E)$ are explicitly determined for every $\E \in \mathrm{Bun}_n \mathbb{P}_{\Fq}^{1},$ $1 \leq r \leq n$ and $x$ a closed point in $\mathbb{P}_{\Fq}^{1}.$
Namely, fixed an positive integer $n,$ in \cite{roberto-graphs} we develop an algorithm to describe $\mathcal{G}_{x,r}$ for every $x$ closed point in $\mathbb{P}_{\Fq}^{1}$ and $1 \leq r \leq n.$ By Remark \ref{propconnection} we are able to determine $\pi^{\mathrm{vec}}(\mathcal{K}_{x}^{\oplus r} \E)$.

\end{remark}

\section{Graphs of Hecke operators for elliptic curves} \label{section3}

\noindent
We return our attention to an elliptic curve $X$ defined over $\mathbb{F}_q$. As we have seen,  describing the graph $\mathcal{G}_{x,r}$ for $X$ is equivalent to calculating the products $\pi^{\mathrm{vec}} (\mathcal{K}_{x}^{r} \text{ } \E)$ in the Hall algebra of $X$  (Remark \ref{propconnection}) where $\E$ runs through $\mathrm{Bun}_n X.$ The strategy is to write these product in terms of elements in the \textit{twisted spherical Hall algebras} and to use the structure of these algebras to compute the products. 

We shall need the usual notions of $\nu$-integers. If $\nu \neq \pm 1$, then we set 
$$[s]_{\nu} := \frac{\nu^s - \nu^{-s}}{\nu - \nu^{-1}}.$$
We shall mostly use $[s]:=[s]_{v}$ where $v^2 = q^{-1}.$ 

Let us recall some properties of the classical Hall algebra, i.e.\ the Hall algebra of finite modules over a discrete valuation ring. For a deeper discussion of classical Hall algebras, we refer the reader to \cite{macdonald}, chapters I, II and III.  

A \textit{partition} is any finite sequence $\lambda = (\lambda_1, \ldots, \lambda_m)$ of non-negative integers in decreasing order: $\lambda_1 \geq \cdots \geq \lambda_m.$ The non-zero $\lambda_i$ are called the \textit{parts} of $\lambda.$ The number of parts is the \textit{length} of $\lambda,$ denote by $l(\lambda);$ and the sum of the parts $\sum_{i=1}^{m}\lambda_i$ is the \textit{weight} of $\lambda,$ denote by $|\lambda|.$ We denote by $\lambda = (1^{d_1}2^{d_2} \ldots m^{d_m})$ the partition that has exactly $d_i$ parts equal to $i.$  

For a finite field $\mathbf{k}$, we denote by $u$ a square root of $(\#\mathbf{k})^{-1}.$ Denote by $\mathcal{A}_{\mathbf{k}}$ the category of finite modules over the discrete valuation ring $\mathrm{R} := \mathbf{k}[[t]].$ For every $s \in \N$ there exists a unique indecomposable module of length $s$ (up to isomorphism), denoted by $I_{(s)}$, which is  the quotient $\mathrm{R}\big/ t^s \mathrm{R}.$ For a partition $\lambda = (\lambda_1, \ldots, \lambda_m)$, we write by $I_{\lambda} := I_{(\lambda_1)} \oplus \cdots \oplus I_{(\lambda_m)}.$ The collection $\{I_{\lambda}\}_{\lambda}$ where $\lambda$ runs over all partitions is a complete collection of representatives for the isomorphism classes of objects of $\mathcal{A}_{\mathbf{k}}.$ Let us denote by $\Lambda :=\Z[e_1,e_2, \ldots]$ the Macdonald's ring of symmetric function, where 
$$e_r = \sum_{i_1 < \cdots < i_r} x_{i_1} \cdots x_{i_r}$$ 
are the elementary symmetric functions in infinite many variables. 
We denote by $\Lambda_t$ the Macdonald's ring of symmetric functions over $\C[t^{\pm 1}]$ i.e.\ 
$$\Lambda_t = \C[t^{\pm 1}] \otimes \Z[e_1,e_2, \ldots].$$ 
Furthermore, let us denote by $p_{r}$ the power-sum symmetric function given by
\[p_r = \sum x_{i}^{r}.\]
For a partition $\lambda=(\lambda_1, \ldots, \lambda_m)$, we denote by $P_{\lambda}$ the Hall-Littlewood symmetric functions (see \cite{macdonald} III.2 for the definition). We write $p_{\lambda}:= p_{\lambda_1}  \cdots p_{\lambda_m}$ and $e_{\lambda} := e_{\lambda_1}  \cdots e_{\lambda_m}$. 

\begin{proposition}[\cite{macdonald}] \label{prop-macdonald}The assignment $I_{(1^m)} \mapsto u^{m(m-1)} e_m$ extends to an algebra isomorphism $\Psi_{\mathbf{k}}: \mathbf{H}_{\mathcal{A}_{\mathbf{k}}} \rightarrow \Lambda_{t}|_{t=u^2}.$ Moreover,
\begin{enumerate}
\item[(i)]  $\Psi_{\mathbf{k}}^{-1}(p_m) = \sum_{|\lambda|=m} n_u(l(\lambda)-1)I_{\lambda} \quad \text{ where }\quad  n_u(l)= \prod_{i=1}^{l}(1-u^{-2i}),$ and  \\
\item[(ii)]  $\Psi_{\mathbf{k}}(I_{(\lambda)}) = u^{2 n(\lambda)} P_{\lambda} \quad \text{ where }\quad n(\lambda)= \sum (i-1)\lambda_i.$
\end{enumerate}
Here, $l(\lambda)$ is the length of $\lambda$ and $|\lambda|$ its weight. 

\begin{proof} The proof of the isomorphism and statement (ii) can be found in \cite{macdonald} III.3 (3.4). The  statement (i) is \cite{macdonald}, III.7 Ex. 2. 
\end{proof}
\end{proposition}

Let $x$ be a closed point of $X$. Consider the category $\mathrm{Tor}_x$ of torsion sheaves on $X$ supported at $x$. We have an equivalence of categories $\mathrm{Tor}_x \cong \mathcal{A}_{\kappa(x)},$ which provides us with an isomorphism of algebras $\Psi_{\kappa(x)}: \mathbf{H}_{\mathrm{Tor}_x} \rightarrow \Lambda_{t}|_{t=v^{2|x|}},$ where $ v^2 = q^{-1}.$

For $m \in \N$ and $x \in |X|,$ we define the elements $T_{(0,m),x} \in \mathbf{H}_X$ by

\[T_{(0,m),x} := \begin{cases}
0  & \text{if} \quad |x| \; \ndiv  m \\
\frac{[m]}{m} |x| \Psi_{\kappa(x)}^{-1}( p_{\frac{m}{|x|}})  & \text{if} \quad  |x| \; | \; m.
 \end{cases} \]
Using Proposition \ref{prop-macdonald}, we can rewrite $T_{(0,m),x}$ explicitly as follows
\[T_{(0,m),x}= \frac{[m]|x|}{m} \sum_{|\lambda|=m/|x|} n_{u_x}(l(\lambda)-1)\mathcal{K}_{x}^{(\lambda)},\]
for $|x| \; | \; m.$

For any $\mu \in \Q \cup \{\infty \}$, we consider the subspace $\mathbf{H}_{X}^{(\mu)} \subset \mathbf{H}_X$ that is spanned by classes $\{[\mathcal{F}] | \F \in \mathsf{C}_{\mu}\}.$ Since the category $\mathsf{C}_{\mu}$ is stable under extensions, $\mathbf{H}_{X}^{(\mu)}$ is a subalgebra of $\mathbf{H}_X.$ The exact equivalence $\epsilon_{\mu,\nu}$ from Theorem \ref{atiyahtheorem} gives rise to an algebra isomorphism $\epsilon_{\mu,\nu} : \mathbf{H}_{X}^{(\nu)} \stackrel{\sim}{\rightarrow} \mathbf{H}_{X}^{(\mu)}.$ 

Note that $T_{(0,m),x} \in \mathbf{H}_{X}^{(\infty)}.$ For every $\mathbf{v} \in \mathbf{Z}$ with slope $\mu,$ we define the element $T_{\mathbf{v},x}$ by $\epsilon_{\mu,\infty}(T_{(0,m),x}).$ Namely, there is $f \in \mathrm{SL}_2(\Z)$ such that $f(0,m)=\mathbf{v}$ and thus $T_{\mathbf{v},x} := f \cdot T_{(0,m),x}$ (cf. paragraph \ref{atiyahparagraph}).

With the definitions of paragraph of Notations (Section \ref{section1}) in mind, we define:

\begin{definition} For a character $\widetilde{\rho} \in \widetilde{\mathrm{Pic}^{0}(X_n)}$ and a closed point $x \in X$, we define 
$$\widetilde{\rho}(x) := \frac{1}{n} \sum_{i=0}^{n-1}\rho((\mathrm{Fr}_{X,n})^{i}(\mathcal{O}_{X_n}(x'))) = \frac{1}{d} \sum_{i=0}^{d-1} \rho((\mathrm{Fr}_{X,n})^{i}(\mathcal{O}_{X_n}(x')))$$
where $x' \in X_n$ is a closed point that sits above $x$ and $d = \gcd(|x|,n).$
\end{definition}

\begin{definition} For a character $\widetilde{\rho} \in \widetilde{X}$ and for a point $\mathbf{v} \in \mathbf{Z}$ we define the element
$$T_{\mathbf{v}}^{\widetilde{\rho}} := \sum_{x \in |X|} \widetilde{\rho}(x) T_{\mathbf{v},x}.$$ 
\end{definition}

\begin{remark} Since the above definition is an average over all possible choices for the representative of $\widetilde{\rho}$  and the point $x'$, it is well defined. \end{remark}

\begin{proposition} \label{inversionlemma}  Let  $\mathbf{v}$ be a point in $\mathbf{Z}$ with $d := \gamma(\mathbf{v}).$ If $y \in |X|$ is such that $|y|\;|\;d$, then
$$T_{\mathbf{v},y} = |y| N_{d}^{-1} \sum_{\widetilde{\rho} \in \widetilde{\mathrm{Pic}^{0}(X_d)}} \widetilde{\rho}(-y) T_{\mathbf{v}}^{\widetilde{\rho}}.$$

\begin{proof} Recall that for a finite abelian group $G$,
\[\sum_{\rho \in \widehat{G}} \rho(g) = \begin{cases}  0 &\text{ if } g \neq 1_{G} \\
 \# G &\text{ if } g=1_{G}.
 \end{cases}
 \]
 By definition
$$T_{\mathbf{v}}^{\widetilde{\rho}} = \sum_{x \in |X|} \widetilde{\rho}(x) T_{\mathbf{v},x} = \sum_{\substack{ x \in |X| \\ |x|\;|\;d}} \widetilde{\rho}(x) T_{\mathbf{v},x}.$$
Summing up all orbits in $\widetilde{\mathrm{Pic}^{0}(X_d)},$ we obtain
\begin{eqnarray*} 
\sum_{\widetilde{\rho} \in \widetilde{\mathrm{Pic}^{0}(X_d)}} \widetilde{\rho}(-y) T_{\mathbf{v}}^{\widetilde{\rho}} &=& \sum_{\widetilde{\rho} \in \widetilde{\mathrm{Pic}^{0}(X_d)}} \sum_{\substack{ x \in |X| \\ |x|\;|\;d}} \widetilde{\rho}(x-y) T_{\mathbf{v},x} \\
& = &\sum_{\substack{ x \in |X| \\ |x|\;|\;d}} \sum_{\widetilde{\rho} \in \widetilde{\mathrm{Pic}^{0}(X_d)}} \frac{1}{|x|} \sum_{i=0}^{|x|-1} \rho((\mathrm{Fr}_{X,d})^{i}(\mathcal{O}_{X_d}(x'-y'-x_0)))T_{\mathbf{v},x} \\
& = &  \sum_{\substack{ x \in |X| \\ |x|\;|\;d}} \frac{1}{|x|} \sum_{\rho \in \widehat{\mathrm{Pic}^{0}(X_d)}} \rho(\mathcal{O}_{X_d}( x'-y'-x_0)) T_{\mathbf{v},x}\\
&=& |y|^{-1} N_d T_{\mathbf{v},y} 
\end{eqnarray*}
where $x',y' \in X_d$ with $x'$  above $x$ and $y'$ above $y.$
\end{proof}

\end{proposition}

In the following, we aim to write $\mathcal{K}_{x}^{\oplus r}$ and $\E$ (for $\E \in \mathrm{Bun}_n X$) in terms of the elements $T_{\mathbf{v}}^{\widetilde{\rho}}$ as the above definition. The next proposition tell us that this is possible. 

\begin{proposition}[\cite{dragos} Proposition 3.4] The set $\big\{ T_{\mathbf{v}}^{\widetilde{\rho}} \; \big| \; \mathbf{v} \in \mathbf{Z}, \widetilde{\rho} \in \widetilde{\mathrm{Pic}^{0}(X_d)}, d=\gamma(\mathbf{v}) \big\},$  generates the Hall algebra $\mathbf{H}_X.$ 

\end{proposition}

\begin{definition} We call a character $\rho \in \widehat{\mathrm{Pic}^0(X_n)}$ \textit{primitive} of degree $n$ if its orbit under the Frobenius $\mathrm{Fr}_{X,n}$ is of maximal cardinal, i.e.\ if it is of cardinal $n.$ Equivalently, $\rho$ is primitive if there does not exist a character $\chi \in \widehat{\mathrm{Pic}^0(X_m)}, m < n,\; m|n,$ such that $\rho = \mathrm{Norm}_{m}^{n}(\chi)$ (\cite{dragos} Lemma 3.8). 
\end{definition}

\begin{remark} Applying the above definition we see that $\rho \in \widehat{\mathrm{Pic}^0(X_n)}$ is primitive or there is a primitive character $\chi \in \widehat{\mathrm{Pic}^0(X_m)}$ with $m|n, \; m < n$ such that $\rho = \mathrm{Norm}_{m}^{n}(\chi)$ (c.f. \cite{dragos} Corollary 3.9). \end{remark}

Let $\mathcal{P}_n$ be the set of primitive characters of degree $n$ modulo the action of Frobenius and $\mathcal{P}:= \coprod_{n \geq 1}\mathcal{P}_n$.  

Next, we are able to define the twisted spherical Hall algebras. This subalgebras of $\mathbf{H}_X$  will play an important role in the algorithm (section \ref{sectionthealgorithm}), to calculate the graphs of Hecke operators. 

\begin{definition}\label{deftwistedsphericalalgebra} Let $n \geq 1$ and $\widetilde{\rho} \in \mathcal{P}_n$ be a primitive character. We define the algebra $\mathbf{U}_{X}^{\widetilde{\rho}},$ called the twisted spherical Hall algebra of $X$ and $\widetilde{\rho},$ as the subalgebra of $\mathbf{H}_X$ generated by  
\[\Big\{ T_{n\mathbf{v}}^{\mathrm{Norm}_{n}^{n \gamma(\mathbf{v})}(\widetilde{\rho})} \Big|\; \mathbf{v} \in \mathbf{Z} \Big\}.\]
\end{definition}

For $n=1$ and the trivial character $\widetilde{\rho}=1 \in \widetilde{\mathrm{Pic}^0(X)}$, the above definition  specializes to the spherical Hall algebra $\mathbf{U}_X$ introduced by Burban and Schiffmann in section 4 of \cite{olivier-elliptic1} and considered implicitly by Kapranov in  sections 3.8 and 5 of \cite{kapranov}.  
We also will need the combinatorial description of the twisted spherical Hall algebra as discovered in section 5 of \cite{olivier-elliptic1}  and generalized in section 3 of \cite{dragos}.

For $\mathbf{v},\mathbf{w} \in \mathbf{Z}$ which are not proportional we denote by $\epsilon_{\mathbf{v},\mathbf{w}}:= \mathrm{sign}(\det(\mathbf{v},\mathbf{w})) \in \{\pm 1\}$ and by $\Delta_{\mathbf{v},\mathbf{w}}$ the triangle formed by the vectors $\mathbf{o},\mathbf{v},\mathbf{v+w}$ where $\mathbf{o}=(0,0).$ 

We recall Pick's formula (\cite{pick}), which will be useful for us. For any pair of linear independents points $\mathbf{v},\mathbf{w} \in \Z^2,$
\begin{equation}
|\det(\mathbf{v},\mathbf{w})| = \gamma(\mathbf{v}) + \gamma(\mathbf{w}) + \gamma(\mathbf{v}+\mathbf{w}) -2 + 2\;\#(\Delta_{\mathbf{v},\mathbf{w}} \cap \Z^2).
\end{equation}

\begin{definition} \label{defpathalgebra}Fix $\sigma, \overline{\sigma} \in \C^{*}$ with $\sigma, \overline{\sigma} \not\in \{\pm 1\}$ and set $\nu:= (\sigma\overline{\sigma})^{-\frac{1}{2}}$ and 
\[ c_i(\sigma,\overline{\sigma}) := (\sigma^{i/2} - \sigma^{-i/2})(\overline{\sigma}^{i/2}-\overline{\sigma}^{-i/2})[i]_{\nu}\big/i \]
Let $\textbf{E}_{\sigma, \overline{\sigma}}^{n}$ be the $\C$-algebra generated by $\{t_{\mathbf{v}} \; |\; \mathbf{v} \in \mathbf{Z}\}$ modulo the following relations:
\begin{enumerate}
\item If $\mathbf{o},\mathbf{v},\mathbf{v'}$ are collinear then $[t_{\mathbf{v}},t_{\mathbf{v'}}]=0.$ \\
\item If $\mathbf{v},\mathbf{w}$ are such that $\gamma(\mathbf{v})=1$ and that $\Delta_{\mathbf{v},\mathbf{w}}$ has no interior lattice point, then
\[[t_{\mathbf{w}}, t_\mathbf{v}] = \epsilon_{\mathbf{v},\mathbf{w}} c_{n \gamma(\mathbf{w})}(\sigma, \overline{\sigma}) \frac{\theta_{\mathbf{v}+\mathbf{w}}}{n(\nu^{-1} - \nu)}\]
where the elements $\theta_{\mathbf{z}}, \mathbf{z} \in \mathbf{Z}$ are defined by the following generating series:
\[\sum_{i \geq 0}\theta_{i\mathbf{z}_0} s^i = \exp\left( n(\nu^{-1}-\nu) \sum_{i \geq 1} t_{i \mathbf{z}_{0}} s^i \right)\]
for any $\mathbf{z}_0 \in \mathbf{Z}$ such that $\gamma(\mathbf{z_0})=1.$
 \end{enumerate}
Note that $\theta_{\mathbf{z}}=n(\nu^{-1}-\nu)t_{\mathbf{z}}$ whenever $\gamma(\mathbf{z})=1.$ 
\end{definition}

Following section 4 of \cite{olivier-elliptic1},  we give a geometric interpretation for the algebra $\mathbf{E}_{\sigma, \overline{\sigma}}^{n}.$ By a path in $\mathbf{Z}$ we understand a sequence $\mathbf{p}=(\mathbf{v}_1, \ldots, \mathbf{v}_m)$ of non-zero elements of $\mathbf{Z}$, which we represent graphically as the polygonal line in $\mathbf{Z}$ that joins the points $\mathbf{o},\mathbf{v}_1,\mathbf{v}_1 + \mathbf{v}_2, \ldots, \mathbf{v}_1+ \cdots + \mathbf{v}_m.$  Let $\widehat{\mathbf{v}\mathbf{w}} \in [0,2\pi)$ denote the angle between the segments $\mathbf{o}\mathbf{v}$ and $\mathbf{o}\mathbf{w}.$ We will call a path $\mathbf{p}=(\mathbf{v}_1, \ldots, \mathbf{v}_m)$ convex if $\widehat{\mathbf{v}_1\mathbf{v}_2} \leq \widehat{\mathbf{v}_1\mathbf{v}_3} \leq \cdots \leq \widehat{\mathbf{v}_1\mathbf{v}_m} < 2\pi.$ Put $L_0 := (0,-1)$ and let $\mathbf{Conv'}$ be the collection of all convex paths $\mathbf{p}=(\mathbf{v}_1, \ldots, \mathbf{v}_m)$ satisfying $\widehat{\mathbf{v}_1 L_0} \geq \cdots \geq \widehat{\mathbf{v}_m L_0}.$ Two convex paths $\mathbf{p}=(\mathbf{x}_1, \ldots, \mathbf{x}_r)$ and $\mathbf{q}=(\mathbf{y}_1, \ldots, \mathbf{y}_s)$ in $\mathbf{Conv'}$ are said to be \textit{equivalent} if $\{\mathbf{x}_1, \ldots, \mathbf{x}_r\}=\{\mathbf{y}_1, \ldots, \mathbf{y}_s\}$ i.e.\ $\mathbf{p}$ is obtained by permuting several segments of $\mathbf{q}$ of the same slope. For example the path $\mathbf{p}=((2,-2),(1,0),(2,0),(1,1))$ is equivalent to the path $\mathbf{q}=((2,-2),(2,0),(1,0),(1,1)).$ 

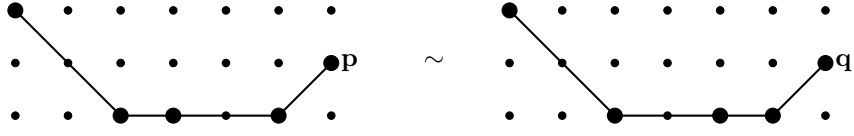
\begin{center}

\begin{figure}[h]
\centering

\begin{minipage}[b]{0.45\linewidth}
\[
 \beginpgfgraphicnamed{tikz/figequiv1}
  \begin{tikzpicture}[>=latex, scale=0.7]

   \foreach \x in {-1, -0.5 ,0,...,2}{
      \foreach \y in {0,0.5,1}{
        \node[draw,circle,inner sep=1pt,fill] at (2*\x,2*\y) {};
      }
    }
\node[draw,circle,inner sep=2pt,fill] at (-2,2) {}; 
\node[draw,circle,inner sep=2pt,fill] at (-0,0) {};
\node[draw,circle,inner sep=2pt,fill] at (1,0) {};
\node[draw,circle,inner sep=2pt,fill] at (3,0) {};
\node[draw,circle,inner sep=2pt,fill] at (4,1) {};       
\draw [ thick] (-2,2) -- (0,0) -- (1,0) -- (3,0)-- (4,1);
\draw [ thick] (4,1) node[right]{$\mathbf{p}$};
   \end{tikzpicture}
 \endpgfgraphicnamed
\]
\end{minipage} 
\begin{minipage}[b]{0.05\linewidth}
$$\sim$$
\vspace*{0.20cm}
\end{minipage}
\begin{minipage}[b]{0.45\linewidth}
\[
 \beginpgfgraphicnamed{tikz/figequiv1}
  \begin{tikzpicture}[>=latex, scale=0.7]

   \foreach \x in {-1, -0.5 ,0,...,2}{
      \foreach \y in {0,0.5,1}{
        \node[draw,circle,inner sep=1pt,fill] at (2*\x,2*\y) {};
      }
    }
\node[draw,circle,inner sep=2pt,fill] at (-2,2) {}; 
\node[draw,circle,inner sep=2pt,fill] at (-0,0) {};
\node[draw,circle,inner sep=2pt,fill] at (2,0) {};
\node[draw,circle,inner sep=2pt,fill] at (3,0) {};
\node[draw,circle,inner sep=2pt,fill] at (4,1) {};       
\draw [ thick] (-2,2) -- (0,0) -- (1,0) -- (3,0)-- (4,1);
\draw [ thick] (4,1) node[right]{$\mathbf{q}$};
   \end{tikzpicture}
 \endpgfgraphicnamed
\] 
\end{minipage}
\end{figure}

\end{center}

We denote by $\mathbf{Conv}$ the set of equivalence classes of paths in $\mathbf{Conv'}$ and we will call the elements of $\mathbf{Conv}$ simply paths. We introduce the positive paths $\mathbf{Conv}^{+}$ and the negatives paths $\mathbf{Conv}^{-}$ as the paths $\mathbf{p}=(\mathbf{v}_1, \ldots, \mathbf{v}_m) \in \mathbf{Conv}$ such that $\widehat{\mathbf{v}_m L_0} \geq \pi$ and   $\widehat{\mathbf{v}_1 L_0} < \pi,$ respectively. By concatenating paths we obtain an identification $\mathbf{Conv} \equiv \mathbf{Conv}^{+} \times \mathbf{Conv}^{-}.$

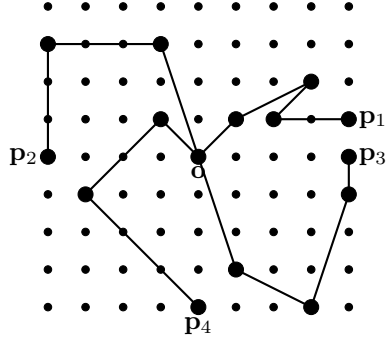
\[
 \beginpgfgraphicnamed{tikz/figconv}
  \begin{tikzpicture}[>=latex, scale=0.5]

\node[draw,circle,inner sep=2pt,fill] at (2,2) {};
\node[draw,circle,inner sep=2pt,fill] at (3,3) {};
\node[draw,circle,inner sep=2pt,fill] at (5,4) {};
\node[draw,circle,inner sep=2pt,fill] at (4,3) {};
\node[draw,circle,inner sep=2pt,fill] at (6,3) {};
\node[draw,circle,inner sep=2pt,fill] at (1,5) {};
\node[draw,circle,inner sep=2pt,fill] at (-2,5) {};
\node[draw,circle,inner sep=2pt,fill] at (-2,2) {};
\node[draw,circle,inner sep=2pt,fill] at (3,-1) {};
\node[draw,circle,inner sep=2pt,fill] at (5,-2) {};
\node[draw,circle,inner sep=2pt,fill] at (6,1) {};
\node[draw,circle,inner sep=2pt,fill] at (6,2) {};
\node[draw,circle,inner sep=2pt,fill] at (1,3) {};
\node[draw,circle,inner sep=2pt,fill] at (-1,1) {};
\node[draw,circle,inner sep=2pt,fill] at (2,-2) {};
    
   \foreach \x in {-1, -0.5 ,0,...,3}{
      \foreach \y in {-1, -0.5,0,...,3}{
        \node[draw,circle,inner sep=1pt,fill] at (2*\x,2*\y) {};
      }
    }
\draw [ thick] (2,2) -- (3,3) -- (5,4) -- (4,3)-- (6,3);
\draw [ thick] (2,2) -- (1,5) -- (-2,5) -- (-2,2);
\draw [ thick] (2,2) -- (3,-1) -- (5,-2) -- (6,1) -- (6,2);
\draw [ thick] (2,2) -- (1,3) -- (-1,1) --(2,-2);
\draw [ thick] (2,-2) node[below]{$\mathbf{p}_4$};
\draw [ thick] (6,2) node[right]{$\mathbf{p}_3$};
\draw [ thick] (-2,2) node[left]{$\mathbf{p}_2$};
\draw [ thick] (6,3) node[right]{$\mathbf{p}_1$};
\draw [ thick] (2,2) node[below]{$\mathbf{o}$};
\node [right, align=flush left,text width=5cm] at (-2.5,-4.5)
        { $\mathbf{p}_1$ is not convex \\ $\mathbf{p}_2 \in \mathbf{Conv}^{-}, \; \mathbf{p}_3 \in \mathbf{Conv}^{+}$ \\   $\mathbf{p}_4$ is convex, but $\mathbf{p}_4 \not\in \mathbf{Conv}$ 
           
        };
   \end{tikzpicture}   
 \endpgfgraphicnamed 
\]

Fix an integer $n \geq 1.$ With a path $\mathbf{p}=(\mathbf{v}_1, \ldots, \mathbf{v}_m) \in \mathbf{Conv},$ we associate the element $t_{\mathbf{p}} \in \mathbf{E}_{\sigma, \overline{\sigma}}^{n}$ defined by
\[ t_{\mathbf{p}} := t_{\mathbf{v}_1} \cdots t_{\mathbf{v}_m},\]
which is a well defined element of $ \mathbf{E}_{\sigma, \overline{\sigma}}^{n}$ due to relation (1).

\begin{lemma}[\cite{olivier-elliptic1} Lemma 5.6] \label{lemmaconvex}The algebra $\mathbf{E}_{\sigma, \overline{\sigma}}^{n}$ is equal to $\bigoplus_{\mathbf{p} \in \mathbf{Conv}^{+}} \C t_{\mathbf{p}}.$

\end{lemma}

Next we associate with each coherent sheaf a convex path in $\mathbf{Z}.$ By Theorem \ref{atiyahtheorem}, the HN-filtration splits, hence for every vector bundle $\E \in \mathrm{Bun}_n X$ we can write
\[\E = \E_1 \oplus \cdots \oplus \E_m\]
where $\E_i \in \mathsf{C}_{\mu_i},$ $\overline{\E_i} = (n_i,d_i)$ and $\mu_1 < \cdots < \mu_m.$ We associate with $\E$ the convex path $\mathbf{p}(\E)=(\mathbf{v}_1, \ldots, \mathbf{v}_m) \in \mathbf{Conv}^{+}$ where $\mathbf{v}_i := (n_i,d_i).$

Let $\mathbf{P}(x,r,\E)$ the polygon determined by the region below the convex path given by $\E$ and above the path given by $((0,-r|x|), \mathbf{v}_1, \ldots, \mathbf{v}_m).$ 

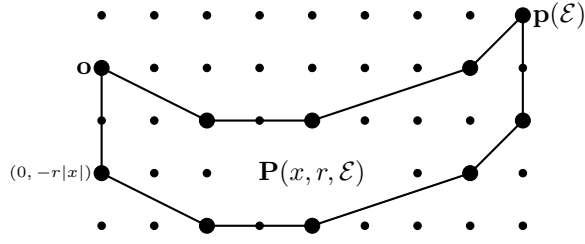
\[
 \beginpgfgraphicnamed{tikz/fig2}
  \begin{tikzpicture}[>=latex, scale=0.7]

   \foreach \x in {-1, -0.5 ,0,...,3}{
      \foreach \y in {-1, -0.5,0,0.5,1}{
        \node[draw,circle,inner sep=1pt,fill] at (2*\x,2*\y) {};
      }
    }
\draw [ thick] (-2,1) -- (0,0) -- (2,0) -- (5,1)-- (6,2);
\draw [ thick] (-2,-1) -- (0,-2) -- (2,-2) -- (5,-1)-- (6,0);
\draw [ thick] (-2,1) -- (-2,-1);
\draw [ thick] (6,0) -- (6,2);
\draw [ thick] (6,2) node[right]{$\mathbf{p}(\E)$};
\draw [ thick] (-2,-1) node[left]{{\tiny $(0,-r|x|)$}};
\draw [ thick] (-2,1) node[left]{$\mathbf{o}$};
\draw [ thick] (2,-1) node[ outer sep=2pt,fill=white]{$\mathbf{P}(x,r,\E)$};
\node[draw,circle,inner sep=2pt,fill] at (-2,1) {};
\node[draw,circle,inner sep=2pt,fill] at (0,0) {};
\node[draw,circle,inner sep=2pt,fill] at (2,0) {};
\node[draw,circle,inner sep=2pt,fill] at (5,1) {};
\node[draw,circle,inner sep=2pt,fill] at (6,2) {};
\node[draw,circle,inner sep=2pt,fill] at (-2,-1) {};
\node[draw,circle,inner sep=2pt,fill] at (0,-2) {};
\node[draw,circle,inner sep=2pt,fill] at (2,-2) {};
\node[draw,circle,inner sep=2pt,fill] at (5,-1) {};
\node[draw,circle,inner sep=2pt,fill] at (6,0) {};
   \end{tikzpicture}
 \endpgfgraphicnamed
\]

\begin{theorem}\label{propconvexpaths} Let $\E,\E' \in \mathrm{Bun}_n X$ and $x \in |X|.$ Suppose that $m_{x,r}(\E,\E') \neq 0.$ Then $\mathbf{p}(\E')$ is contained in the polygon $\mathbf{P}(x,r,\E)$.

\begin{proof} Let $\E,\E'$ be rank $n$ vector bundles on $X$ such that
$$0 \longrightarrow \E' \longrightarrow \E \longrightarrow \mathcal{K}_{x}^{\oplus r} \longrightarrow 0,$$
i.e.\ $m_{x,r}(\E,\E') \neq 0.$ Denote by $\E = \E_1 \oplus \cdots \oplus \E_s$, where $\overline{\E_i}:= \mathbf{v}_i$ and $\mu(\mathbf{v}_1) < \cdots < \mu(\mathbf{v}_s)$, the HN-decomposition of $\E$. Let $\E' = \E_1' \oplus \cdots \oplus \E_t'$, where $\overline{\E_i'}:= \mathbf{v}_i'$ and $\mu(\mathbf{v}_1') < \cdots < \mu(\mathbf{v}_t')$, be HN-decomposition of $\E'$.

Observe that the statement of the theorem is equivalent to say that $\mathbf{p}(\E)$ is contained in the polygon delimited by  $\mathbf{p}(\E')$ and $((0,r|x|), \mathbf{v}_1', \ldots, \mathbf{v}_t').$

First suppose that $\E'$ is semistable, i.e.\ $\E'=\E_1'$. Since $\E$ must appears in the product $\mathcal{K}_{x}^{\oplus r}\E'$ in the Hall algebra $\mathbf{H}_X$, it follows from Lemma \ref{lemma4.6} that $\mathbf{p}(\E)$ is a convex path in the triangle $\Delta_{(\mathbf{v}_1',(0,r|x|))}$ whose vertices are the origin $\mathbf{o}$, $\mathbf{v}_1'$ and $(0,r|x|)$. 
   
The general case follows by concatenation of  paths. Namely, we may write $\E'$ as the product  $v^{\sum_{i<j}\langle \E_i',\E_j' \rangle} \E_1' \cdots \E_t'$ in $\mathbf{H}_X$, thus $\E$ must appears in the product
\[ v^{\sum_{i<j}\langle \E_i',\E_j' \rangle} \mathcal{K}_{x}^{\oplus r} \E_1' \cdots \E_t'.\]
Applying the previous step inductively, we have the desired. 
\end{proof}
\end{theorem}

In the particular case of $n=2$ and $\E = \Line_1 \oplus \Line_2$ with $\deg \Line_2 - \deg\Line_1 > |x|$, we obtain a complete list of possible neighbours of $\E$ in $\mathcal{G}_{x,1},$ for every $x \in |X|.$  

\begin{corollary} Fix $x \in |X|$ and $n=2.$ Let $\Line_1, \Line_2$ be line bundles on X with $\deg(\Line_2) > \deg(\Line_1) +|x|$. Let us denote $\E = \Line_1 \oplus \Line_2,$ then 
$$\mathcal{V}_{x,1}(\E) = \big\{ (\E, \Line_1(-x) \oplus  
\Line_2,\; m_1), (\E,\Line_1 \oplus \Line_2(-x), \;m_2) \big\}$$
where $m_1 = m_{x,1}\big(\E, \Line_1(-x) \oplus  \Line_2\big)$ and $m_2 = \big(\E, \Line_1 \oplus  \Line_2(-x)\big).$

\begin{proof} By the last theorem there are only decomposable vector bundles on $\mathcal{V}_{x,1}(\E).$ All decomposable vector bundles on $\mathcal{V}_{x,1}(\E)$ have this form, c.f. Proposition 3.5  in \cite{roberto-graphs}.  
\end{proof}

\end{corollary}

Later (Proposition \ref{neighbousoftotaldec}) we will be able to show that $m_1=1$ and $m_2=q_x$ in the above corollary. This recovers the special case for $n=2$ in \cite{oliver-elliptic}.

Next we state the Theorem 5.2 of \cite{dragos}, which yields the relationship between the algebras $\mathbf{U}_{X}^{\widetilde{\rho}}$ and $\mathbf{E}_{\sigma, \overline{\sigma}}^{k}$ defined before. Recall that $N_i$ denotes the number of rational points of $X$ over $\mathbb{F}_{q^i}$ and that $v=q^{-1/2}.$ By the Hasse-Weil theorem (see e.g. \cite{biblia}, Appendix C), there exist conjugate algebraic numbers $\sigma, \overline{\sigma}$, satisfying $\sigma \overline{\sigma} =q$ such that 
\[N_i = q^i + 1 - (\sigma^i + \overline{\sigma}^{i})\]
for any $i \geq 1.$  These numbers $\sigma, \overline{\sigma}$ are the eigenvalues of the Frobenius automorphism acting on $H^1(\overline{X}, \overline{\Q_l}).$ Note that 
$$c_{i}(\sigma, \overline{\sigma}) = \frac{v^i [i] \#X(\mathbb{F}_{q^i})}{i},$$
cf. Definition \ref{defpathalgebra}. For $\sigma, \overline{\sigma}$ fixed as before, we denote  $c_{i}(\sigma, \overline{\sigma})$ simply by $c_i.$

\begin{theorem}[\cite{dragos} Theorem 5.2] \label{dragosmaintheorem}Let $\sigma, \overline{\sigma}$ as in the Hasse-Weil theorem. For any two integers $n,m \geq 1$ and any two different primitive characters $\widetilde{\rho} \in \mathcal{P}_n$ and $\widetilde{\sigma} \in \mathcal{P}_m$, we have
\begin{enumerate}

\item The twisted spherical Hall algebra $\mathbf{U}_{X}^{\widetilde{\rho}}$ is isomorphic to $\mathbf{E}_{\sigma,\overline{\sigma}}^{n}.$ This isomorphism is given by
\[t_{\mathbf{v}} \longmapsto T_{n \mathbf{v}}^{\mathrm{Norm}_{n}^{n \gamma(\mathbf{v})}(\widetilde{\rho})}.\]

\item The algebras $\mathbf{U}_{X}^{\widetilde{\rho}}$ and $\mathbf{U}_{X}^{\widetilde{\sigma}}$ commute with each other. 

\item The Hall algebra $\mathbf{H}_X$ decomposes into a commutative restricted tensor product
\[\mathbf{H}_X \cong \bigotimes_{\widetilde{\rho} \in \mathcal{P}}\text{$'$ } \mathbf{U}_{X}^{\widetilde{\rho}} \]
\end{enumerate}

\end{theorem}


\noindent
\textbf{Notation.}\label{remarknotation} By Atiyah's classification of vector bundles over $X$, see Theorem \ref{atiyahtheorem}, an irreducible vector bundle $\E$ is semistable. Let $\mu$ be its slope. Then its image via $\epsilon_{\infty,\mu}$ is an irreducible torsion sheaf on  $X.$ Hence $\E$ is completely determined by (i) its image $\overline{\E}=(n,d)$ in the numerical Grothendieck group $\Z^2$; (ii) a closed point $x \in |X|$ which is the support of $\epsilon_{\infty,\mu}(\E)$; and (iii) a weight $\ell$ that determines the unique irreducible torsion sheaf with support in $x$ as $\mathcal{K}_{x}^{(\ell)}.$  
We denote $\E$  by $\E_{(x,\ell)}^{(n,d)}.$\\

\section{The algorithm}\label{sectionthealgorithm}\label{sectionalgorithm}

\noindent
In the following,  we explain the main theorems that provide the basis for our algorithm to calculate the structure constants $h_{\mathcal{K}_{x}^{\oplus r}\E'}^{\E}.$

\begin{theorem}[\cite{olivier-elliptic1}]\label{th1} Let $\F = \F_1 \oplus \cdots \oplus \F_s$ be a coherent sheaf on $X$.  If $\F_i \in \mathsf{C}_{\mu_i}$ for $i=1, \ldots,n,$ and $\mu_1 < \cdots < \mu_s,$ then 
\[[\F] = v^{\sum_{i<j}\langle \F_i, \F_j \rangle} [\F_1] \cdots [\F_s]\]
in the Hall algebra of $X.$

\begin{proof} This theorem follows from Lemma 2.4 and (2.6) in \cite{olivier-elliptic1}.  
\end{proof}
\end{theorem}

\begin{theorem} \label{th2} If $\F \in \mathsf{C}_{\mu},$ then 
\[\F = \prod_{i=1}^{n} \sum_{j=1}^{n_i} \sum_{\widetilde{\rho}_{ij_1} } \cdots \sum_{\widetilde{\rho}_{ij_{k_j}} } a(\widetilde{\rho}_{ij_1}, \ldots, \widetilde{\rho}_{ij_{k_j}}) \;T_{\mathbf{v}_{ij_1}}^{\widetilde{\rho}_{ij_1}} \cdots T_{\mathbf{v}_{ij_{k_j}}}^{\widetilde{\rho}_{ij_{k_j}}}. \]
for some $a(\widetilde{\rho}_{ij_1}, \ldots, \widetilde{\rho}_{ij_{k_j}}) \in \C$  where $\widetilde{\rho}_{ij_k}$ runs through $ \widetilde{\mathrm{Pic}^{0}(X_{\gamma(\mathbf{v}_{ij_k})})}$ for $i=1, \ldots,n$ and $j=1, \ldots, n_i.$
Moreover, $\mu(\mathbf{v}_{ij_k})=\mu$ and $\sum_{i=1}^{n} \sum_{k=1}^{k_j} \mathbf{v}_{ij_{k}} = \overline{\F}.$

\begin{proof} By Theorem \ref{atiyahtheorem}, $\F$ corresponds to a torsion sheaf $\mathcal{T}$. Since $\mathrm{Tor}(X)$ is the direct product of blocks $\prod_{x \in |X|} \mathrm{Tor}_x$ where $\mathrm{Tor}_x$ is the category of torsion sheaves supported at $x,$ we may write 
$$\mathcal{T} = \mathcal{K}_{x_1}^{(\lambda_1)} \oplus \cdots \oplus \mathcal{K}_{x_n}^{(\lambda_n)},$$
where $\lambda_1, \ldots, \lambda_n$ are partitions and $x_i \neq x_j$ if $i \neq j$. Since the $\mathrm{Ext}$ of two torsion sheaves with disjoint support is trivial, we have 
\[\mathcal{T} =   \mathcal{K}_{x_1}^{(\lambda_1)}  \cdots  \mathcal{K}_{x_n}^{(\lambda_n)} \]
in the Hall algebra of $X.$ Via Proposition \ref{prop-macdonald}, $\mathcal{K}_{x_i}^{(\lambda_i)}$ corresponds to the Hall-Littlewood symmetric function $P_{\lambda_i}(\underline{x},q_{x_i}^{-1}).$ From \cite{macdonald}, I.2.12, the power-sums generate the Macdonald's ring of symmetric functions (as $\C[q_{x_i}^{\pm1}]$-algebra). Thus we can write
$$P_{\lambda_i}(\underline{x},q_{x_i}^{-1}) = \sum_{j=1}^{n_i} a_{ij}\; p_{\lambda_{ij}}$$
for certain $a_{ij} \in \C$ and partitions $\lambda_{ij}$  with $|\lambda_{ij}|=|\lambda_i|$. Hence 
$$\mathcal{K}_{x_i}^{(\lambda_i)} = \sum_{j=1}^{n_i} b_{ij}\; T_{(0,m_{ij_1}),x_i} \cdots T_{(0,m_{ij_{k_j}}),x_i}$$
where $\sum_{k=1}^{k_j} m_{ij_k} = |x_i||\lambda_i|,$ for $j=1, \ldots, n_i.$ Thus
$$\mathcal{T} = \prod_{i=1}^{n} \sum_{j=1}^{n_i} b_{ij} \; T_{(0,m_{ij_1}),x_i} \cdots T_{(0,m_{ij_{k_j}}),x_i},$$
with $b_{ij} \in \C$. Using Atiyah's theorem again, we get
\[\F = \prod_{i=1}^{n} \sum_{j=1}^{n_i} b_{ij} \; T_{\mathbf{v}_{ij_{1}},x_i} \cdots T_{\mathbf{v}_{ij_{k_j}},x_i} \]
where $\gamma(\mathbf{v}_{ij_{k}}) = m_{ij_{k}},$ $\mu(\mathbf{v}_{ij_{k}})=\mu$ and
$\sum_{i=1}^{n} \sum_{k=1}^{k_j} \mathbf{v}_{ij_{k}} = \overline{\F}.$
The desired follows from Proposition \ref{inversionlemma}. 
\end{proof}
\end{theorem}

The next result is a particular case of the previous theorem, which will be of later use.  

\begin{corollary} \label{co2} Let $x \in |X|$ and $r \geq 1$ an integer. Then
\[\mathcal{K}_{x}^{\oplus r} = \sum_{i=1}^{n} \sum_{\widetilde{\rho}_{i_1}} \cdots \sum_{\widetilde{\rho}_{i_k} } b_{ij}\; T_{(0,m_{i_1})}^{\widetilde{\rho}_{i j_1}}  \cdots T_{(0,m_{i_k})}^{\widetilde{\rho}_{i j_k} }\]
for some $ b_{ij} \in \C$ where $\widetilde{\rho}_{i_k}$ runs through the elements in $\widetilde{\mathrm{Pic}^{0}(X_{m_{i_k}})}$ for $i=1, \ldots,n.$ 

\begin{proof} Since $\mathcal{K}_{x}^{\oplus r} \in \mathsf{C}_{\infty},$ this follows from previous theorem.  
\end{proof}
\end{corollary}

\begin{example} Let $x$ be a degree $d$ closed point in $X.$ By definition
$ T_{(0,d),x} = [d] \mathcal{K}_{x}^{}$  and thus $\mathcal{K}_{x}^{} = [d]^{-1} T_{(0,d),x}.$ From Proposition \ref{inversionlemma}, 
\[T_{(0,d),x}= d N_{d}^{-1} \sum_{\widetilde{\rho} \in \mathrm{P}_d} \widetilde{\rho}(-x) T_{(0,d)}^{\widetilde{\rho}}.\]
Therefore, $\mathcal{K}_{x}^{} = [d]^{-1} d N_{d}^{-1} \sum_{\widetilde{\rho} \in \mathrm{P}_d} \widetilde{\rho}(-x) T_{(0,d)}^{\widetilde{\rho}}. $
\end{example}

\begin{remark}  Via Proposition \ref{prop-macdonald}, $\mathcal{K}_{x}^{\oplus r}$ corresponds to the elementary symmetric functions $e_r$  and these functions have an explicitly description in terms of power-sums given by Newton's formula, see \cite{macdonald} I.(2.11'). These formulas are useful for the calculation of our graphs.  We state these formulas in the following Step-4 of our algorithm.  
\end{remark}

\begin{theorem} \label{th3} Let $\E \in \mathrm{Bun}_n X,$ then
\[\E =  \sum_{\widetilde{\rho}_{ij_k}}^{} a_{ij_k} \;T_{\mathbf{v}_{1j_1}}^{\widetilde{\rho}_{1j_1}} \cdots T_{\mathbf{v}_{1j_{k_j}}}^{\widetilde{\rho}_{1 j_{k_j}}} \cdots T_{\mathbf{v}_{sj_1}}^{\widetilde{\rho}_{s j_1}} \cdots T_{\mathbf{v}_{s j_{k_j}}}^{\widetilde{\rho}_{s j_{k_j}}},\]
for some $a_{ij_k} \in \C$ where $\widetilde{\rho}_{ij_k}$ runs through the elements in $\widetilde{\mathrm{Pic}^{0}(X_{\gamma(\mathbf{v}_{ij_k})})}$ and the path given by $ \big( \mathbf{v}_{1j_1}, \ldots, \mathbf{v}_{1j_{k_j}}, \ldots, \mathbf{v}_{sj_1}, \ldots, \mathbf{v}_{sj_{k_j}}  \big)$ is a convex path which defines the same polygonal line as $\mathbf{p}(\E).$  

\begin{proof} Since the Harder-Narasinham filtration splits (c.f. Theorem \ref{atiyahtheorem}), we may write 
$\E = \E_1 \oplus \cdots \oplus \E_s$ with $\E_i \in \mathsf{C}_{\mu_i}$ for $i=1, \ldots,s,$ and 
$\mu_1 < \cdots < \mu_s.$ By Theorem \ref{th1},
\[ \E = v^{\sum_{i<j}\langle \E_i, \E_j \rangle} \E_1 \cdots \E_s\]
in $\mathbf{H}_X.$  The claim follows from Theorem \ref{th2}. 
\end{proof}
\end{theorem}

Recall from section \ref{section3} that any $\mu \in \Q \cup \{\infty\}$ yields the subalgebra $\mathbf{H}_{X}^{(\mu)} \subset \mathbf{H}_{X}$ of $\mathbf{H}_X$ linearly spanned by $\big\{ \F \big| \F \in \mathsf{C}_{\mu} \big\}.$ Moreover, the exact equivalence $\epsilon_{\mu_1, \mu_2}$ defined in Theorem \ref{atiyahtheorem} gives rise to an algebra isomorphism
$\epsilon_{\mu_1, \mu_2} : \mathbf{H}_{X}^{(\mu_2)} \rightarrow \mathbf{H}_{X}^{(\mu_1)}.$ For $\mu_1 \leq \mu_2$, let $\vec{\bigotimes}_{\mu_1 \leq \mu \leq \mu_2} \mathbf{H}_{X}^{(\mu)}$ stand for the (restricted) tensor product of spaces $\mathbf{H}_{X}^{(\mu)}$ with $\mu_1 \leq \mu \leq \mu_2$, ordered from left to right in increasing order, that is, for the vector spanned by elements of the form $a_{\nu_1} \otimes \cdots \otimes a_{\nu_r}$ with $a_{\nu_i} \in \mathbf{H}_{X}^{(\nu_i)}$ and $\mu_1 \leq \nu_1 < \cdots < \nu_r \leq \mu_2.$

Let $\mathsf{C}[\mu_1, \mu_2]$ be the full subcategory of sheaves whose HN-decomposition contains only slopes $\mu \in [\mu_1, \mu_2].$ This category is exact and, in particular, stable under extensions.  

\begin{lemma}[\cite{olivier-elliptic1} Lemma 2.5]\label{lemma4.6} For any $\mu_1 \leq \mu_2$ the Hall algebra of the exact category $\mathsf{C}[\mu_1, \mu_2]$ is a subalgebra of $\mathbf{H}_X$ isomorphic to  $\vec{\bigotimes}_{\mu_1 \leq \mu \leq \mu_2} \mathbf{H}_{X}^{(\mu)}$ (via the multiplication map).

\end{lemma}

\begin{theorem} \label{th4} Let $x \in |X|$, $r \geq 1$ an integer and $\E \in \mathrm{Bun}_n X$, then
\[ \mathcal{K}_{x}^{\oplus r} \; \E =  \sum_{i=1}^{m} a_{i}\; T_{\mathbf{v}_{i_1}}^{\widetilde{\rho}_{i_1}} \cdots T_{\mathbf{v}_{i_{\ell}}}^{\widetilde{\rho}_{i_{\ell}}}\]
for some $a_i \in \C$ where $\widetilde{\rho}_{i_j}$ runs through the elements in $\widetilde{\mathrm{Pic}^{0}(X_{\gamma(\mathbf{v}_{i_j})})},$  and where  $\mathbf{v}_{i_j} \in \mathbf{Z}$ are such that $ \big( \mathbf{v}_{i_1}, \ldots, \mathbf{v}_{i_{\ell}} \big)$ defines a convex path in $\Delta_{(\rk(\E),\deg(\E)),(0,r|x|)}$ for all $i=1, \ldots,m.$

\begin{proof} Considering the HN-decomposition $\E = \E_1 \oplus \cdots \oplus \E_s$ with $\E_i \in \mathsf{C}_{\mu_i}$ and $\mu_1 < \cdots < \mu_s,$ we obtain that $\mathcal{K}_{x}^{\oplus r}, \E \in \mathsf{C}[\mu_1, \infty].$ Since $\vec{\bigotimes}_{\mu_1 \leq \mu \leq \infty} \mathbf{H}_{X}^{(\mu)}$ is a subalgebra of $\mathbf{H}_X,$ the product $ \mathcal{K}_{x}^{\oplus r} \; \E$ still belong to $\vec{\bigotimes}_{\mu_1 \leq \mu \leq \infty} \mathbf{H}_{X}^{(\mu)}$ and thus we may write
\[  \mathcal{K}_{x}^{\oplus r} \; \E = \sum_{i=1}^{m'} a_{i}' \; a_{\mu_{i_1}} \cdots a_{\mu_{i_t}} \]
with $a_i' \in \C$, $a_{\mu_{i_j}} \in \mathbf{H}_{X}^{(\mu_{i_j})}$ and $\mu_1 \leq \mu_{i_1} < \cdots < \mu_{i_t} \leq \infty$ for all $i=1, \ldots, m'.$ Hence 
\[ \mathcal{K}_{x}^{\oplus r} \; \E = \sum_{i=1}^{m''}a_{i}'' \; \E_{\mu_{i_1}} \cdots \E_{\mu_{i_t}}\]
with $a_{i}'' \in \C$ and $\E_{\mu_{i_j}} \in \mathsf{C}_{\mu_{i_j}}.$ Since $\mu_1 \leq \mu_{i_1} < \cdots < \mu_{i_t} \leq \infty$ for all $i=1, \ldots, m''$ and $\E_{\mu_{i_1}} \cdots \E_{\mu_{i_t}} \in \mathrm{Ext}(\mathcal{K}_{x}^{\oplus r},\E),$ it follows from Lemma \ref{lemma4.6} that $\big(\overline{\E}_{\mu_{i_1}}, \ldots, \overline{\E}_{\mu_{i_t}} \big)$ defines a convex path in $\Delta_{(\rk(\E),\deg(\E)),(0,r|x|)}.$ 

Since each $\E_{\mu_{i_j}}$ is semistable we can use  Theorem \ref{th1} to write $\E_{\mu_{i_j}}$ in terms of the generators of the twisted spherical Hall algebras. The theorem follows from rearranging the terms.  
\end{proof}
\end{theorem}

\paragraph*{\textbf{The algorithm}} \label{thealgorithm} Fix $r,n \geq 1$ integers and $x \in |X|$. Let $\E \in \mathrm{Bun}_n X$. As explained in Remark \ref{propconnection}, it suffices for the determinate of the graph $\mathcal{G}_{x,r}$ to calculate the products $\pi^{\mathrm{vec}}\big(\mathcal{K}_{x}^{\oplus r} \; \E\big) = \pi^{\mathrm{vec}}[\mathcal{K}_{x}^{\oplus r}, \E]$ for every $\E \in \mathrm{Bun}_n X.$ Our aim is to write the product $\mathcal{K}_{x}^{\oplus r} \; \E$ as in  Theorem \ref{th4} and then use the definition of the elements in the twisted spherical Hall algebras and the structure of the Hall algebra $\mathbf{H}_X$ to calculate this product explicitly. In the following steps, we outline how to perform such calculations. \\

\paragraph*{\textbf{Step 1}} The Harder-Narashimhan decomposition tells us that every vector bundle $\E$ can be written as follows
\[\E = \E_1 \oplus \cdots \oplus \E_s \]
where each $\E_i \in \mathsf{C}_{\mu_i}$ for $i=1, \ldots,s,$ and $\mu_1 < \cdots < \mu_s.$  \\

\paragraph*{\textbf{Step 2}} We write each $\E_i$ in terms of elements in the twisted spherical Hall algebras. Namely, via Atiyah's theorem (\ref{atiyahtheorem}) $\E_i$ corresponds to a torsion sheaf of degree $\gcd(\rk(\E_i), \deg(\E_i)).$ Such a torsion sheaf is the direct sum of torsions sheaves with disjoint support. Let $\mathcal{T}_x$ be a torsion sheaf with support in $x.$ By Proposition \ref{prop-macdonald},  $\mathcal{T}_x$ corresponds to the Hall-Littlewood symmetric function $P_{\lambda}$ for some partition $\lambda$. Writing $P_{\lambda}$ as linear combination of $p_{\lambda_i}$ for some partitions $\lambda_1, \ldots, \lambda_m,$ we are able to write the $\E_i$ as sums of products of elements in the twisted spherical Hall algebras as stated in  Theorem \ref{th2}.   \\ 

\paragraph*{\textbf{Step 3}} By Theorem \ref{th1}, 
\[\E = v^{\sum_{i<j}\langle \E_i, \E_j \rangle} \E_1 \cdots \E_s. \] 
By the previous step, we may write $\E$ as sums of products of elements in the spherical Hall algebras, cf.\ Theorem \ref{th3}. \\

\paragraph*{\textbf{Step 4}} By Corollary \ref{co2} the previous step applies also to $\mathcal{K}_{x}^{\oplus r}.$ 
By the definition of $T_{\mathbf{v}}^{\widetilde{\rho}}$,  we only have to known how to write the elementary symmetric functions in terms of power-sums (cf. Proposition \ref{prop-macdonald}). These are the well-known Newton formulas: 
\[ m e_m = \sum_{i=1}^{m} (-1)^{i-1} p_i e_{m-i}. \] 
Therefore, we may write $\mathcal{K}_{x}^{\oplus r}$ as sums of products of elements in the twisted spherical Hall algebras with slope $\infty.$ \\

\paragraph*{\textbf{Step 5}} Since we know how to write $\E$ and $\mathcal{K}_{x}^{\oplus r}$ as elements in the twisted spherical Hall algebras, the product $\mathcal{K}_{x}^{\oplus r} \; \E$ can be written as sums of products of elements in the twisted spherical Hall algebras. By Corollary \ref{co2} and Theorem \ref{th3}, we obtain an expression of the form
\[ \mathcal{K}_{x}^{\oplus r} \; \E =  \sum_{i=1}^{m} a_{i}\; T_{(0,m_{i_1})}^{\widetilde{\rho}_{i_1}}  \cdots T_{(0,m_{i_k})}^{\widetilde{\rho}_{i_k} }  \;T_{\mathbf{v}_{i_1}}^{\widetilde{\rho}_{i_1}} \cdots T_{\mathbf{v}_{i_{\ell}}}^{\widetilde{\rho}_{i_{\ell}}}\] 
The problem here is that these products of elements in the spherical Hall algebras are not in increasing order of slopes. Since the algebras $\mathbf{U}_{X}^{\widetilde{\rho}}$ and $\mathbf{U}_{X}^{\widetilde{\sigma}}$ commute for different primitive characters, we are left with the problem to write these products in increasing order of slopes in a twisted spherical Hall algebra $\mathbf{U}_{X}^{\widetilde{\rho}}$ for a fixed $\widetilde{\rho} \in \mathcal{P}_n$. \\

\paragraph*{\textbf{Step 6}} In order to write  products in $\mathbf{U}_{X}^{\widetilde{\rho}}$ in increasing order of slopes, we use the isomorphisms from Theorem \ref{dragosmaintheorem} to reduce the problem to  calculations in $\mathbf{E}_{\sigma, \overline{\sigma}}^{n}.$ 
Lemma \ref{lemmaconvex} tells us that we can write any element in $\mathbf{E}_{\sigma, \overline{\sigma}}^{n}$ as a linear combination of elements $t_{\mathbf{p}}$ where $\mathbf{p}$ ranges through convex paths in $\mathbf{Conv}^{+}$. If $\mathbf{x}, \mathbf{y} \in \mathbf{Z}$ are such that there are not interior lattice points in $\Delta_{\mathbf{x}, \mathbf{y}}$, (as in the proof of Lemma \ref{lemmaconvex}) then either
\[\gamma(\mathbf{x}) = \gamma(\mathbf{y}) = \gamma(\mathbf{x}+\mathbf{y}) =2 \quad \text{or}\quad  \gamma(\mathbf{x})=1 \quad \text{ or }\quad  \gamma(\mathbf{y})=1.\]
In th first case, we may assume up to the $\mathrm{SL}_2(\Z)$-action that $\mathbf{x}=(0,2)$ and $\mathbf{y}=(2,0)$, applying the Definition \ref{defpathalgebra} item (2) yields 
\[\big[t_{\mathbf{x}},t_{\mathbf{y}}\big] = c \;t_{(1,1)}^{2} + c_2 \Big( \frac{c_2}{c_1} - 2 \Big) t_{(2,2)},\]
where $c \in \C.$ For the last two cases, the relation (2) in Definition \ref{defpathalgebra} directly yields that we can write the bracket $\big[t_{\mathbf{x}},t_{\mathbf{y}}\big]$ as linear combination of $t_{\mathbf{p}}$ where $\mathbf{p}$ runs through the convex paths in $\Delta_{\mathbf{x},\mathbf{y}}.$

On the other hand, if $\mathbf{Z} \cap \Delta_{\mathbf{x},\mathbf{y}}$ is not empty, we can always write the bracket $\big[t_{\mathbf{x}},t_{\mathbf{y}}\big]$ in terms of brackets $\big[t_{\mathbf{x}'},t_{\mathbf{y}'}\big]$ such that  $\Delta_{\mathbf{x}',\mathbf{y}'}$ has no interior lattice points. In order to do that, we can choose $\mathbf{z} \in  \Delta_{\mathbf{x},\mathbf{y}}$ such that $\Delta_{\mathbf{z},\mathbf{x}}$ has no interior lattice points and 
$ \deg(\mathbf{z}) = \deg(\mathbf{x-z})=1. $
Hence, using the relation $(2)$ in the Definition \ref{defpathalgebra},
\[  \big[ t_{\mathbf{z}}, t_{\mathbf{x-z}} \big] = c_1 \frac{\theta_{\mathbf{x}}}{v^{-1}-v} = c_1 t_{\mathbf{x}} +t \]
for some $t$ belonging to the subalgebra $\left\langle t_{\mathbf{x}_0}, \ldots, t_{(\deg(\mathbf{x})-1)\mathbf{x}_0} \right\rangle,$ where $\mathbf{x}_0 = \frac{\mathbf{x}}{\deg(\mathbf{x_0})}.$ Therefore,
\[ c_1\big[t_{\mathbf{x}}, t_{\mathbf{y}} \big] = \big[ t_{\mathbf{x-z}}, [t_{\mathbf{y}}, t_{\mathbf{z}}] \big] + \big[ t_{\mathbf{z}}, [t_{\mathbf{x-z}}, t_{\mathbf{y}}] \big] + \big[ t_{\mathbf{y}}, t \big]. \]
The \textit{Claim} in \cite[pag. 30] {olivier-elliptic1}asserts that this process ends in a finite number of steps.   

 Therefore,  we can always write $\big[t_{\mathbf{x}},t_{\mathbf{y}}\big]$ as linear combination of $t_{\mathbf{p}}$ where $\mathbf{p}$ runs through the set of convex paths in $\Delta_{\mathbf{x},\mathbf{y}}.$\\

\paragraph*{\textbf{Step 7}} By Lemma \ref{lemmacommutator}, $\pi^{\mathrm{vec}}\big(\mathcal{K}_{x}^{\oplus r} \E\big) = \pi^{\mathrm{vec}}\big[\mathcal{K}_{x}^{\oplus r}, \E\big]$. Aiming to use the relations of $\mathbf{E}_{\sigma, \overline{\sigma}}^{n}$ as in the previous step, we consider $\big[\mathcal{K}_{x}^{\oplus r}, \E\big]$. For $\mathbf{x}, \mathbf{x}_1, \ldots, \mathbf{x}_{\ell} \in \mathbf{Z}$ with $\mu(\mathbf{x})=\infty$ and $\mu(\mathbf{x}_1) \leq \cdots \leq \mu(\mathbf{x}_{\ell})$ the  formula
\[\big[ t_\mathbf{x}, t_{\mathbf{x_1}} \cdots t_{\mathbf{x_\ell}}\big] = \sum_{i=1}^{\ell} t_{\mathbf{x}_1} \cdots t_{\mathbf{x}_{i-1}} \big[ t_{\mathbf{x}} , t_{\mathbf{x}_i} \big] t_{\mathbf{x}_{i+1}} \cdots t_{\mathbf{x}_{\ell}}\] 
combined with $ t_{\mathbf{z}} t_{\mathbf{w}} = \big[t_{\mathbf{z}}, t_{\mathbf{w}}  \big] + t_{\mathbf{w}} t_{\mathbf{z}} $ for $\mu(\mathbf{z}) > \mu(\mathbf{w})$ and with the fact (previous step) that we can calculate $\big[t_{\mathbf{z}}, t_{\mathbf{w}}  \big]$ using the relation (2) of Definition \ref{defpathalgebra} results in finitely many steps in
  \[\big[ t_\mathbf{x}, t_{\mathbf{x_1}} \cdots t_{\mathbf{x_\ell}}\big] = c_1 t_{\mathbf{p}_1} + \cdots + c_s t_{\mathbf{p}_s}\]  
for $\mathbf{p}_1, \ldots, \mathbf{p}_s \in \mathbf{Conv}^{+}.$ Therefore, 
\[ \big[\mathcal{K}_{x}^{\oplus r}, \E\big] = \sum_{i=1}^{m} a_{i}\; T_{\mathbf{v}_{i_1}}^{\widetilde{\rho}_{i_1}} \cdots T_{\mathbf{v}_{i_{\ell}}}^{\widetilde{\rho}_{i_{\ell}}} \]
in increasing order of slopes i.e.\ $\mu(\mathbf{v}_{i_1}) < \cdots < \mu(\mathbf{v}_{i_{\ell}})$, for $i=1, \ldots,m$. By construction (or Theorem \ref{th4}) $\mathbf{p}_i = (\mathbf{v}_{i_1}, \ldots, \mathbf{v}_{i_{\ell}} )$ runs through the convex paths in $\Delta_{(\rk(\E),\deg(\E)),(0,r|x|)}$.\\

\paragraph*{\textbf{Step 8}}  Next we replace the $T_{\mathbf{v}_{i_j}}^{\widetilde{\rho}_{i_j}}$ by their definition. Using that $\sum_{\rho} \rho(x)$ is zero unless $x=x_0$, we can write the product  $\mathcal{K}_{x}^{\oplus r} \; \E$ as linear combination of $ T_{\mathbf{v}_{i_1}, x_{i_1}} \cdots T_{\mathbf{v}_{i_n},x_{i_n}},$ where $\mathbf{p}_i = (\mathbf{v}_{i_1}, \ldots, \mathbf{v}_{i_n} )$ runs through the convex paths in the triangle defined by  $\mathcal{K}_{x}^{\oplus r}$ and $\E$. \\

\paragraph*{\textbf{Step 9}} Finally, since $T_{\mathbf{v}_i, x_i} T_{\mathbf{v}_j,x_j} = T_{\mathbf{v}_i,x_i} \oplus T_{\mathbf{v}_j,x_j}$ if either $\mu(\mathbf{v}_i) < \mu(\mathbf{v}_j)$ or  $\mu(\mathbf{v}_i) = \mu(\mathbf{v}_j)$ and $x_i \neq x_j,$ our problem reduces to solve the case $T_{\mathbf{v}_i,x}  T_{\mathbf{v}_j,x}$ with $\mu(\mathbf{v}_i) = \mu(\mathbf{v}_j).$\\

\paragraph*{\textbf{Step 10}} In order to solve the previous step consider  $T_{\mathbf{v}_i,x} , T_{\mathbf{v}_j,x} \in \mathbf{H}_{\mathrm{Tor}_x}$ with $\mu(\mathbf{v}_i) = \mu(\mathbf{v}_j).$ By Proposition \ref{prop-macdonald} once again we can express the product $T_{\mathbf{v}_i,x}  T_{\mathbf{v}_j,x}$ as a product of power-sums. Since the Hall-Littlewood functions form a $\C[u_x^{\pm 1}]$-basis for the Macdonald ring of symmetric functions we can write the product of these two power-sums as a linear combination of Hall-Littlewood functions. Since the Hall-Littlewood functions $P_{\lambda}$ correspond to $\mathcal{K}_{x}^{(\lambda)}$, we have an explicit description for the product  $T_{\mathbf{v}_i,x} T_{\mathbf{v}_j,x},$ and therefore an explicit description for the product $\pi^{\mathrm{vec}}\big( \mathcal{K}_{x}^{\oplus r} \; \E \big).$  \\

\begin{remark} In fact, our algorithm does not only calculate the products $\mathcal{K}_{x}^{\oplus r} \E$, but also a more general class of products in the Hall algebra of an elliptic curve. For instance the products of any torsion sheaf by a vector bundle. 
\end{remark}

\section{Calculating structure constants}\label{section-structureconstants}

In this section we will apply our algorithm to calculate some structure constants for the graphs $\mathcal{G}_{x,r}.$ Our results cover some cases of particular interest and are meant to exemplify the way the algorithm works. An explicit description of the graphs $\mathcal{G}_{x,r}$ for all $n$, all $x$ and all $r$ seems out of reach. But we cover the case of $|x|=1, \; n=2,\; r=1$ in completion in next section.   

\begin{theorem} \label{propstable} Fix an integer $n \geq 1$. Let  $x \in |X|$ with $|x|=1.$ Let $\E = \E_{(y,1)}^{(n,d)}$ with $n|d$ be a stable vector bundle on $X.$ Then $m_{x,1}\big(\E',\E \big) \neq 0$ if and only if $\E' \cong \E_{(z,1)}^{(n,d+1)}$ where $z = x + \mathrm{Norm}_{1}^{n}(y'),$ and $y' \in |X_n|$ sits above $y.$ In this case $m_{x,1}\big(\E',\E \big) =1.$  Graphically,
\[
  \beginpgfgraphicnamed{tikz/fig3}
  \begin{tikzpicture}[>=latex, scale=2]
        \vertex[circle,fill,label={below:$\E_{(z,1)}^{(n,d+1)}$}](00) at (0,0) {};
        \vertex[circle,fill,label={below:$\E_{(y,1)}^{(n,d)}$}](10) at (2,0) {};
    \path[-,font=\scriptsize]
    (00) edge[->-=0.8] node[pos=0.2,auto,black] {\normalsize $1$} (10)
    ;
   \end{tikzpicture}
 \endpgfgraphicnamed
\] 
are all incoming arrows at $\E_{(y,1)}^{(n,d)}$ in $\mathcal{G}_{x,1}.$

\begin{proof} Via Atiyah's classification, $\E = \E_{(y,1)}^{(n,d)}$ corresponds to $\mathcal{K}_y$ where $|y|=n.$ 
We aim to calculate $\pi^{\mathrm{vec}}(\mathcal{K}_{y} \; \E)$  following the algorithm from section \ref{sectionthealgorithm}. 

First note that $T_{(0,n),y} = [n]\mathcal{K}_{y}^{}.$ For $\rho \in \widehat{\mathrm{Pic}^{0}(X_n)}$ consider
\[T_{(0,n)}^{\widetilde{\rho}} = \sum_{\substack{x' \in |X|\\ |x'|\;| \;n}} \widetilde{\rho}(x') T_{(0,n),x'} \]
and by Proposition \ref{inversionlemma}
\[T_{(0,n),y} = n N_{n}^{-1} \sum_{\widetilde{\rho} \in \mathrm{P}_n} \widetilde{\rho}(-y) T_{(0,n)}^{\widetilde{\rho}}.\]   
Hence
\[\mathcal{K}_{y}^{} = [n]^{-1} n N_{n}^{-1} \sum_{\widetilde{\rho} \in \mathrm{P}_n} \widetilde{\rho}(-y) T_{(0,n)}^{\widetilde{\rho}}\]
and
\[\E = [n]^{-1}n N_{n}^{-1} \sum_{\widetilde{\rho} \in \mathrm{P}_n} \widetilde{\rho}(-y) T_{(n,d)}^{\widetilde{\rho}}.\]
Therefore
\begin{eqnarray*} 
\pi^{\mathrm{vec}}(\mathcal{K}_{x} \; \E) & = & \big[\mathcal{K}_{x}, \E\big] \\
& = & \Big[N_{1}^{-1} \sum_{\sigma \in \widehat{\mathrm{Pic}^{0}(X)}} \sigma(-x) T_{(0,1)}^{\sigma}, [n]^{-1}n N_{n}^{-1} \sum_{\widetilde{\rho} \in \mathrm{P}_n} \widetilde{\rho}(-y) T_{(n,d)}^{\widetilde{\rho}}\Big]\\
& = & N_{1}^{-1} [n]^{-1} n N_{n}^{-1} \sum_{\substack{\sigma \in \mathrm{P}_1 \\ \widetilde{\rho} \in \mathrm{P}_n}} \sigma(-x)\rho(-y) \big[ T_{(0,1)}^{\sigma}, T_{(n,d)}^{\widetilde{\rho}} \big] \\
& = & N_{1}^{-1} [n]^{-1} n N_{n}^{-1} c_n \sum_{\substack{\sigma \in \mathrm{P}_1  \\ \rho = \mathrm{Norm}_{1}^{n}(\sigma)}} \sigma(-x-\mathrm{Norm}_{1}^{n}(y'))\; T_{(n,d+1)}^{\sigma} \\
& = & N_{1}^{-1} v^{n}\sum_{\sigma \in \mathrm{P}_1}  \sum_{x' \in |X|} \sigma(x'-x-\mathrm{Norm}_{1}^{n}(y')) T_{(n,d+1),x'} \\
& = & v^n \;T_{(n,d+1),z} \\
& = & v^n \;\E_{(z,1)}^{(n,d+1)},
\end{eqnarray*}
where $z = x + \mathrm{Norm}_{1}^{n}(y')$ with $y' \in |X_n|$ sits above $y.$ 
\end{proof}
\end{theorem}

Let us stress that the equality $z = x + \mathrm{Norm}_{1}^{n}(y')$ in the above theorem, as well as, throughout the paper, means that  
$$ z- x + \mathrm{Norm}_{1}^{n}(y') \equiv \mathcal{O}_X$$
via $\mathrm{Pic}^m X \equiv \mathrm{Pic}^0X,$ where $m = \deg( z- x + \mathrm{Norm}_{1}^{n}(y'))$, as explained in section \ref{section1}.

\begin{theorem} \label{neighbousoftotaldec}Let $x$ be a closed point in $X$ and $ \Line_1, \ldots, \Line_n$  line bundles on $X$. If $\E' = \Line_1 \oplus \cdots \oplus \Line_n$ with $\deg(\Line_i) > \deg(\Line_{i-1}) +|x|$ for $i=1, \ldots,n$, then 
\[\mathcal{V}_{x,1}(\E') = \bigcup_{k=1}^{n} \big\{ \big(\E',\Line_1 \oplus \cdots \oplus \Line_k (- x) \oplus \cdots \oplus \Line_n, \; q_{x}^{k-1} \big)               \big\}.\]

\begin{proof} Following the algorithm in section \ref{sectionthealgorithm}, the proposition will follow from calculating $ \pi^{\mathrm{vec}}(\mathcal{K}_{x}^{} \E)=[\mathcal{K}_{x}^{},\E].$ We have 
$$\E = \E_{(x_1,1)}^{(1,d_1)} \oplus \cdots \oplus \E_{(x_n,1)}^{(1,d_n)} = v^{\sum_{i<j}\langle (1,d_i),(1,d_j)\rangle} \E_{(x_1,1)}^{(1,d_1)}  \cdots \E_{(x_n,1)}^{(1,d_n)}. $$

Denote by $d$ the degree of $x.$ By definition, $T_{(0,d),x} = [d]\mathcal{K}_{x}^{}.$ By Proposition \ref{inversionlemma},
\[\mathcal{K}_{x}^{} = [d]^{-1} d N_{d}^{-1} \sum_{\widetilde{\sigma} \in \widetilde{\mathrm{Pic}^0(X_d)}} \widetilde{\sigma}(-x) T_{(0,d)}^{\widetilde{\sigma}}.\]
By Atiyah's classification (Theorem \ref{atiyahtheorem}) $\E_{(x_i,1)}^{(1,d_i)}$ corresponds to $\mathcal{K}_{x_i}^{}$ where $|x_i|=1$ for $i=1, \ldots,n.$ Since $T_{(0,1),x_i} = \mathcal{K}_{x_i}^{}$, we have for $i=1, \ldots,n$ that
\[\mathcal{K}_{x_i}^{} = N_{1}^{-1} \sum_{\widetilde{\rho_i} \in \widetilde{\mathrm{Pic}^0(X)}} \widetilde{\rho_i}(-x_i) T_{(0,1)}^{\widetilde{\rho_i}}\]
by Proposition \ref{inversionlemma}. Thus 
\[\E_{(x_i,1)}^{(1,d_i)} =  N_{1}^{-1} \sum_{\widetilde{\rho_i} \in \widetilde{\mathrm{Pic}^0(X)}} \widetilde{\rho_i}(-x_i) T_{(1,d_i)}^{\widetilde{\rho_i}},\quad \text{ for $i=1, \ldots,n$}.\]   
Therefore
\begin{eqnarray*} 
 \pi^{\mathrm{vec}}(\mathcal{K}_{x}^{} \E) 
& = &  \Big[ [d]^{-1} d N_{d}^{-1} \sum_{\widetilde{\sigma} \in \mathrm{P}_d} \widetilde{\sigma}(-x) T_{(0,d)}^{\widetilde{\sigma}}, \\
& & v^{\sum_{i<j}\langle (1,d_i),(1,d_j)\rangle} N_{1}^{-n} \sum_{i=1}^{n}\sum_{ \widetilde{\rho_i} \in \mathrm{P}_1} \widetilde{\rho_1}(-x_1) \cdots \widetilde{\rho_n}(-x_n)T_{(1,d_1)}^{\widetilde{\rho_1}} \cdots T_{(1,d_n)}^{\widetilde{\rho_n}} \Big] \\
& = &  [d]^{-1} d N_{d}^{-1}v^{\sum_{i<j}\langle (1,d_i),(1,d_j)\rangle} N_{1}^{-n}  \sum_{i=1}^{n}\sum_{ \substack{ \widetilde{\rho_i} \in \mathrm{P}_1 \\ \widetilde{\sigma} \in \mathrm{P}_d}} \widetilde{\sigma}(-x) \widetilde{\rho_1}(-x_1) \cdots \widetilde{\rho_n}(-x_n)\\ & &\Big[ T_{(0,d)}^{\widetilde{\sigma}}, T_{(1,d_1)}^{\widetilde{\rho_1}} \cdots T_{(1,d_n)}^{\widetilde{\rho_n}} \Big] \\
&=&  [d]^{-1} d N_{d}^{-1}v^{\sum_{i<j}\langle (1,d_i),(1,d_j)\rangle} N_{1}^{-n} \Big( \sum_{i=1}^{n}\sum_{ \substack{ \widetilde{\rho_i} \in \mathrm{P}_1 \\ \widetilde{\sigma} \in \mathrm{P}_d}} \widetilde{\sigma}(-x) \widetilde{\rho_1}(-x_1) \cdots \widetilde{\rho_n}(-x_n)\\ & & \Big[T_{(0,d)}^{\widetilde{\sigma}}, T_{(1,d_1)}^{\widetilde{\rho_1}}\Big]T_{(2,d_2)}^{\widetilde{\rho_2}} \cdots T_{(1,d_n)}^{\widetilde{\rho_n}} + \cdots + \\
& & \sum_{i=1}^{n}\sum_{ \substack{ \widetilde{\rho_i} \in \mathrm{P}_1 \\ \widetilde{\sigma} \in \mathrm{P}_d}} \widetilde{\sigma}(-x) \widetilde{\rho_1}(-x_1) \cdots \widetilde{\rho_n}(-x_n) T_{(1,d_1)}^{\widetilde{\rho_1}} \cdots T_{(1,d_{n-1})}^{\widetilde{\rho_{n-1}}} \Big[T_{(0,d)}^{\widetilde{\sigma}}, T_{(1,d_n)}^{\widetilde{\rho_n}} \Big]\Big).   
\end{eqnarray*}
Observe that $\big[T_{(0,d)}^{\widetilde{\sigma}}, T_{(1,d_i)}^{\widetilde{\rho_i}}\big]=0 $ unless $\sigma = \mathrm{Norm}_{1}^{d}\rho_i.$ When $\sigma = \mathrm{Norm}_{1}^{d}\rho_i$, we have by Theorem \ref{dragosmaintheorem} that
\[\big[T_{(0,d)}^{\widetilde{\sigma}}, T_{(1,d_i)}^{\widetilde{\rho_i}}\big] = c_d\; T_{(1,d_i+d)}^{\widetilde{\rho_i}}.\]
This yields
\begin{eqnarray*} 
\pi^{\mathrm{vec}}(\mathcal{K}_{x}^{} \E) & = & c_d [d]^{-1} d N_{d}^{-1}v^{\sum_{i<j}\langle (1,d_i),(1,d_j)\rangle} N_{1}^{-n} \cdot\\
&\cdot & \Big( \sum_{i=1}^{n}\sum_{ \substack{ \widetilde{\rho_i} \in \mathrm{P}_1 \\ \widetilde{\sigma} =  \mathrm{Norm}_{1}^{d}\rho_i }} \widetilde{\rho_1}(-x_1-\mathrm{Norm}_{1}^{d}(x'))\widetilde{\rho}_2(-x_2) \cdots \widetilde{\rho_n}(-x_n) T_{(1,d_i)}^{\widetilde{\rho_i}} T_{(2,d_2)}^{\widetilde{\rho_2}} \cdots T_{(1,d_n)}^{\widetilde{\rho_n}} + \\
& \cdots & + \Big( \sum_{i=1}^{n}\sum_{ \substack{ \widetilde{\rho_i} \in \mathrm{P}_1 \\ \widetilde{\sigma} = \mathrm{Norm}_{1}^{d}\rho_n}}
\widetilde{\rho_1}(-x_1) \cdots  \widetilde{\rho}_{n-1}(-x_{n-1}) \widetilde{\rho_n}(-x_n- \mathrm{Norm}_{1}^{d}(x'))\Big) \cdot \\ 
& & \cdot\; T_{(1,d_1)}^{\widetilde{\rho_1}} \cdots T_{(1,d_{n-1})}^{\widetilde{\rho}_{n-1}} T_{(1,d_n +d)}^{\widetilde{\rho_n}}\Big)\\
& = & N_{1}^{-n} v^{\sum_{i<j}\langle (1,d_i),(1,d_j)\rangle + d}\\
& & \Big( \sum_{i=1}^{n}\sum_{ \substack{ \widetilde{\rho_i} \in \mathrm{P}_1 \\ \widetilde{\sigma} =  \mathrm{Norm}_{1}^{d}\rho_i }} \sum_{x_{i}^{1} \in |X|}\widetilde{\rho_1}(x_{1}^{1}-x_1-\mathrm{Norm}_{1}^{d}(x'))\widetilde{\rho}_2(x_{2}^{1}-x_2) \cdots \widetilde{\rho_n}(x_{n}^{1}-x_n)\\ 
& & T_{(1,d_1+d),x_{1}^{1}} T_{(2,d_2),x_{2}^{1}} \cdots T_{(1,d_n),x_{n}^{1}} + \cdots \\
&  +  & \sum_{i=1}^{n}\sum_{ \substack{ \widetilde{\rho_i} \in \mathrm{P}_1 \\ \widetilde{\sigma} = \mathrm{Norm}_{1}^{d}\rho_n}} \sum_{x_{i}^{n} \in |X|}\widetilde{\rho_1}(x_{1}^{n}-x_1) \cdots  \widetilde{\rho}_{n-1}(x_{n-1}^{n}-x_{n-1}) \widetilde{\rho_n}(x_{n}^{n}-x_n - \mathrm{Norm}_{1}^{d}(x'))\\ 
& & T_{(1,d_1),x_{1}^{n}} \cdots T_{(1,d_{n-1}),x_{n-1}^{n}} T_{(1,d_n +d),x_{n}^{n}}\Big).  \\
& = & N_{1}^{-n} v^{\sum_{i<j}\langle (1,d_i),(1,d_j)\rangle + d} \Big( N_{1}^{n} \sum_{k=1}^{n} v^{c_k} \E_{(x_{1}^{k},1)}^{(1,d_1)} \oplus \cdots \oplus \E_{(x_{k}^{k},1)}^{(1,d_k + d)} \oplus \cdots \oplus \E_{(x_{n}^{k},1)}^{(1,d_n)} \Big)
\end{eqnarray*}
where $x' \in |X_d|$ sits above $x$, $x_{i}^{j} = x_i$ if $i \neq j$ and $x_{i}^{i} = x_i \oplus \mathrm{Norm}_{1}^{d}(x')$ for $i=1, \ldots,n$ and $$c_k=-\sum_{j=k+1}^{n}\langle (1,d_k + d),(1,d_j)\rangle - \sum_{i=1}^{k-1}\langle (1,d_i),(1,d_k+d) \rangle - \sum_{\substack{i<j \\ i,j \neq k}}\langle (1,d_i),(1,d_j)\rangle.$$
The multiplicities are determined by the formula 
\[\sum_{i<j}\langle (1,d_i),(1,d_j)\rangle + d + c_k = d(n-2(k-1)).\]
and by the multiplication for $v^{-nd}$, cf.\ Remark \ref{propconnection}. 
\end{proof}
\end{theorem}

We end this section with a special case of edges in the graphs of Hecke operators. We proof the case $n=3$ but taking the coefficient of $s^n$ in the power series 
\[\exp\left( m(v^{-1}-v) \sum_{i \geq 1} t_{i \mathbf{z}_{0}} s^i \right)\]
(recall Definition \ref{defpathalgebra})the same proof can be applied for any stable bundle $\E_{(x',1)}^{(n,d)}$ where $n$ is a prime number and $d \equiv -1 (\mathrm{mod}\; n).$

\begin{theorem} Let $x$ be a degree one closed point on $X.$ Let  $\E := \E_{(x',1)}^{(3,d)}$ with $d \equiv -1 (\mathrm{mod}\; 3).$ Then
\[m_{x,1}\big(\E_{(y,1)}^{(3,d+1)} ,\E\big)= q^2 + q +1, \quad \quad m_{x,1}\big(\E_{(z,3)}^{(3,d+1)} ,\E\big)= q^2  \]
\[m_{x,1}\big(\E_{(y',1)}^{(1,d')} \oplus \E_{(y'',1)}^{(2,2d')} ,\E\big)= q^2 -1, \quad \quad
m_{x,1}\big(\E_{(z',1)}^{(1,d')} \oplus \E_{(z'',2)}^{(2,2d')} ,\E\big)= q^2 -q \]
\[ \text{ and } \quad m_{x,1}\big(\E_{(y_1,1)}^{(1,d')} \oplus \E_{(y_2,1)}^{(1,d')} \oplus \E_{(y_3,1)}^{(1,d')}  ,\E\big)= q^2 -2q+1,\]
where $d'=(d+1)/3$, $|y|=3$ with $\mathrm{Norm}_{1}^{3}(w) = x+x'$ where $w \in |X_3|$ sits above $y$; $|z|=1$ with $3z=x+x'$; $|y'|=1, |y''|=2$ with $y' +\mathrm{Norm}_{1}^{2}(w'')=x+x'$ where $w'' \in |X_2|$ sits above $y''$; $|z'|=|z''|=1$ with $z'+2z''=x+x';$ and $|y_1|=|y_2|=|y_3|=1, y_1\neq y_2 \neq y_3$ with $y_1+y_2+y_3 =x+x'.$  

\begin{proof} Observe that by Theorem \ref{atiyahtheorem}, $\E$ corresponds to $\mathcal{K}_{x'}^{}.$ It follows from the definition that 
$\mathcal{K}_{x'}^{} = T_{(0,1),x'}$, and by Proposition \ref{inversionlemma}
\[\E = N_{1}^{-1} \sum_{\widetilde{\rho} \in \mathrm{P}_1} \widetilde{\rho}(-x') T_{(3,d)}^{\widetilde{\rho}}.\]
Thus,
\[\pi^{\mathrm{vec}}(\mathcal{K}_{x}^{}\;\E) = N_{1}^{-2} \sum_{\sigma \in \mathrm{P}_1}\sigma(-x-x')\big[ T_{(0,1)}^{\sigma}, T_{(3,d)}^{\sigma}\big].\]
Theorem \ref{dragosmaintheorem} and Definition \ref{defpathalgebra} now yields
\[\big[ T_{(0,1)}^{\sigma}, T_{(3,d)}^{\sigma}\big] = c_1 (v^{-1}-v)^{-1} \theta_{(3,d+1)},\]
where
\[\theta_{(3,d+1)} = (v^{-1}-v) T_{(3,d+1)}^{\mathrm{Norm}_{1}^{3}\sigma} + (v^{-1}-v)^{2} T_{(1,d')}^{\sigma}\; T_{(2,2d')}^{\mathrm{Norm}_{1}^{2}\sigma} + \frac{(v^{-1}-v)^{3}}{6}\big( T_{(1,d')}^{\sigma}\big)^3 .\]
Let $w \in |X_3|$ sits above $y$ and $w'' \in |X_2|$ sits above $y''.$ Hence, 
\begin{eqnarray*}
v^{-3}\pi^{\mathrm{vec}}(\mathcal{K}_{x}^{}\;\E) &=&v^{-2}  \sum_{\substack{y; |y||3 \\ \mathrm{Norm}_{1}^{3}(w)=x+x'}}  T_{(3,d+1),y}\\ & + &v^{-2}(v^{-1}-v) \sum_{\substack{ y',y'' \\ |y'|=1,\; |y''|=2 \\ y' + \mathrm{Norm}_{1}^{2}(w'')=x+x'}} T_{(1,d'),y'} T_{(1,2d'),y''} \\
&+&\frac{v^{-2}(v^{-1}-v)^{2}}{6} \sum_{\substack{y_1,y_2,y_3 \\ |y_1|=|y_2|=|y_3|=1 \\ y_1+y_2+y_3=x+x'}}\; T_{(1,d'),y_1} T_{(1,d'),y_2} T_{(1,d'),y_3}, 
\end{eqnarray*}
where $d'=(d+1)/3.$ Let us calculate each term separated. 

Observe that $T_{(3,d+1),y} \neq 0$ if, and only if, $|y|=1$ or $|y|=3.$ If $|y|=3$, then 
\[v^{-2} \;T_{(3,d+1),y} = v^{-2}[3]\; \E_{(y,1)}^{(3,d+1)} = (q^2 + q +1)\; \E_{(y,1)}^{(3,d+1)}.\]
If $|y|=1$, then 
\begin{eqnarray*}
v^{-2} \;T_{(3,d+1),y} &=& 3^{-1}(v^{-4}+ v^{-2}+1)\; \E_{(y,3)}^{(3,d+1)} + 3^{-1}(1-v^{-6})\; \E_{(y,1)}^{(1,d')} \oplus \E_{(y,2)}^{(2,2d')}\\ 
&+& 3^{-1}(1-v^{-4}-v^{-6}+ v^{-10})\; \E_{(y,1)}^{(1,d')} \oplus \E_{(y,1)}^{(1,d')} \oplus \E_{(y,1)}^{(1,d')}. 
\end{eqnarray*}

Next we calculate the second term. We have $T_{(1,d'),y'}= \E_{(y',1)}^{(1,d')}$ and
\[ T_{(2,2d'),y''} =  \begin{cases} [2] \;\E_{(y'',1)}^{(2,2d')} & \text{ if } |y''|=2 \\
\frac{[2]}{2} \left( \E_{(y'',2)}^{(2,2d')} + (1-v^{-2}) \; \E_{(y'',1)}^{(1,d')} \oplus \E_{(y'',1)}^{(1,d')} \right) & \text{ if } |y''|=1.
\end{cases}  \]
If $y' \neq y''$, then $T_{(1,d'),y'} T_{(2,2d'),y''} =T_{(1,d'),y'} \oplus T_{(2,2d'),y''}$, and this case follows from above calculation. 
If $y' = y''$, the product $T_{(1,d'),y'} T_{(2,2d'),y''}$ is equivalent to the product $\frac{[2]}{2} p_1 p_2$ of power-sums in the Macdonald's ring. Since
\[ p_1 p_2 = (v^6 -1)\;P_{(1^3)} + v^2 \; P_{(2,1)} + P_{(3)}, \]
for  $y' = y''$ we obtain 
\[T_{(1,d'),y'} T_{(2,2d'),y''} = \frac{[2]}{2} \left((v^6 -1)v^{-6}\;\Big(\E_{(y',1)}^{(1,d')}\Big)^{\oplus 3} +  \E_{(y',2)}^{(2,2d')} \oplus \E_{(y',1)}^{(1,d')} + \E_{(y',3)}^{(3,d+1)} \right). \]
Therefore the second terms yields
\[(v^{-4}-1)\; \E_{(y',1)}^{(1,d')} \oplus\E_{(y'',1)}^{(2,2d')}   \]
if $y'\neq y'' \text{ and } |y''|=2$;
\[\frac{1}{2}(v^{-4}-1) \; \E_{(y',1)}^{(1,d')} \oplus \E_{(y'',2)}^{(2,2d')} + \frac{1}{2}(-v^{-6} + v^{-4}+ v^{-2}-1) \E_{(y',1)}^{(1,d')} \oplus\Big( \E_{(y'',1)}^{(1,d')}  \Big)^{\oplus 2}  \]
if $ y'\neq y'' \text{ and } |y''|=1;$ and
\[ \frac{1}{2}(-v^{-10}+ v^{-6} + v^{-4}-1)\;\Big(\E_{(y',1)}^{(1,d')}\Big)^{\oplus 3} + \frac{1}{2}(v^{-4}-1)\; \E_{(y',2)}^{(2,2d')} \oplus \E_{(y',1)}^{(1,d')} + \frac{1}{2}(v^{-4}-1)\; \E_{(y',3)}^{(3,d+1)} \]
 for $y'=y''$.
 
We separated the third term $T_{(1,d'),y_1} T_{(1,d'),y_2} T_{(1,d'),y_3}$ in three cases. The first one is when $y_1 \neq y_2 \neq y_3$, in this case we have 
\[T_{(1,d'),y_1} T_{(1,d'),y_2} T_{(1,d'),y_3} = \E_{(y_1,1)}^{(1,d')} \oplus  \E_{(y_2,1)}^{(1,d')} \oplus  \E_{(y_3,1)}^{(1,d')}.\]
If $y_i = y_j$ and $y_j \neq y_k$ for $\{i,j,k\}=\{1,2,3\},$ we use the following equation 
\[p_{1}^{2} = (v^2 +1 ) P_{(1^2)} + P_{(2)}\]
to conclude that the third term yields
\[\frac{1}{2}(v^{-6}-v^{-4}-v^{-2}+1) \; \Big( \E_{(y_1,1)}^{(1,d')} \Big)^{\oplus 2} \oplus \E_{(y_2,1)}^{(1,d')}    + \frac{1}{2}(v^{-4} -2v^{-2} +1) \; \E_{(y_1,2)}^{(2,2d')} \oplus \E_{(y_2,1)}^{(1,d')}.\]
The last case is when $y_1 = y_2 =y_3.$ Since 
\[p_{1}^{3} = (v^6 + 2 v^4 + 2v^2 +1) P_{(1^3)} + (v^2 +2) P_{(2,1)} + P_{(3)},\]
we conclude that the third case is equals to
\[\frac{1}{6}(v^{-10} - v^{-6} - v^{-4}+1)\; \Big(\E_{(y_1,1)}^{(1,d')}\Big)^{\oplus 3} + 
\frac{1}{6}(2v^{-6}-3v^{-4}+1) \;\E_{(y_1,2)}^{(2,2d')} \oplus \E_{(y_1,1)}^{(1,d')} \;+\]
\[+ \frac{1}{6}(v^{-4}-2 v^{-2}+1) \;\E_{(y_1,3)}^{(3,d+1)}.\]
Putting all together we have the desired.  
\end{proof}
\end{theorem}

\section{The case of rank 2}\label{rank2}

\noindent
We shall use our algorithm to describe the graphs of unramified Hecke operators for $n=2$, $r=1$ and  a closed point $x$ of degree one.

We start with an explicit description for the vertices of $\mathcal{G}_{x,1}$ for $n=2.$ In \cite{oliver-graphs} and \cite{oliver-elliptic}, Lorscheid describe explicitly the set $\mathbb{P}\mathrm{Bun}_2 X$ of isomorphism classes of $\mathbb{P}^{1}$-bundles over $X$ (\cite{biblia}, Ex. II.7.10). We derive a description of $\mathrm{Bun}_2 X$ from this by determining the of representatives for each classes in $\mathbb{P}\mathrm{Bun}_2 X$.

We call a vector bundles $\E$ indecomposable if for every decomposition $\E = \E_1 \oplus \E_2$ into two subbundles $\E_1, \E_2$ one factor is trivial and the other is isomorphic to $\E.$ The Krull-Schmidt theorem (c.f. \cite{atiyah-krull} Theorem 2) holds for the category of vector bundles over $X$, i.e.\ every vector bundle on $X$  has a unique decomposition into a direct sum of indecomposable subbundles, up to permutation of factors. 

An extension of scalars $\mathbb{F}_{q^i}F/F$, or geometrically $\pi: X_{i} := X \times_{\Spec\mathbb{F}_q} \Spec \mathbb{F}_i \rightarrow X$, defines the inverse image or the \textit{constant extension} of vector bundles 
\begin{eqnarray*}
\pi^{*}: \mathrm{Bun}_n X &\longrightarrow & \mathrm{Bun}_n X_i . \\
\E & \longmapsto & \pi^{*} \E 
\end{eqnarray*}
The isomorphism  classes of rank $n$ bundles that after extension of constants to $\mathbb{F}_{q^i}$ become isomorphic to  $\pi^{*}\E$ are classified by $H^1(\mathrm{Gal}(\mathbb{F}_{q^i}/\Fq), \mathrm{Aut}(\E \otimes \mathbb{F}_{q^i}))$, c.f.\ \cite{arason} Section 1. Since the algebraic group $\mathrm{Aut}(\E \otimes \mathbb{F}_{q^i})$ is an open subvariety of the connected algebraic group $\mathrm{End}(\E \otimes \mathbb{F}_{q^i}),$ it is itself a connected algebraic group. As a consequence of Lang's theorem (\cite{lang} Corollary to Theorem 1), we have  $H^1(\mathrm{Gal}(\mathbb{F}_{q^i}/\Fq), \mathrm{Aut}(\E \otimes \mathbb{F}_{q^i})) =1.$ We deduce that $\pi^{*}$ is injective. In particular, one can consider the constant extension to the geometric curve $\overline{X} = X \times_{\mathrm{Spec}\mathbb{F}_q} \mathrm{Spec}\overline{\Fq}$, where $\overline{\Fq}$ is an algebraic closure of $\Fq.$ Then two vector bundles are isomorphic if and only if they are geometrically isomorphic, i.e.\ their constant extensions to $\overline{X}$ are isomorphic. We can therefore think $\mathrm{Bun}_n X$ as a subset of $\mathrm{Bun}_n X_i$ and $\mathrm{Bun}_n \overline{X}.$ 

On the other hand, $\pi: X_i \rightarrow X$ defines the direct image or the \textit{trace} of vector bundles
\begin{eqnarray*}
\pi_{*}: \mathrm{Bun}_n X_i &\longrightarrow & \mathrm{Bun}_{ni} X. \\
\E & \longmapsto & \pi_{*} \E
\end{eqnarray*}
For $\E \in  \mathrm{Bun}_n X$ we have that $\pi_{*}\pi^{*}\E = \E^{\oplus i}$ and for $\E \in  \mathrm{Bun}_n X_i$ that $\pi^{*}\pi_{*}\E = \bigoplus \E^{\tau}$
where $\tau$ ranges over $\mathrm{Gal}(\mathbb{F}_{q^i} /\Fq )$ and $\E^{\tau}$ is defined by the stalks $\E_{x}^{\tau} = \E_{\tau^{-1}(x)}$.
We call a vector bundle geometrically indecomposable if its extension to $\overline{X}$ is indecomposable. 
In \cite[Theorem 1.8]{arason}, it is shown that every indecomposable vector bundle over $X$ is the trace
of a geometrically indecomposable bundle over some constant extension $X_i$ of $X.$

We are interested in the case $i=n=2$.The set $\mathrm{Bun}_2 X$ is the disjoint union of the set of classes of decomposable rank 2 bundles and the set of classes of indecomposable bundles. We denote theses sets by $\mathrm{Bun}_{2}^{\mathrm{dec}}X$ and $\mathrm{Bun}_{2}^{\mathrm{ind}}X,$ respectively. Let $\mathrm{Bun}_{2}^{\mathrm{gi}}X \subset \mathrm{Bun}_{2}^{\mathrm{ind}}X$ be the subset of classes of geometrically indecomposable vector bundles. Since the rank is $2$, the complement $\mathrm{Bun}_{2}^{\mathrm{tr}}X = \mathrm{Bun}_{2}^{\mathrm{ind}}X - \mathrm{Bun}_{n}^{\mathrm{gi}}X$ consists of classes  of traces $ \pi_{*}\Line '$ where $\Line' \in \mathrm{Pic} X_2.$ Moreover, $\pi_{*}\Line '= \pi_{*}\Line ''$ if, and only if, $\Line '' \in \{\Line ', \; (\Line')^{\tau}\},$ where $\tau$ generates $\mathrm{Gal}(\mathbb{F}_{q^2}|\mathbb{F}_q)$ and $\pi_{*}\Line '$ decomposes if, and only if, $\Line' \in \mathrm{Pic}X$ (cf. Proposition 6.4 of \cite{oliver-graphs}). Thus we have a disjoint union
\[\mathrm{Bun}_{2}X = \mathrm{Bun}_{2}^{\mathrm{dec}}X \; \amalg \; \mathrm{Bun}_{2}^{\mathrm{tr}}X \; \amalg \;  \mathrm{Bun}_{2}^{\mathrm{gi}}X.\]
Namely, 
\[ \mathrm{Bun}_{2}^{\mathrm{dec}}X = \big\{\Line_{1} \oplus \Line_{2} \; \big| \; \Line_1,\Line_2 \in \mathrm{Pic}X\big\}\]
and 
\[\mathrm{Bun}_{2}^{\mathrm{tr}}X = \big\{\Line \otimes \pi_{*}\Line ' \; \big|\;\Line \in \mathrm{Pic} X,\; \Line' \in \mathrm{Pic}^{0}X_2 - \mathrm{Pic}^{0}X\big\}.\]
Next we determine $\mathrm{Bun}_{2}^{\mathrm{gi}}X.$ Since $g=1$ and $\omega_X = \mathcal{O}_X$, the Riemann-Roch theorem reduces to 
\[\dim_{\Fq} \Gamma(\Line) - \dim_{\Fq} \Gamma(\Line^{-1}) = \deg\Line.\]
Since $\Gamma(\Line)$ is non-zero if and only if $\Line$ is associated to an effective divisor (\cite{biblia} Proposition II.7.7(a)), we obtain: 
\[ \dim_{\Fq} \Gamma(\Line) = \begin{cases} 
 0 \quad \quad \hspace{0.2cm} \text{ if $\deg\Line \leq 0$ and $\Line \not\cong \mathcal{O}_X,$}\\
1  \quad \quad \hspace{0.3cm}  \text{if $\Line \cong \mathcal{O}_x,$ and} \\
\deg\Line \quad \text{if $\deg\Line > 0.$}
\end{cases}
\]
The units $\Fq^{*}$ operate by multiplication on the $\Fq-$vector space 
\[\mathrm{Ext}(\Line, \Line ') \cong \mathrm{Hom}(\Line, \Line' \omega_{X}^{\vee})\]
The multiplication of a morphism $\Line \rightarrow \Line ' \omega_{X}^{\vee}$ by an $a \in \Fq^{*}$ is nothing else but multiplying the generic stalk $\Line_{\eta}$ by $a^{-1}$ and all stalks $(\Line ' \omega_{X}^{\vee})_x$ at a closed point $x$ by $a.$ This induces automorphisms on both $\Line$ and $\Line ' \omega_{X}^{\vee}$, respectively. Thus, two elements of $\mathrm{Ext}^1(\Line, \Line ')$ that are $\Fq^{*}$-multiplies of each other define the same bundle on $\mathrm{Bun}_2X.$ We get a well-defined map
\[\mathbb{P}\mathrm{Ext}(\Line, \Line ') \longrightarrow \mathrm{Bun}_2 X\]
where the projective space $\mathbb{P}\mathrm{Ext}(\Line, \Line')$ is defined as the empty set when $\mathrm{Ext}(\Line, \Line')$ is trivial.

Let $\Line \in \mathrm{Pic}X.$ Serre duality asserts that $\mathrm{Ext}(\Line, \Line) \cong \mathrm{Hom}(\Line, \Line) \cong \Gamma(\mathcal{O}_{X})$ is one-dimensional. Thus $\mathbb{P}\mathrm{Ext}(\Line, \Line)$ contains only one element. By the above discussion, this determine a rank $2$-bundle $\E(\Line)$. For a place $x$ of degree one, the $\Fq$-vector space $\mathrm{Ext}(\Line(x),\Line) \cong \mathrm{Hom}(\Line,\Line(x)) \cong \Gamma( \mathcal{O}_{X}(x))$ is also one-dimensional and defines a rank $2$ bundle $\E_{x}(\Line).$ 

\begin{proposition} \[\mathrm{Bun}_{2}^{\mathrm{gi}}X = \big\{\E_{x}(\Line) \big| x \in |X|, |x|=1 \text{ and } \Line \in \mathrm{Pic}X\big\} \cup \big\{ \E(\Line) \big| \Line \in \mathrm{Pic}X\big\}.\]

\begin{proof} The proof follows from Proposition 7.1.4 of \cite{oliver-thesis} by removing the action of $\mathrm{Pic}X$. 
\end{proof}
\end{proposition}

Our next task is to prove the following theorem which describes the graph $\mathcal{G}_{x,1}$ for  a degree one place $x$ and $n=2.$ We consider $x$ to be the fixed point at $X$ as in section \ref{section1}. Hence, we have the bijection 
\[X(\Fq) \equiv \mathrm{Pic}^0 X \equiv \mathrm{Pic}^d X \]
which maps $z$ to $z-x$ and $z+(d-1)x,$ respectively, as $z$ varies through all degree one places while $x$ is fixed. 

\begin{theorem} \label{theoremcase2}Let $x$ be a closed point of $X$ of degree one. Then the edges of $\mathcal{G}_{x,1}$ for $n=2$ are given by the following list.

\begin{enumerate}
\item $\mathcal{V}_{x,1}(\E) = \Big\{ (\E, \Line_1(-x)\oplus \Line_2,\;1),(\E,\Line_1\oplus\Line_2(-x),\;q) \Big\} $ \\for $\E= \Line_1\oplus \Line_2$ and  either $\deg(\Line_2) - \deg(\Line_1) >1$ or $\deg(\Line_2) - \deg(\Line_1) =1$ and $\Line_2 \otimes \Line_1^{\vee} \not\simeq \mathcal{O}_{X}(x).$ \\

\item $\mathcal{V}_{x,1}(\E) = \Big\{ (\E, \Line_1(-x)\oplus \Line_2,\;1),(\E,\Line_1\oplus \Line_2(-x),\;1),  (\E, \E(\Line_1), \;q-1) \Big\} $\\ 
for $\E = \Line_1\oplus \Line_2, \;\deg(\Line_2) - \deg(\Line_1) =1$  and  $\Line_2(-x) \simeq \Line_{1}$. \\

\item $\mathcal{V}_{x,1}(\E) = \Big\{ (\E, \Line_1(-x)\oplus\Line_2, 1),(\E,\Line_1\oplus\Line_2(-x),1),(\E, \E_{x'}(\Line_2(-x)), q-1) \Big\}$ \\
for $\E = \Line_1 \oplus \Line_2, \; \deg(\Line_1)= \deg(\Line_2)$, $\Line_1 \not\simeq \Line_2$ and $x'$ is a degree one point such that $\Line_1 \otimes \Line_{2}^{-1} \simeq \mathcal{O}_X(x'-x)$.\\

\item $\mathcal{V}_{x,1}(\E) = \Big\{ (\E, \Line(-x) \oplus \Line, q+1) \Big\}$
for $\E = \Line \oplus \Line,  \; \Line \in \mathrm{Pic}X$.\\

\item $\mathcal{V}_{x,1}(\E(\Line )) = \Big\{ ( \E(\Line),\E_{x}(\Line(-x)),q), ( \E(\Line), \Line \oplus \Line(-x) ,1) \Big\} $\\
for $\Line \in \mathrm{Pic}X.$ \\ 

\item  $\mathcal{V}_{x,1}(\E_{x'}(\Line)) =  \bigcup_{\substack{ \Line'=\mathcal{O}_{X_2}(y-x) \in \mathrm{Pic}^{0}X_2 -   \mathrm{Pic}^{0}X \\ \mathrm{Norm}_{1}^{2}(y)=x'-x }}\Big\{ (\E_{x'}(\Line), \Line \otimes \pi_{*}(\Line'),1) \Big\} \\ 
\hspace*{2.3cm} \bigcup \Big\{ \;(\E_{x'}(\Line), \E(\Line), 1) \;|\; x'-x \in 2\mathrm{Pic}^0 X \Big\}\\ 
\hspace*{2.3cm}\bigcup_{\substack{ \Line_i \otimes \Line^{-1} \simeq \mathcal{O}_X(x_i - x), \; i=1,2 \\ x_1 \neq x_2, \; \mathcal{O}_X (x_1-x) \otimes \mathcal{O}_X(x_2-x) =\mathcal{O}_X( x'-x)}} \Big\{  (\E_{x'}(\Line), \Line_1\oplus \Line_2, 1)\Big\}$ \\
where $x'$ is a degree one point and $ \;\Line \in \mathrm{Pic}X.$ \\

\item $\mathcal{V}_{x,1}(\Line \otimes \pi_{*}(\Line')) = \Big\{(\Line \otimes \pi_{*}(\Line'),  \E_{x'}(\Line(-x)), q+1)\Big\}$ \\
where $\Line'=\mathcal{O}_{X_2}(y-x) \in \mathrm{Pic}^{0}X_2, \; \Line \in \mathrm{Pic}X$ and $x' =   \mathrm{Norm}_{1}^{2}(y) -x.$
\end{enumerate}

The identities between closed points (as throughout the article) are taken as divisors such that the equivalence on $\mathrm{Pic}^0X$ is the trivial line bundle.  
\end{theorem}

The rest of this section is dedicated to the proof of Theorem \ref{theoremcase2}. Part (1) follows from Proposition \ref{neighbousoftotaldec} and  Lemma \ref{lemma2}. Proposition \ref{neighbousoftotaldec}, Lemma \ref{lemma4} and Lemma \ref{lemma3} proof part (2). Part (3) follows from lemmas \ref{lemma1} and \ref{lemma5}. Part (4) follows from Lemma \ref{lemma6}. Lemma \ref{lemma5} and Lemma \ref{lemma6} prove (5). Part (6) follows from Proposition \ref{propstable}, Lemma \ref{lemma2} and Lemma \ref{lemma3}. Finally, (7) follows from Lemma \ref{lemma5}.

\begin{lemma}\label{lemma1} Let $\E = \E_{(x_1,1)}^{(1,d-1)} \oplus \E_{(x_2,1)}^{(1,d)}$ with $x_1 \neq x_2$. Then
\[m_{x,1}\big(\E_{(y_1,1)}^{(1,d)} \oplus \E_{(y_2,1)}^{(1,d)}, \E\big)=1 \quad 
\text{ and } \quad
m_{x,1}\big(\E_{(y_1',1)}^{(1,d)} \oplus \E_{(y_2',1)}^{(1,d)}, \E\big)=q\]
where $y_1 = x + x_1,\; y_2 = x_2$, $y_1'=x_1$ and  $y_2'=x + x_2.$

\begin{proof} In the Hall algebra $\mathbf{H}_X$ of $X$, we have $\E = v \; \E_{(x_1,1)}^{(1,d-1)} \E_{(x_2,1)}^{(1,d)}.$ Via Atiyah's classification, $\E_{(x_1,1)}^{(1,d-1)}$ corresponds to $\mathcal{K}_{x_1}^{}$  and $\E_{(x_2,1)}^{(1,d)}$ corresponds to $\mathcal{K}_{x_2}^{}.$ By Proposition \ref{inversionlemma}, $\mathcal{K}_{x_i}^{} = T_{(0,1),x_i} = N_{1}^{-1} \sum_{\widetilde{\rho}_i \in  \mathrm{P}_1} \widetilde{\rho}_i(-x_i) T_{(0,1)}^{\widetilde{\rho}_i}$, for $i=1,2.$ Thus,
\[ \E_{(x_1,1)}^{(1,d-1)} = N_{1}^{-1} \sum_{\widetilde{\rho}_1 \in   \mathrm{P}_1} \widetilde{\rho}_1(-x_1) T_{(1,d-1)}^{\widetilde{\rho}_1}, \; \E_{(x_2,1)}^{(1,d)} = N_{1}^{-1} \sum_{\widetilde{\rho}_2 \in   \mathrm{P}_1} \widetilde{\rho}_2(-x_2) T_{(1,d)}^{\widetilde{\rho}_2}\]
and 
\[\E = v N_{1}^{-2} \sum_{\widetilde{\rho}_1, \widetilde{\rho}_2 \in\mathrm{P}_1} \widetilde{\rho}_1(-x_1) \widetilde{\rho}_2(-x_2)\; T_{(1,d-1)}^{\widetilde{\rho}_1} T_{(1,d)}^{\widetilde{\rho}_2}. \]
Therefore, 
\begin{eqnarray*} \pi^{\mathrm{vec}}(\mathcal{K}_{x}^{}\; \E)& = &\Big[N_{1}^{-1} \sum_{\widetilde{\rho} \in  \mathrm{P}_1} \widetilde{\rho}(-x) T_{(0,1)}^{\widetilde{\rho}},  v N_{1}^{-2} \sum_{\widetilde{\rho}_1, \widetilde{\rho}_2 \in \mathrm{P}_1} \widetilde{\rho}_1(-x_1) \widetilde{\rho}_2(-x_2)\; T_{(1,d-1)}^{\widetilde{\rho}_1} T_{(1,d)}^{\widetilde{\rho}_2}\Big]\\
& = & v N_{1}^{-3} \sum_{\widetilde{\rho} ,\widetilde{\rho}_1, \widetilde{\rho}_2 \in \mathrm{P}_1} \widetilde{\rho}(-x) \widetilde{\rho}_1(-x_1) \widetilde{\rho}_2(-x_2)\Big[T_{(0,1)}^{\widetilde{\rho}}, \; T_{(1,d-1)}^{\widetilde{\rho}_1} T_{(1,d)}^{\widetilde{\rho}_2} \Big]\\
& = &  v N_{1}^{-3} \Big( \sum_{\widetilde{\rho} \neq \widetilde{\rho}_2} \widetilde{\rho}(-x-x_1) \widetilde{\rho}_2(-x_2) \Big[ T_{(0,1)}^{\widetilde{\rho}}, T_{(1,d-1)}^{\widetilde{\rho}} \Big] T_{(1,d)}^{\widetilde{\rho}_2} \\
& + &   \sum_{\widetilde{\rho} \neq \widetilde{\rho}_1} \widetilde{\rho}(-x-x_2) \widetilde{\rho}_1(-x_1) T_{(1,d-1)}^{\widetilde{\rho}_1}\Big[ T_{(0,1)}^{\widetilde{\rho}}, T_{(1,d)}^{\widetilde{\rho}} \Big] \\
& + & \sum_{\widetilde{\rho}} \widetilde{\rho}(-x-x_1-x_2)  \Big[ T_{(0,1)}^{\widetilde{\rho}}, T_{(1,d-1)}^{\widetilde{\rho}} T_{(1,d)}^{\widetilde{\rho}} \Big] \Big).\\
\end{eqnarray*}
By Pick's formula, $\Delta_{(1,d-1),(0,1)} \cap \mathbf{Z} = \Delta_{(1,d),(0,1)} \cap \mathbf{Z} = \emptyset.$ It follows from Definition \ref{defpathalgebra} and Theorem \ref{dragosmaintheorem} that 
\[\Big[ T_{(0,1)}^{\widetilde{\rho}}, T_{(1,d-1)}^{\widetilde{\rho}} \Big] = c_1 T_{(1,d)}^{\widetilde{\rho}} ,\quad \Big[ T_{(0,1)}^{\widetilde{\rho}}, T_{(1,d)}^{\widetilde{\rho}} \Big] = c_1 T_{(1,d+1)}^{\widetilde{\rho}} \quad \text{ and }\]
\[\Big[ T_{(0,1)}^{\widetilde{\rho}}, T_{(1,d-1)}^{\widetilde{\rho}} T_{(1,d)}^{\widetilde{\rho}} \Big] = c_1 \Big( T_{(1,d)}^{\widetilde{\rho}} T_{(1,d)}^{\widetilde{\rho}} + T_{(1,d-1)}^{\widetilde{\rho}} T_{(1,d+1)}^{\widetilde{\rho}} \Big).\] 
Thus 
\begin{eqnarray*} \pi^{\mathrm{vec}}(\mathcal{K}_{x}^{}\; \E)& = & v c_1 N_{1}^{-3} \Big( \sum_{\widetilde{\rho},\widetilde{\rho}_2 \in \mathrm{P}_1} \widetilde{\rho}(-x-x_1)\widetilde{\rho}_2(-x_2) T_{(1,d)}^{\widetilde{\rho}} T_{(1,d)}^{\widetilde{\rho}_2} \\
& + & \sum_{\widetilde{\rho},\widetilde{\rho}_1 \in \mathrm{P}_1} \widetilde{\rho}(-x-x_2)\widetilde{\rho}_1(-x_1) T_{(1,d-1)}^{\widetilde{\rho}_1} T_{(1,d+1)}^{\widetilde{\rho}} \Big)\\
& = & v c_1 N_{1}^{-3} \Big( N_{1}^{2}\; \E_{(y_1,1)}^{(1,d)} \oplus \E_{(y_2,1)}^{(1,d)} + v^{-2} N_{1}^{2}\; \E_{(y_1',1)}^{(1,d-1)} \oplus \E_{(y_2',1)}^{(1,d+1)} \Big) \\
& = & v^2 \; \E_{(y_1,1)}^{(1,d)} \oplus \E_{(y_2,1)}^{(1,d)} +  \E_{(y_1',1)}^{(1,d-1)} \oplus \E_{(y_2',1)}^{(1,d+1)}
\end{eqnarray*}
where $y_1 = x + x_1,\; y_2 = x_2$, $y_1'=x_1$ and $ y_2'=x + x_2.$ The multiplicities of the edges follows from Remark \ref{propconnection}, namely by multiplication by $v^{-2}$. 
\end{proof}
\end{lemma}

\begin{lemma}\label{lemma2} Let $\E = \E_{(x_1,1)}^{(1,d)} \oplus \E_{(x_2,1)}^{(1,d)}$ and $x_1 \neq x_2.$ Then
\[ m_{x,1}\big(\E_{(x',1)}^{(2,2d+1)}, \E\big)=1, \quad
m_{x,1}\big(\E_{(y_1,1)}^{(1,d/2)} \oplus \E_{(y_2,1)}^{(1,d/2+ 1)}, \E \big)=q \quad \text{and} \] 
\[m_{x,1}\big(\E_{(y_1',1)}^{(1,d/2)} \oplus \E_{(y_2',1)}^{(1,d/2+1)}, \E \big)=q\]
where 
\[ x' =  x+x_1+x_2, \quad y_1 =x_2,\quad y_2 = x+x_1, \quad  y_1'=x_1 \quad \text{ and} \quad y_2'=x+x_2.\]

\begin{proof} Since $x_1 \neq x_2$, $\E = N_{1}^{-2} \sum_{\widetilde{\rho}_1,\widetilde{\rho}_2 \in \mathrm{P}_1} \widetilde{\rho}_1(-x_1) \widetilde{\rho}(-x_2) T_{(1,d)}^{\widetilde{\rho}_1} T_{(1,d)}^{\widetilde{\rho}_2}.$ Thus
\begin{eqnarray*} \pi^{\mathrm{vec}}(\mathcal{K}_{x}^{}\; \E)& = & N_{1}^{-3} \sum_{\widetilde{\rho},\widetilde{\rho}_1\widetilde{\rho}_2 \in \mathrm{P}_1} \widetilde{\rho}(-x)\widetilde{\rho}_1(-x_1)\widetilde{\rho}_2(-x_2) \Big[T_{(0,1)}^{\widetilde{\rho}} ,T_{(1,d)}^{\widetilde{\rho}_1} T_{(1,d)}^{\widetilde{\rho}_2} \Big] \\
\end{eqnarray*}
From $ \Big[T_{(0,1)}^{\widetilde{\rho}} ,T_{(1,d)}^{\widetilde{\rho}} T_{(1,d)}^{\widetilde{\rho}} \Big] = c_{1}^{2} \;T_{(2,d+1)}^{\widetilde{\rho}} + \;2 c_1 \;T_{(1,d)}^{\widetilde{\rho}} T_{(1,d+1)}^{\widetilde{\rho}}$ and \\$ \Big[T_{(0,1)}^{\widetilde{\rho}} ,T_{(1,d)}^{\widetilde{\rho}} \Big] = c_1 T_{(1,d+1)}^{\widetilde{\rho}}$, it follows that 
\begin{eqnarray*}  \pi^{\mathrm{vec}}(\mathcal{K}_{x}^{}\; \E)& = & c_1  N_{1}^{-3}  \sum_{\widetilde{\rho}\neq \widetilde{\rho}_2} \widetilde{\rho}(-x-x_1)\widetilde{\rho}_2(-x_2) T_{(1,d)}^{\widetilde{\rho}_2} T_{(1,d+1)}^{\widetilde{\rho}} \\
& + & c_1 N_{1}^{-3} \sum_{\widetilde{\rho}\neq \widetilde{\rho}_1} \widetilde{\rho}(-x-x_2)\widetilde{\rho}_1(-x_1) T_{(1,d)}^{\widetilde{\rho}_1} T_{(1,d+1)}^{\widetilde{\rho}} \\
& + & c_1 N_{1}^{-3} \sum_{\widetilde{\rho}}\widetilde{\rho}(-x-x_1-x_2) \big(c_1 T_{(2,2d+1)}^{\widetilde{\rho}} + 2 T_{(1,d)}^{\widetilde{\rho}} T_{(1,d+1)}^{\widetilde{\rho}} \big)  \\
& = & \E_{(y_1,1)}^{(1,d/2)} \oplus \E_{(y_2,1)}^{(1,d/2+ 1)} + \E_{(y_1',1)}^{(1,d/2)} \oplus \E_{(y_2',1)}^{(1,d/2+1)} + v^{2} \E_{(x',1)}^{(2,2d+1)}
\end{eqnarray*}
where $x' =  x+x_1+x_2$, $y_1 =x_2,\; y_2 = x+x_1$  $y_1'=x_1,$ and $ y_2'=x+x_2.$
\end{proof}
\end{lemma}

\begin{lemma} \label{lemma3}Let $\E = \E_{(x',2)}^{(2,2d)}$ with $|x'|=1.$ Then
\[m_{x,1}\big(\E_{(y_1,1)}^{(1,d)} \oplus \E_{(y_2,1)}^{(1,d+1)}, \E \big)=q-1 \quad 
\text{and} \quad 
m_{x,1}\big(\E_{(y',1)}^{(2,2d+1)}, \E \big)=1\]
where $y_1 = x', y_2 = x+x'$ and $y'-x = \mathrm{Norm}_{1}^{2}(x'') = 2x'$, with $x'' \in |X_2|$ sits above $x'.$ 

\begin{proof} Observe that $\E$ corresponds to $\mathcal{K}_{x'}^{(2)}$ via Atiyah's Theorem \ref{atiyahtheorem}, and by definition $\mathcal{K}_{x'}^{(2)} = \frac{2}{[2]} T_{(0,2),x'} - (1-v^{-2}) \mathcal{K}_{x}^{\oplus 2}.$ By Newton's formula,
\[\mathcal{K}_{x'}^{\oplus 2} = \frac{v^2}{2}\Big( T_{(0,1),x'}^{2} - \frac{2}{[2]} T_{(0,2),x'}\Big)\]
and thus
\[\E = \ell_1 \; T_{(2,2d),x'} - \ell_2\; T_{(1,d),x'}^{2}\]
where $\ell_1 := \frac{2}{[2]} + (1-v^{-2}) \frac{v^{2}}{[2]}$ and $\ell_2 := (1-v^{-2}) \frac{v^2}{2}.$ Hence
\begin{eqnarray*}
\pi^{\mathrm{vec}}(\mathcal{K}_{x} \; \E) & = & \ell_1 N_{2}^{-1} N_{1}^{-1} \sum_{\substack{\widetilde{\rho} \in \mathrm{P}_1 \\ \widetilde{\chi} \in \mathrm{P}_2}} \widetilde{\rho}(-x)\widetilde{\chi}(-x') \Big[T_{(0,1)}^{\widetilde{\rho}}, T_{(2,2d)}^{\widetilde{\chi}} \Big] \\
& - & \ell_2 N_{1}^{-3} \sum_{\widetilde{\rho},\widetilde{\rho}_1, \widetilde{\rho}_2 \in \mathrm{P}_1} \widetilde{\rho}(-x)\widetilde{\rho}_1(-x')\widetilde{\rho}_2(-x') \Big[ T_{(0,1)}^{\widetilde{\rho}}, T_{(1,d)}^{\widetilde{\rho}_1} T_{(1,d)}^{\widetilde{\rho}_2} \Big].
\end{eqnarray*}
If $\widetilde{\chi} \neq \mathrm{Norm}_{1}^{2}\widetilde{\rho}$, then $\Big[T_{(0,1)}^{\widetilde{\rho}}, T_{(2,2d)}^{\widetilde{\chi}} \Big]=0$ by Theorem \ref{dragosmaintheorem}. Otherwise 
\[\Big[T_{(0,1)}^{\widetilde{\rho}}, T_{(2,2d)}^{\widetilde{\chi}} \Big] = c_2 T_{(2,2d+1)}^{\widetilde{\rho}}.\]
By the above observation and a calculation similar to that in the proof of Lemma \ref{lemma1}, we have 
\begin{eqnarray*}
\pi^{\mathrm{vec}}(\mathcal{K}_{x} \; \E) & = & \ell_1 v^2 \frac{[2]}{2} \;\E_{(y',1)}^{(2,2d+1)} 
 - 2 \ell_2 \;\E_{(y_1,1)}^{(1,d)} \oplus \E_{(y_2,1)}^{(1,d+1)} - \ell_2 v^{2}\; \E_{(y',1)}^{(2,2d+1)} \\
& = & v^2 \; \E_{(y',1)}^{(2,2d+1)} + (q-1)v^{2} \;\E_{(y_1,1)}^{(1,d)} \oplus \E_{(y_2,1)}^{(1,d+1)}
\end{eqnarray*}
where $y_1 = x', y_2 = x+x'$ and $y'-x=\mathrm{Norm}_{1}^{2}(x'')=2x'$, with $x'' \in |X_2|$ sits above $x'.$ 
\end{proof}
\end{lemma}

\begin{lemma} \label{lemma4}Let $\E = \E_{(x',1)}^{(1,d)} \oplus \E_{(x',1)}^{(1,d)},$ then
$m_{x,1}\big(\E_{(y_1,1)}^{(1,d)} \oplus\E_{(y_2,1)}^{(1,d+1)}  , \E\big)=1$
where $y_1 = x,\; y_2=x+x'.$

\begin{proof}By Atiyah's classification and Newton's formulas, we have 
$$\E ={\tiny \frac{v^2}{2}}\Big( N_{1}^{-2} \big(\sum_{\widetilde{\rho}_1 \in \mathrm{P}_1} \widetilde{\rho}_1(-x') T_{(1,d)}^{\widetilde{\rho}_1} \big)^{2} - {\tiny\frac{2}{[2]}} N_{2}^{-1} \sum_{ \widetilde{\chi} \in \mathrm{P}_2} \widetilde{\chi}(-x') T_{(2,2d)}^{\widetilde{\chi}} \Big).$$
Thus, 
\begin{eqnarray*}
\pi^{\mathrm{vec}}(\mathcal{K}_{x} \; \E) & = & \frac{v^2}{2} N_{1}^{-3} \sum_{\widetilde{\rho},\widetilde{\rho}_1 \in \mathrm{P}_1} \widetilde{\rho}(-x)\widetilde{\rho}_1(-2x') \Big[T_{(0,1)}^{\widetilde{\rho}}, T_{(1,d)}^{\widetilde{\rho}_1} T_{(1,d)}^{\widetilde{\rho}_1} \Big]\\
& + & v^{2} N_{1}^{-3} \sum_{\substack{\widetilde{\rho},\widetilde{\rho}_1,\widetilde{\rho}_2 \in \mathrm{P}_1 \\ \widetilde{\rho}_1 \neq \widetilde{\rho}_2}} \widetilde{\rho}(-x)\widetilde{\rho}_1(-x')\widetilde{\rho}_2(-x') \Big[ T_{(0,1)}^{\widetilde{\rho}}, T_{(1,d)}^{\widetilde{\rho}_1} T_{(1,d)}^{\widetilde{\rho}_2} \Big] \\
& - & \frac{v^2}{[2]} N_{1}^{-1} N_{2}^{-1} \sum_{\substack{ \widetilde{\rho} \in \mathrm{P}_1 \\ \widetilde{\chi} \in \mathrm{P}_2}} \widetilde{\rho}(-x)\widetilde{\chi}(-x') \Big[T_{(0,1)}^{\widetilde{\rho}}, T_{(2,2d)}^{\widetilde{\chi}} \Big].
\end{eqnarray*}
We have $\Big[T_{(0,1)}^{\widetilde{\rho}}, T_{(2,2d)}^{\widetilde{\chi}} \Big] =0$ unless $\widetilde{\chi} = \mathrm{Norm}_{1}^{2}(\widetilde{\rho}),$ in such case $\Big[T_{(0,1)}^{\widetilde{\rho}}, T_{(2,2d)}^{\widetilde{\chi}} \Big] = c_2 T_{(2,2d+1)}^{\widetilde{\rho}}$ (see Definition \ref{defpathalgebra} and Theorem \ref{dragosmaintheorem}). Therefore
\[ \pi^{\mathrm{vec}}(\mathcal{K}_{x} \; \E)  = v^{2} \; \E_{(y_1,1)}^{(1,d)} \oplus\E_{(y_2,1)}^{(1,d+1)}\]
where $y_1 = x,\; y_2=x+x'.$
\end{proof}
\end{lemma}

\begin{lemma} \label{lemma5}Let $\E= \E_{(x',1)}^{(2,d)}$ with $\gcd(2,d)=1.$ Then 
\[ m_{x,1}(\E_{(y',2)}^{(2,d+1)},  \E)=q, \quad 
m_{x,1}\big(\E_{(y,1)}^{(2,d+1)}, \E\big)=q+1, \] 
\[\text{and} \quad m_{x,1}\big(\E_{(y_1,1)}^{(1,(d+1)/2)} \oplus \E_{(y_2,1)}^{(1,(d+1)/2)}, \E \big)=q-1  \]
where $|y|=2,$ $\mathrm{Norm}_{1}^{2}(z) =  x+ x'$ with $z \in |X_2|$ sits above $y$. And $y_1 \neq y_2,\; y_1+ y_2 = x+x'$ and $|y'|=1,$ $2 y'=x+x'.$

\begin{proof} In the Hall algebra $\mathbf{H}_X$, $\E = N_{1}^{-1} \sum_{\widetilde{\rho}_1 \in \mathrm{P}_1}\widetilde{\rho}_1(-x') T_{(2,d)}^{\widetilde{\rho}_1},$ thus
\begin{eqnarray*}
\pi^{\mathrm{vec}}(\mathcal{K}_{x} \; \E) & = & \Big[ N_{1}^{-1} \sum_{\widetilde{\rho} \in \mathrm{P}_1} \widetilde{\rho}(-x) T_{(0,1)}^{\widetilde{\rho}}, \; N_{1}^{-1} \sum_{\widetilde{\rho}_1 \in \mathrm{P}_1}\widetilde{\rho}_1(-x') T_{(2,d)}^{\widetilde{\rho}_1} \Big] \\
& = & N_{1}^{-2} \sum_{\widetilde{\rho},\widetilde{\rho}_1 \in \mathrm{P}_1} \widetilde{\rho}(-x)\widetilde{\rho}_1(-x') \Big[ T_{(0,1)}^{\widetilde{\rho}}, \; T_{(2,d)}^{\widetilde{\rho}_1} \Big] \\
& = & N_{1}^{-2} \sum_{\widetilde{\rho} \in \mathrm{P}_1} \widetilde{\rho}(-x-x') \Big[ T_{(0,1)}^{\widetilde{\rho}}, \; T_{(2,d)}^{\widetilde{\rho}} \Big] 
\end{eqnarray*}
By Theorem \ref{dragosmaintheorem}, Pick's formula and Definition \ref{defpathalgebra}
\[ \Big[ T_{(0,1)}^{\widetilde{\rho}}, \; T_{(2,d)}^{\widetilde{\rho}} \Big]  = c_1 \;\frac{\theta_{(2,d+1)}}{v^{-1}-v},\]
where $\theta_{(2,d+1)} = (v^{-1}-v)\; T_{(2,d+1)}^{\mathrm{Norm}_{1}^{2}\widetilde{\rho}} + \frac{(v^{-1}-v)^{2}}{2} \big( T_{(1,(d+1)/2)}^{\widetilde{\rho}} \big)^{2}.$ 
Hence
\begin{eqnarray*}
\pi^{\mathrm{vec}}(\mathcal{K}_{x} \; \E) & = & c_1 \; N_{1}^{-2} \Big(  \sum_{ \tiny\substack{ y' \in |X| \\
\widetilde{\rho} \in \mathrm{P}_1} } \widetilde{\rho}(\mathrm{Norm}_{1}^{2}(y'')-x-x') T_{(2,d+1),y'} \\
& + & \frac{\tiny (v^{-1}-v)}{2} \sum_{\tiny\substack{y_1, y_2 \in |X| \\ \widetilde{\rho} \in \mathrm{P}_1}} \widetilde{\rho}(y_1+y_2-x-x') T_{(1,(d+1)/2), y_1} \;T_{(1,(d+1)/2),y_2} \Big).
\end{eqnarray*}
where $y'' \in |X_2|$ sits above $y'.$

Let us calculate which term separated. Since $T_{(2,d+1),y'}= [2]\E_{(y',1)}^{(2,d+1)}$ if $|y'|=2$ and 
$T_{(2,d+1),y'}= \frac{[2]}{2} \Big((1-v^{-2})\; \E_{(y',1)}^{(1,(d+1)/2)} \oplus \E_{(y',1)}^{(1,(d+1)/2)} + \E_{(y',2)}^{(2,d+1)}\Big)$ if $|y'|=1$, the first term yields
\[\sum_{\tiny\substack{y \in |X|, |y|=2 \\ \mathrm{Norm}_{1}^{2}(z)=x+x'}}(v^2 +1) \E_{(y,1)}^{(2,d+1)} + \sum_{\tiny \substack{ y' \in |X|, \;|y'|=1 \\ \mathrm{Nomr}_{1}^{2}(z') = x+x'}} \frac{v^2 +1}{2} \Big( (1-v^{-2}) \E_{(y',1)}^{(1,\frac{d+1}{2})} \oplus \E_{(y',1)}^{(1,\frac{d+1}{2})} + \E_{(y',2)}^{(2,d+1)} \Big),\]
where $z,z' \in |X_2|,$ $z$ sits above $y$ and $z'$ sits above $y'.$

About the second term, 
\begin{eqnarray*}
& &c_1 \; N_{1}^{-2} \; \frac{\tiny (v^{-1}-v)}{2} \sum_{\tiny\substack{y_1, y_2 \in X \\ \widetilde{\rho} \in \mathrm{P}_1}} \widetilde{\rho}(y_1+y_2-x-x') T_{(1,(d+1)/2), y_1} \;T_{(1,(d+1)/2),y_2} \\
& = & (1-v^{2})\sum_{\tiny\substack{y_1 \neq y_2 \\ y_1 + y_2 = x+x'}} \E_{(y_1,1)}^{(1,\frac{d+1}{2})} \oplus \E_{(y_2,1)}^{(1,\frac{d+1}{2})} + \frac{1-v^2}{2} \sum_{\tiny\substack{ y' \in X \\ 2y'=x+x'}} \E_{(y',1)}^{(1,\frac{d+1}{2})} \; \E_{(y',1)}^{(1,\frac{d+1}{2})}. 
\end{eqnarray*}
Rest us calculate the product $\E_{(y',1)}^{(1,\frac{d+1}{2})}\;\E_{(y',1)}^{(1,\frac{d+1}{2})}.$ Via the equivalence $\mathsf{C}_{\frac{d+1}{2}} \equiv \mathsf{C}_{\infty}$ and the isomorphism $\mathbf{H}_{\mathrm{Tor}_{y'}} \cong \Lambda_{v^{2|x|}}$, the previous product corresponds to the product $e_1 e_1$ of elementary symmetric function in $\Lambda_{v^{2|x|}}.$ By Newton's formula 
\[e_1 e_1 = 2 e_2 + p_2 = (v^{-2}+1) e_2 + P_{(2)}\]
and  $\mathcal{K}_{y'}^{}\mathcal{K}_{y'}^{} = (v^{-2}+1)\mathcal{K}_{y'}^{\oplus 2} + \mathcal{K}_{y'}^{(2)}.$ Thus,
\[\E_{(y',1)}^{(1,\frac{d+1}{2})}\;\E_{(y',1)}^{(1,\frac{d+1}{2})} = (v^{-2}+1)\; \E_{(y',1)}^{(1,\frac{d+1}{2})}\oplus \E_{(y',1)}^{(1,\frac{d+1}{2})} + \E_{(y',2)}^{(2,d+1)}.\] 
Therefore the second term yields
\[(1-v^{2})\sum_{\tiny\substack{y_1 \neq y_2 \\ y_1 + y_2 = x+x'}} \E_{(y_1,1)}^{(1,\frac{d+1}{2})} \oplus \E_{(y_2,1)}^{(1,\frac{d+1}{2})} + \frac{v^{-2}-v^{2}}{2} \sum_{\tiny\substack{ y' \in |X| \\ 2y'=x+x'}}\E_{(y',1)}^{(1,\frac{d+1}{2})}\oplus \E_{(y',1)}^{(1,\frac{d+1}{2})} + \frac{1-v^2}{2}  \sum_{\tiny\substack{ y' \in |X|  \\ 2y'=x+x'}}\E_{(y',2)}^{(2,d+1)}.  \]

Putting all together follows the proof. 
\end{proof}
\end{lemma}

\begin{lemma} \label{lemma6} Let $\E = \E_{(x',1)}^{(1,d-1)} \oplus \E_{(x',1)}^{(1,d)}$. Then
$$ m_{x,1}\big(\E_{(x',1)}^{(1,d)} \oplus \E_{(x',1)}^{(1,d)}, \E\big)=q+1, \quad\;
m_{x,1}\big(\E_{(x',2)}^{(2,2d)}, \E\big)=1 \quad $$
$$\text{ and } \quad m_{x,1}\big(\E_{(x',1)}^{(1,d-1)} \oplus \E_{(x',1)}^{(1,d+1)}, \E\big)=q.$$

\begin{proof} By a similar calculation as in the proof of Lemma \ref{lemma1}, we have
\[\pi^{\mathrm{vec}}(\mathcal{K}_{x}^{} \;\E) = v c_1 N_{1}^{-1}\Big( T_{(1,d),x'} T_{(1,d),x'} + T_{(1,d-1),x'} T_{(1,d+1),x'}\Big)\]
since we are considering $x$ as the neutral element of $X(\Fq).$ By Proposition \ref{prop-macdonald}, $T_{(1,d),x'}$ corresponds to the elementary symmetric function $e_1,$ thus the product $T_{(1,d),x'} T_{(1,d),x'} $ corresponds to the product $e_1 e_1$, which  can be written as $e_1 e_1 = 2 e_2 + p_2$ by Newton's formula. Hence
$$T_{(1,d),x'} T_{(1,d),x'}  = (v^{-2}+1)\; \E_{(x',1)}^{(1,d)} \oplus \E_{(x',1)}^{(1,d)} + \E_{(x',2)}^{(2,2d)}$$   
and therefore 
\[ \pi^{\mathrm{vec}}(\mathcal{K}_{x}^{} \;\E) = (1+ v^{2})\;\E_{(x',1)}^{(1,d)} \oplus \E_{(x',1)}^{(1,d)} + v^{2}\; \E_{(x',2)}^{(2,2d)} + \E_{(x',1)}^{(1,d-1)} \oplus \E_{(x',1)}^{(1,d+1)}. \]
The lemma follows after multiplying with $v^{-2}.$
\end{proof}
\end{lemma}

This concludes the proof of Theorem \ref{theoremcase2}. 
We end this section with an explicit example. 

\begin{example}  The easiest example are given by elliptic curves with only one rational point. 
There are up to a isomorphism three such elliptic curves (see for example \cite{serre} 2.4.4 and Ex. 3 of 2.4): $X_2$ over $\mathbb{F}_2$ defined by the Weierstrass equation $y^2 + y = x^3 + x + 1,$ $X_{3}$ over $\mathbb{F}_3$ defined by the Weierstrass equation $y^2 = x^3 + 2x +2 $ and $X_{4}$ over $\mathbb{F}_4$ defined by the Weierstrass equation $y^2 + y = x^3 + \alpha$ where $\mathbb{F}_4 = \mathbb{F}_2(\alpha)$. Observe that in this example $X_n$ does not denote $X \times_{\mathrm{Spec}\; \Fq} \mathrm{Spec}\; \mathbb{F}_{q^n}.$ Since the class number is $1$, 
$$\mathrm{Bun}_{2}^{\mathrm{dec}} X_{q} = \{ \mathcal{O}_{X_{q}}(n) \oplus \mathcal{O}_{X_{q}}(m) \; |\; n,m \in \Z ,\; n \geq m\}$$ and  
$$\mathrm{Bun}_{2}^{\mathrm{gi}}X_q = \{ \E_x(\mathcal{O}_{X_{q}}(n)), \; \E(\mathcal{O}_{X_{q}}(n)) \; | \; n \in \Z \}$$ for $q \in \{ 2,3,4\}.$ One calculates that $\mathrm{Pic}^{0}\big( X_{2} \otimes \mathbb{F}_4 \big) \cong \Z/5\Z$, $\mathrm{Pic}^{0}\big( X_{3} \otimes \mathbb{F}_9\big) \cong \Z/7\Z$ and $\mathrm{Pic}^{0} \big( X_{4} \otimes \mathbb{F}_{16}\big) \cong \Z /9\Z,$ thus 
$$\mathrm{Bun}_{2}^{\mathrm{tr}}X_{q} = \{ \mathcal{O}_{X_{q}}(n)\pi_{*}(\Line_{1}'), \cdots, \mathcal{O}_{X_{q}}(n)\pi_{*}(\Line_{q}') \; | \; n \in \Z \}$$ 
where  $\mathrm{Pic}^{0}\big( X_{q} \otimes \mathbb{F}_{q^2} \big)= \{\Line_{0}' , \ldots, \Line_{2q}'\}.$ By Theorem \ref{theoremcase2} we obtain the following figure.    


\begin{sidewaysfigure}

\centering

\vspace{13cm}


$
  \beginpgfgraphicnamed{tikz/figlast}
     \begin{tikzpicture}[ >=latex, scale=1.8]
    
        \vertex[circle,fill,label={left:$$}] (t1n) at (0,2) {}; 
        \vertex[circle,fill,label={left:$$}](t2n) at (0,1) {}; 

  \vertex[circle,fill,label={left:$$}] at (0,0.5) {}; 
        \vertex[circle,fill,label={left:$$}] at (0,-0.5) {}; 

  \vertex[circle,fill,label={left:$$}] at (0,-1) {}; 
        \vertex[circle,fill,label={left:$$}] at (0,-2) {} ;

        { traces};

 \path[-,font=\scriptsize]
(2,0.75) edge[->-=0.8] node[pos=0.2,auto,black,swap] {\tiny $1$} (0,2)
(2,0.75) edge[->-=0.8] node[pos=0.3,auto,black, swap] {\tiny  $1$} (0,1);

 \path[-,font=\scriptsize]
(2,-0.75) edge[->-=0.8] node[pos=0.2,auto,black,swap] {\tiny $1$} (0,0.5)
(2,-0.75) edge[->-=0.8] node[pos=0.3,auto,black, swap] {\tiny  $1$} (0,-0.5);

 \path[-,font=\scriptsize]
(2,-2.25) edge[->-=0.8] node[pos=0.2,auto,black,swap] {\tiny $1$} (0,-1) 
 (2,-2.25) edge[->-=0.8] node[pos=0.3,auto,black, swap] {\tiny  $1$} (0,-2);

     \vertex[circle,fill,label={above:$$}] (sn-1) at (2,2.25) {};
        \vertex[circle,fill,label={above:$$}] at (2,0.75) {}; 
          \vertex[circle,fill,label={above:$$}] at (2,-0.75) {}; 
        \vertex[circle,fill,label={below:$$}](11) at (2,-2.25) {};

 \path[-,font=\scriptsize]
(0,2) edge[->-=0.8] node[pos=0.3, auto,black] {\tiny $q+1$} (2,2.25)
 (0,1) edge[->-=0.8] node[ pos=0.3,auto,black] {\tiny $q+1$} (2,2.25);

\path[-,font=\scriptsize]
 (0,0.5) edge[->-=0.8] node[pos=0.3,auto,black] {\tiny $q+1$} (2,0.75)
(0,-0.5) edge[->-=0.8] node[pos=0.3,auto,black] {\tiny  $q+1$} (2,0.75);

 \path[-,font=\scriptsize]
(0,-1) edge[->-=0.8] node[pos=0.3,auto,black] {\tiny $q+1$} (2,-0.75)
(0,-2) edge[->-=0.8] node[pos=0.3,auto,black] {\tiny  $q+1$} (2,-0.75);

 \path[-,font=\scriptsize]
(2,-2.25) edge[->-=0.8] node[pos=0.2,auto,black] {\tiny $1$} (4,-1.5)
(4,-1.5) edge[->-=0.8] node[pos=0.2,auto,black,swap] {\tiny  $q$} (2,-0.75);
 \path[-,font=\scriptsize]
(2,-0.75) edge[->-=0.8] node[pos=0.2,auto,black] {\tiny $1$} (4,0)
(4,0) edge[->-=0.8] node[pos=0.2,auto, black,swap] {\tiny  $q$} (2,0.75);
 \path[-,font=\scriptsize]
(2,0.75) edge[->-=0.8] node[pos=0.2,auto,black] {\tiny $1$} (4,1.5)
(4,1.5) edge[->-=0.8] node[pos=0.2,auto,black,swap] {\tiny  $q$} (2,2.25);

        \vertex[circle,fill,label={above:$$}] at (4,1.5) {};
          \vertex[circle,fill,label={above:$$}] at (4,0) {};
        \vertex[circle,fill,label={below:$$}](11) at (4,-1.5) {};

        \vertex[circle,fill,label={left:$$}] at (5,1.5) {};
          \vertex[circle,fill,label={left:$$}] at (5,0) {};
        \vertex[circle,fill,label={left:$$}](11) at (5,-1.5) {}; 
        

   \vertex[circle,fill,label={left:$$}] at (6,2.25) {};
        \vertex[circle,fill,label={left:$$}] at (6,0.75) {};
          \vertex[circle,fill,label={left:$$}] at (6,-0.75) {};
        \vertex[circle,fill,label={left:$$}](11) at (6,-2.25) {};

 \path[-,font=\scriptsize]
(6,-2.25) edge[->-=0.8] node[pos=0.3,auto,black,swap] {\tiny $1$} (5,-1.5)
(5,-1.5) edge[->-=0.8] node[below=24pt ,pos=0.6,auto,black] {\tiny  $q+1$} (6,-0.75);

 \path[-,font=\scriptsize]
(6,-0.75) edge[->-=0.8] node[pos=0.3,auto,black,swap] {\tiny $1$} (5,0)
(5,0) edge[->-=0.8] node[below=24pt ,pos=0.6,auto,black] {\tiny  $q+1$} (6,0.75);

 \path[-,font=\scriptsize]
(6,0.75) edge[->-=0.8] node[pos=0.3,auto,black,swap] {\tiny $1$} (5,1.5)
(5,1.5) edge[->-=0.8] node[below=24pt ,pos=0.6,auto,black] {\tiny  $q+1$} (6,2.25);

 \path[-,font=\scriptsize]
(6,-2.25) edge[->-=0.8] node[below=18pt, pos=0.4,auto,black,swap] {\tiny $q-1$} (4,-1.5)
(4,-1.5) edge[->-=0.8] node[pos=0.2,auto,black] {\tiny  $1$} (6,-0.75);
 \path[-,font=\scriptsize]
(6,-0.75) edge[->-=0.8] node[below=18pt, pos=0.4,auto,black,swap] {\tiny $q-1$} (4,0)
(4,0) edge[->-=0.8] node[pos=0.2,auto,black] {\tiny  $1$} (6,0.75);
 \path[-,font=\scriptsize]
(6,0.75) edge[->-=0.8] node[below=18pt,pos=0.4,auto,black,swap] {\tiny $q-1$} (4,1.5)
(4,1.5) edge[->-=0.8] node[pos=0.2,auto,black] {\tiny  $1$} (6,2.25);      

        \vertex[circle,fill,label={left:$$}] at (8,1.5) {};
          \vertex[circle,fill,label={left:$$}] at (8,0) {};
        \vertex[circle,fill,label={left:$$}](11) at (8,-1.5) {}; 

 \vertex[circle,fill,label={left:$$}] at (10,2.25) {};
        \vertex[circle,fill,label={left:$$}] at (10,0.75) {};
          \vertex[circle,fill,label={left:$$}] at (10,-0.75) {};
        \vertex[circle,fill,label={left:$$}](11) at (10,-2.25) {};

 \path[-,font=\scriptsize]
(10,-2.25) edge[->-=0.8] node[pos=0.2,auto,black, swap] {\tiny $q$} (8,-1.5)
(8,-1.5) edge[->-=0.8] node[pos=0.2,auto,black] {\tiny  $1$} (10,-0.75)
(10,-0.75) edge[->-=0.8] node[pos=0.2,auto,black, swap] {\tiny  $q$} (8,0)
(8,0) edge[->-=0.8] node[pos=0.2,auto,black] {\tiny  $1$} (10,0.75)
(10,0.75) edge[->-=0.8] node[pos=0.2,auto,black,swap] {\tiny  $q$} (8,1.5)
(8,1.5) edge[->-=0.8] node[pos=0.2,auto,black] {\tiny  $1$} (10,2.25)
;

 \path[-,font=\scriptsize]
(6,-2.25) edge[->-=0.8] node[pos=0.2,auto,black] {\tiny  $1$} (8,-1.5);

 \path[-,font=\scriptsize]
(8,-1.5) edge[->-=0.8] node[pos=0.2,auto,black,swap] {\tiny $q$} (6,-0.75)
(6,-0.75) edge[->-=0.8] node[pos=0.2,auto,black] {\tiny  $1$} (8,0);
 \path[-,font=\scriptsize]
(8,0) edge[->-=0.8] node[pos=0.2,auto,black,swap] {\tiny $q$} (6,0.75)
(6,0.75) edge[->-=0.8] node[pos=0.2,auto,black] {\tiny  $1$} (8,1.5);
 \path[-,font=\scriptsize]
(8,1.5) edge[->-=0.8] node[pos=0.2,auto,black,swap] {\tiny $q$} (6,2.25);

\draw (0,1.4) circle (0.015cm); 
\fill  (0,1.4) circle (0.015cm);      
\draw (0,1.5) circle (0.015cm);
\fill  (0,1.5) circle (0.015cm);  
\draw (0,1.6) circle (0.015cm);
\fill  (0,1.6) circle (0.015cm);

\draw (0,0.1) circle (0.015cm); 
\fill  (0,0.1) circle (0.015cm);      
\draw (0,0) circle (0.015cm);
\fill  (0,0) circle (0.015cm);  
\draw (0,-0.1) circle (0.015cm);
\fill  (0,-0.1) circle (0.015cm);

\draw (0,-1.4) circle (0.015cm); 
\fill  (0,-1.4) circle (0.015cm);      
\draw (0,-1.5) circle (0.015cm);
\fill  (0,-1.5) circle (0.015cm);  
\draw (0,-1.6) circle (0.015cm);
\fill  (0,-1.6) circle (0.015cm);

\draw (0,2.5) circle (0.015cm); 
\fill  (0,2.5) circle (0.015cm);      
\draw (0,2.7) circle (0.015cm);
\fill  (0,2.7) circle (0.015cm);  
\draw (0,2.9) circle (0.015cm);
\fill  (0,2.9) circle (0.015cm);

\draw (2,2.5) circle (0.015cm); 
\fill  (2,2.5) circle (0.015cm);      
\draw (2,2.7) circle (0.015cm);
\fill  (2,2.7) circle (0.015cm);  
\draw (2,2.9) circle (0.015cm);
\fill  (2,2.9) circle (0.015cm); 

\draw (4,2.5) circle (0.015cm); 
\fill  (4,2.5) circle (0.015cm);      
\draw (4,2.7) circle (0.015cm);
\fill  (4,2.7) circle (0.015cm);  
\draw (4,2.9) circle (0.015cm);
\fill  (4,2.9) circle (0.015cm); 

\draw (6,2.5) circle (0.015cm); 
\fill  (6,2.5) circle (0.015cm);      
\draw (6,2.7) circle (0.015cm);
\fill  (6,2.7) circle (0.015cm);  
\draw (6,2.9) circle (0.015cm);
\fill  (6,2.9) circle (0.015cm);

\draw (8,2.5) circle (0.015cm); 
\fill  (8,2.5) circle (0.015cm);      
\draw (8,2.7) circle (0.015cm);
\fill  (8,2.7) circle (0.015cm);  
\draw (8,2.9) circle (0.015cm);
\fill  (8,2.9) circle (0.015cm);

\draw (10,2.5) circle (0.015cm); 
\fill  (10,2.5) circle (0.015cm);      
\draw (10,2.7) circle (0.015cm);
\fill  (10,2.7) circle (0.015cm);  
\draw (10,2.9) circle (0.015cm);
\fill  (10,2.9) circle (0.015cm);

\draw (10.3,2.25) circle (0.015cm); 
\fill  (10.3,2.25) circle (0.015cm);      
\draw (10.5,2.25) circle (0.015cm);
\fill  (10.5,2.25) circle (0.015cm);  
\draw (10.7,2.25) circle (0.015cm);
\fill  (10.7,2.25) circle (0.015cm);  

\draw (10.3,0.75) circle (0.015cm); 
\fill  (10.3,0.75) circle (0.015cm);      
\draw (10.5,0.75) circle (0.015cm);
\fill  (10.5,0.75) circle (0.015cm);  
\draw (10.7,0.75) circle (0.015cm);
\fill  (10.7,0.75) circle (0.015cm);

 \draw (10.3,-0.75) circle (0.015cm); 
\fill  (10.3,-0.75) circle (0.015cm);      
\draw (10.5,-0.75) circle (0.015cm);
\fill  (10.5,-0.75) circle (0.015cm);  
\draw (10.7,-0.75) circle (0.015cm);
\fill  (10.7,-0.75) circle (0.015cm);

\draw (10.3,-2.25) circle (0.015cm); 
\fill  (10.3,-2.25) circle (0.015cm);      
\draw (10.5,-2.25) circle (0.015cm);
\fill  (10.5,-2.25) circle (0.015cm);  
\draw (10.7,-2.25) circle (0.015cm);
\fill  (10.7,-2.25) circle (0.015cm);

\draw (0,-2.5) circle (0.015cm); 
\fill  (0,-2.5) circle (0.015cm);      
\draw (0,-2.7) circle (0.015cm);
\fill  (0,-2.7) circle (0.015cm);  
\draw (0,-2.9) circle (0.015cm);
\fill  (0,-2.9) circle (0.015cm) ;  


\draw  [underbrace style] (-0.25,-3) -- (0.25,-3) node [underbrace text style] {traces};

\draw  [underbrace style] (1.8,-3) -- (4.2,-3) node [underbrace text style] {geometrically indecomposable};

\draw [underbrace style] (4.8,-3) -- (10.2,-3) node [underbrace text style] {sum of two lines bundles};

\draw (2,-2.5) circle (0.015cm); 
\fill  (2,-2.5) circle (0.015cm);      
\draw (2,-2.7) circle (0.015cm);
\fill  (2,-2.7) circle (0.015cm);  
\draw (2,-2.9) circle (0.015cm);
\fill  (2,-2.9) circle (0.015cm); 

\draw (4,-2.5) circle (0.015cm); 
\fill  (4,-2.5) circle (0.015cm);      
\draw (4,-2.7) circle (0.015cm);
\fill  (4,-2.7) circle (0.015cm);  
\draw (4,-2.9) circle (0.015cm);
\fill  (4,-2.9) circle (0.015cm); 

\draw (6,-2.5) circle (0.015cm); 
\fill  (6,-2.5) circle (0.015cm);      
\draw (6,-2.7) circle (0.015cm);
\fill  (6,-2.7) circle (0.015cm);  
\draw (6,-2.9) circle (0.015cm);
\fill  (6,-2.9) circle (0.015cm);

\draw (8,-2.5) circle (0.015cm); 
\fill  (8,-2.5) circle (0.015cm);      
\draw (8,-2.7) circle (0.015cm);
\fill  (8,-2.7) circle (0.015cm);  
\draw (8,-2.9) circle (0.015cm);
\fill  (8,-2.9) circle (0.015cm);

\draw (10,-2.5) circle (0.015cm); 
\fill  (10,-2.5) circle (0.015cm);      
\draw (10,-2.7) circle (0.015cm);
\fill  (10,-2.7) circle (0.015cm);  
\draw (10,-2.9) circle (0.015cm);
\fill  (10,-2.9) circle (0.015cm);

    ;
   \end{tikzpicture}
 \endpgfgraphicnamed
$
    
\end{sidewaysfigure}

\end{example}
\newpage

\noindent
\textbf{Acknowledgements:} This article is part of the author's Ph.D. thesis at
IMPA under the supervision of Oliver Lorscheid. He deeply
thanks him for his constant support, patience, encouragement and availability. He also thanks Olivier Schiffmann for hosting him for a term in Paris, Dragos Fratila for fruitful discussions and the reviewer for the carefully reading and helpful comments. Funding: This work was supported by  FAPERJ [grant 200.322/2016]; CAPES [grant  88881.134709/2016-01]; and FAPESP [grant number 2017/21259-3].


\bibliographystyle{plain}

\end{document}